\documentclass[reqno,12pt]{amsart}

\usepackage{amsmath}
\usepackage{amscd}
\usepackage{amssymb}
\usepackage{enumerate}
\usepackage{amsfonts}
\usepackage{graphicx}
\usepackage[all]{xy}
\usepackage{mathrsfs}

\usepackage[hypertex]{hyperref}

\usepackage{bbm}

\newcommand{\e}{\mathbbm{1}}

\numberwithin{equation}{section}

\newcommand{\rar}[1]{\stackrel{#1}{\longrightarrow}}
\newcommand{\xrar}[1]{\xrightarrow{#1}}

\newcommand{\isom}{\rar{\simeq}}

\newcommand{\into}{\hookrightarrow}

\newcommand{\al}{\alpha} \newcommand{\be}{\beta} \newcommand{\ga}{\gamma}
\newcommand{\Ga}{\Gamma} \newcommand{\de}{\delta} \newcommand{\De}{\Delta}
\newcommand{\la}{\lambda} \newcommand{\La}{\Lambda}
\newcommand{\ze}{\zeta} 
\newcommand{\eps}{\epsilon} \newcommand{\sg}{\sigma}
\newcommand{\Sg}{\Sigma} \newcommand{\te}{\theta} 
\newcommand{\om}{\omega} \newcommand{\Om}{\Omega}

\newcommand{\bA}{{\mathbb A}}

\newcommand{\bC}{{\mathbb C}}
\newcommand{\bD}{{\mathbb D}}

\newcommand{\bF}{{\mathbb F}}
\newcommand{\bG}{{\mathbb G}}

\newcommand{\bK}{{\mathbb K}}
\newcommand{\bL}{{\mathbb L}}

\newcommand{\bN}{{\mathbb N}}

\newcommand{\bQ}{{\mathbb Q}}
\newcommand{\bR}{{\mathbb R}}

\newcommand{\bZ}{{\mathbb Z}}

\newcommand{\cC}{{\mathcal C}}

\newcommand{\cE}{{\mathcal E}}
\newcommand{\cF}{{\mathcal F}}

\newcommand{\cH}{{\mathcal H}}

\newcommand{\cK}{{\mathcal K}}
\newcommand{\cL}{{\mathcal L}}
\newcommand{\cM}{{\mathcal M}}

\newcommand{\cO}{{\mathcal O}}

\newcommand{\cS}{{\mathcal S}}

\newcommand{\sA}{{\mathscr A}}

\newcommand{\sD}{{\mathscr D}}

\newcommand{\sG}{{\mathscr G}}
\newcommand{\sH}{{\mathscr H}}

\newcommand{\fF}{{\mathfrak F}}

\newcommand{\fa}{{\mathfrak a}}

\newcommand{\fc}{{\mathfrak c}}

\newcommand{\ff}{{\mathfrak f}}
\newcommand{\fg}{{\mathfrak g}}
\newcommand{\fh}{{\mathfrak h}}

\newcommand{\fm}{{\mathfrak m}}
\newcommand{\fn}{{\mathfrak n}}

\newcommand{\fz}{{\mathfrak z}}

\newcommand{\Gb}{\overline{G}}

\newcommand{\Gh}{\widehat{G}}

\newcommand{\abs}[1]{\vert #1\vert}

\newcommand{\limproj}{\lim\limits_{\longleftarrow}}
\newcommand{\limind}{\lim\limits_{\longrightarrow}}

\newcommand{\indlim}[1]{\lim\limits_{\stackrel{\longrightarrow}{#1}}}

\newcommand{\Ker}{\operatorname{Ker}}

\newcommand{\End}{\operatorname{End}}
\newcommand{\Hom}{\operatorname{Hom}}
\newcommand{\Ext}{\operatorname{Ext}}
\newcommand{\Aut}{\operatorname{Aut}}
\newcommand{\Spec}{\operatorname{Spec}}

\newcommand{\id}{\operatorname{id}}
\newcommand{\Id}{\operatorname{Id}}
\newcommand{\pr}{\mathrm{pr}}
\newcommand{\ad}{\operatorname{ad}}

\newcommand{\Ad}{\operatorname{Ad}}

\newcommand{\Span}{\operatorname{span}}
\newcommand{\Supp}{\operatorname{supp}}
\newcommand{\supp}{\operatorname{supp}}

\newcommand{\Ind}{\operatorname{Ind}}
\newcommand{\tr}{\operatorname{tr}}

\newcommand{\Exp}{\operatorname{Exp}}
\newcommand{\Log}{\operatorname{Log}}

\newcommand{\act}{\operatorname{act}}

\newcommand{\tens}{\otimes}

\newcommand{\st}{\,\big\vert\,}

\newcommand{\chr}{\operatorname{char}}

\newcommand{\sbr}{\smallbreak}
\newcommand{\mbr}{\medbreak}
\newcommand{\bbr}{\bigbreak}

\newcommand{\loccit}{{\em loc.\,cit.}}

\newcommand{\ti}{\widetilde}

\newcommand{\bra}{\langle}
\newcommand{\ket}{\rangle}

\newcommand{\eval}[2]{\bra #1\,\big\vert\,#2\ket}

\newcommand{\Z}{\mathsf{Z}}

\newcommand{\gh}{\widehat{\Ga}}

\newtheorem{conj}{Conjecture}

\newtheorem{thm}{Theorem}[section]
\newtheorem{cor}[thm]{Corollary}
\newtheorem{lem}[thm]{Lemma}

\newtheorem{prop}[thm]{Proposition}
\newtheorem{exer}[thm]{Exercise}

\theoremstyle{remark}
\newtheorem{rem}[thm]{Remark}
\newtheorem{rems}[thm]{Remarks}
\newtheorem{example}[thm]{Example}

\newtheorem{defin}[thm]{Definition}

\newcommand{\cst}{\bC^\times}

\newcommand{\Meas}{\operatorname{Meas}}
\newcommand{\Fun}{\operatorname{Fun}}
\newcommand{\Perm}{\operatorname{Perm}}

\newcommand{\Fr}{\operatorname{Fr}}

\newcommand{\Gal}{\operatorname{Gal}}

\newcommand{\Bt}{\widetilde{B}}

\newcommand{\Vt}{\widetilde{V}}

\newcommand{\mut}{\widetilde{\mu}}

\newcommand{\pit}{\widetilde{\pi}}

\newcommand{\Sch}{\operatorname{\cS{}ch}}

\newcommand{\Schp}{\operatorname{\cS{}ch}^{perf}}

\newcommand{\Nilp}{\operatorname{Nilp}}

\newcommand{\nilp}{\mathfrak{nilp}}

\newcommand{\card}{\operatorname{card}}

\newcommand{\ql}{\overline{\bQ}_\ell}

\renewcommand{\hom}{\sH\!om}
\newcommand{\hm}[1]{\operatorname{Hom}\bigl( #1 \bigr)}

\newcommand{\cs}{CS}

\newcommand{\Heis}{\operatorname{Heis}}

\newcommand{\pih}{\widehat{\Pi}}
\newcommand{\phh}{\widehat{\phi}}
\newcommand{\bzh}{\widehat{\bZ}}

\newcommand{\Perv}{\operatorname{Perv}}

\title[An introduction to character sheaves and the orbit method]{A
motivated introduction to character sheaves and the orbit method
for unipotent groups in positive characteristic}

\author[M.~Boyarchenko and V.~Drinfeld]{Mitya Boyarchenko and Vladimir Drinfeld}

\thanks{Both authors were supported by NSF grant
DMS-0401164. \\ \textit{Address}: Department of Mathematics, University of
Chicago, Chicago, IL 60637. \\ \textit{E-mail}: {\tt mitya@math.uchicago.edu}
(M.B.), \ {\tt drinfeld@math.uchicago.edu} (V.D.)}

\dedicatory{To Joseph Bernstein with admiration and warmest regards}

\begin{document}

\begin{abstract}
This article is based on lectures given by the authors in 2005 and
2006. Our first goal is to present an introduction to the orbit
method with an emphasis on the character theory of finite
nilpotent groups. The second goal (motivated by a recent work of
G.~Lusztig) is to explain several nontrivial aspects of character
theory for finite groups of the form $G(\bF_{q^n})$, where $G$ is
a unipotent algebraic group over $\bF_q$. In particular, we
introduce the notion of a character sheaf for a unipotent group,
and provide a toy model for the representation-theoretic notion of
an $\bL$-packet.
\end{abstract}

\maketitle

\setcounter{tocdepth}{1}

\tableofcontents


\section*{Introduction}

This paper is an expanded collection of notes for several lectures
on the orbit method and geometric representation theory given by
the authors at the University of Chicago in June and October of
2005, and more recently in June of 2006 at the Langlands Program
conference at CIRM, Luminy (cf.~\cite{charsheaves-orbmethod}).
These lectures originated in our attempt to understand \S7 of
George Lusztig's remarkable work \cite{lusztig}, in which he began
investigating character theory and the theory of character sheaves
for unipotent groups over finite fields and pointed out that its
features are similar in spirit to those of the theory of character
sheaves for general reductive groups. In particular, he provided
the first example of a nontrivial $\bL$-packet\footnote{Throughout this article, the term ``$\bL$-packet'' is used as an abbreviation of ``Lusztig packet.''} for representations
of unipotent groups.

\mbr

Let $G$ be a unipotent algebraic group over a finite field $\bF_q$. For each
positive integer $n$, we have a finite nilpotent group $G(\bF_{q^n})$, and it
is natural to ask whether there is a connection between the complex irreducible
characters of this group and some objects related to the algebraic group
$G\otimes_{\bF_q}\overline{\bF}_q$. We hope that the answer is positive, and we
show this for several classes of unipotent groups. We explain some nontrivial
aspects of this connection, and state several questions and conjectures in this
setup. Already at the classical level (i.e., before passing from functions to
sheaves), some of our questions, examples and observations are either new or at
least not part of the ``common knowledge''. For instance, we point out that
even if the nilpotence class of $G$ is less than the characteristic $p$ of the
ground field $\bF_q$ (in which case the orbit method applies), and even if $G$
has exponent $p$, there are important differences from the representation
theory of unipotent groups over a field of characteristic $0$. First of all,
the logarithm map identifies $G$, in general, not with its Lie algebra but with
a certain \textit{Lie ring scheme}\footnote{A {\em Lie ring} is an abelian
group $\fg$ with a biadditive map $[\cdot,\cdot]:\fg\times\fg\to\fg$ satisfying
the Jacobi identity and the identity $[x,x]$=0. A {\em Lie ring scheme} over a
field $k$ is a Lie ring object in the category of $k$-schemes.}, $\fg$, over
$\bF_q$. Second, the stabilizer in $G$ of a point of $\fg^*$ may be
disconnected. As a result, a geometric orbit of $G$ in $\fg^*$ may correspond
not to a single irreducible representation of $G(\bF_{q^n})$ but to a finite
collection of such representations (this is what we call an {\em $\bL$-packet}).
Finally, the $G$-orbits in $\fg^*$ may have odd dimension. It is precisely
these three phenomena that make the subject very interesting from our point of
view.

\mbr

However, the main novelty of our work is a general definition of character
sheaves and $\bL$-packets for unipotent groups in positive characteristic, and a
collection of conjectures about them. On the one hand, we show that our notion
of a character sheaf is compatible with the orbit method, in the sense that if
$G$ is a connected unipotent group over $\bF_q$ of nilpotence class less than
$p$, so that the Lie ring scheme $\fg$ of $G$ and its dual $\fg^*$ are defined,
then the character sheaves on $G\otimes_{\bF_q}\overline{\bF}_q$ are (up to
cohomological shift) the (inverse) Fourier transforms of irreducible
equivariant local systems on the coadjoint orbits for $G$, as one could expect.
On the other hand, the definition of a character sheaf is independent of the
nilpotence class of $G$, and we hope that character sheaves are closely related
to irreducible characters of $G(\bF_{q^n})$ even in the cases when the orbit
method cannot be applied (notably, when $G$ is the group $UL_{N,q}$ of
unipotent upper-triangular matrices of size $N$ over $\bF_q$ where $N\gg p$).
We also expect the theory of character sheaves to have applications to
geometric representation theory beyond the study of irreducible characters, but
these applications lie outside the scope of this article.

\mbr

We aimed to make this paper accessible to as large an audience as
possible. It is essentially self-contained: we carefully state all
the results we use and give all the necessary definitions, even
though we omit a few proofs. For the most part the reader will
only need a superficial knowledge of algebraic groups over
possibly non-algebraically closed fields, as well as a few basic
facts about representations of finite groups \cite{serre-reps}.
Some familiarity with the classical orbit method \cite{kirillov}
would help the reader, but it is not required.

\mbr

Another feature of our presentation is a large collection of
examples that help motivate the main constructions and results of
the paper. Among them we would like to mention the ``fake
Heisenberg groups'', unipotent algebra groups, and maximal
unipotent subgroups of symplectic groups. Most of the proofs
appearing in this work are straightforward.

\subsection*{Organization} Given the size of the paper, we would
like to assure the reader that it does not have to be read in a
linear fashion. For instance, the reader looking for a basic
introduction to the orbit method only needs to read
\S\ref{s:orbmethod}, parts of Appendices \ref{a:reduction},
\ref{a:vergne}, and Appendix \ref{a:counterexamples}, while
referring to Appendices \ref{a:characters} and \ref{a:heisenberg}
if necessary (see also \cite{orbmethod}). After that, for a deeper
understanding of the geometric aspects of character theory in our
setting, we recommend reading \S\S\ref{s:overview},
\ref{s:orbmethod-unipotent} and studying the examples presented in
\S\S\ref{s:algebra-groups}, \ref{s:usp4}. Finally, the theory of
character sheaves for unipotent groups is explained in
\S\ref{s:char-sheaves-orbmethod} and
\S\ref{s:charsheaves-general}. It is also possible to read
\S\S\ref{s:overview}--\ref{s:charsheaves-general} as a
self-contained introduction to this theory (some, but not all, of
this material is presented in \cite{charsheaves-orbmethod}).

\mbr

In more detail, the organization of the paper is as follows. In
\S\ref{s:overview} we describe the setup of the article, formulate
a few questions that arise in the geometric setting, and answer
them for connected commutative algebraic groups. These questions
motivate almost all the rest of the paper. In \S\ref{s:orbmethod}
we review the orbit method for finite nilpotent groups. In
\S\ref{s:orbmethod-unipotent} we explain some nontrivial aspects
of this method for unipotent groups over finite fields, discuss
the $\bL$-indistinguishability phenomenon, and illustrate our
discussion with the example of the ``fake Heisenberg groups''. We
then move on to an introduction to character sheaves in the
unipotent setting. In \S\ref{s:char-sheaves-orbmethod} we give an
\emph{ad hoc} definition of character sheaves for a connected
unipotent group of ``small'' nilpotence class and explain their
relationship with irreducible characters. Afterwards, in
\S\ref{s:charsheaves-general}, we give a general definition of
character sheaves for unipotent groups and show that it agrees
with the \emph{ad hoc} one when the latter is applicable. The main
body of the paper concludes with \S\ref{s:algebra-groups} and
\S\ref{s:usp4} where we introduce several classes of examples that
are interesting both from the geometric and from the
group-theoretic points of view.

\mbr

The appendices are devoted mostly to background material that
could not be explained in our lectures due to the time
constraints. The first four appendices are purely algebraic, the
next three deal with certain geometric aspects, and the last one
discusses some counterexamples. In more detail, Appendix
\ref{a:characters} proves some auxiliary results about irreducible
characters of finite groups, and Appendix \ref{a:heisenberg} gives
a self-contained introduction to the theory of ``Heisenberg
representations'', which is an important special case of the orbit
method. Appendices \ref{a:reduction} and \ref{a:vergne} contain a
proof of the fact that every irreducible representation of a
finite nilpotent group is induced from a $1$-dimensional
representation of a subgroup, as well as a detailed discussion of
two different constructions of polarizations. In particular, we
show that the construction due to M.~Vergne is very natural from
the viewpoint of representation theory. Appendix
\ref{a:equivariant-derived} is a brief summary of what one needs
to know about derived categories of constructible $\ell$-adic
complexes, and equivariant versions of these categories, in order
to read those parts of the paper that are devoted to character
sheaves. In Appendices \ref{a:duality} and \ref{a:fourier-deligne}
we review Serre duality and the Fourier-Deligne transform for
perfect connected commutative unipotent groups over perfect
fields; we also recall the definition of perverse sheaves.
Finally, the goal of Appendix \ref{a:counterexamples} is to show
that certain more-or-less natural conjectures that arise in the
context of the orbit method for finite nilpotent groups are
actually {\em false}; the reader familiar with the classical orbit
method may wish to consult it right away in order to avoid some
possible pitfalls.

\subsection*{Remark} We would like to emphasize that in this paper
we only consider representations of groups over algebraically
closed fields of \emph{characteristic zero}, such as $\bC$ or
$\ql$, whereas the ground field for all our geometric objects must
have \emph{positive characteristic}, in order for our discussion
to be interesting and/or meaningful.

\subsection*{Acknowledgements} We would like to thank David Kazhdan for his
suggestion to publish the notes for our lectures, and for many comments on our article which have significantly improved the quality of the presentation. We are grateful to George
Lusztig for drawing our attention to algebra groups. We are also grateful to
David Kazhdan and Robert Kottwitz for several useful discussions. We are
greatly indebted to Maria Sabitova for catching many errors and misprints in
our text. Last but not least, we thank Masoud Kamgarpour for taking great notes
for the first author's lectures, for helpful comments on the earlier versions
of this paper, and for useful conversations about the various notions of a
Frobenius morphism.

\subsection*{Warning}
We use two different conventions in different parts of the text.
Namely, when we study irreducible characters and consider the
elementary questions of \S\ref{ss:main-geom-questions}, we always
work either with a finite group usually denoted by $\Ga$ (and
sometimes by $\Pi$), or with an algebraic group $G$ over a finite
field $\bF_q$. On the other hand, it is more convenient to
formulate the theory of character sheaves for algebraic groups
over an arbitrary algebraically closed field $k$. In this case, we
denote such a group by $G$, and if $\operatorname{char}k=p$ and
$G$ has an $\bF_q$-structure, we will denote by $G_0$ an algebraic
group over $\bF_q$ such that $G_0\otimes_{\bF_q}k\cong G$. More
generally, in the sections devoted to character sheaves we follow
the standard practice (see, e.g., \cite{bbd}) where the notation
for objects over $\bF_q$ contains a subscript ``$0$'', and
omission of this subscript indicates extension of scalars to $k$.


\section{Elementary questions and examples}\label{s:overview}

\subsection{Algebraic aspects}\label{ss:alg-aspects}
If $\Ga$ is a finite group, we will write $\gh$ for the set of isomorphism
classes of irreducible representations of $\Ga$ over some fixed algebraically
closed field $\overline{E}$ of characteristic $0$. The ``minimal'' problem in
character theory that one would like to solve is to describe the set $\gh$
together with the ``dimension function''
\[
\dim:\gh\rar{}\bN.
\]
Sometimes one would like to know not only the function $\dim$, but also the map
\[
\chi:\gh\rar{}\Fun(\Ga),
\]
where $\Fun(\Ga)$ denotes the space of all $\overline{E}$-valued functions on
$\Ga$, and $\chi$ takes an irreducible representation $\rho$ of $\Ga$ to its
character $\chi_\rho(g)=\tr(g;\rho)$ ($g\in \Ga$). Of course, the function
$\chi$ contains more information than $\dim$, since $\dim\rho=\chi_\rho(1)$,
where $1$ denotes the identity element of $\Ga$. When $\Ga$ is a finite
nilpotent group of ``not too large'' nilpotence class\footnote{We recall that
if $G$ is a (discrete or algebraic or Lie, etc.) group, the {\em nilpotence
class} of $G$ is defined as the smallest integer $n$ such that $G^{(n)}=\{1\}$
(if it exists), where $G^{(n)}$ is the $n$-th iterated commutator of $G$,
defined inductively by $G^{(0)}=G$, $G^{(i+1)}=[G,G^{(i)}]$. A similar
definition applies to Lie rings.} $c$ (the precise condition is that $c!$ is
prime to $\abs{\Ga}$, the order of $\Ga$), Kirillov's orbit method gives a
rather satisfactory solution of the last problem. We recall it in Section
\ref{s:orbmethod}.

\mbr

In this article we would like to study the problems mentioned
above for finite groups of the following special type. Let us fix
a prime $p$ and an algebraic closure $\bF$ of $\bF_p$. For each
$k\in\bN$, let $\bF_{p^k}$ denote the unique subfield of $\bF$ of
order $p^k$. Let $G$ be an algebraic group (i.e., a reduced group
scheme of finite type) over a finite subfield $\bF_q\subset\bF$,
and define $\Ga_n=G(\bF_{q^n})$ for each $n\in\bN$. The finite
groups $\Ga_n$ form an inductive system: if $m\lvert n$, then
$\bF_{q^m}\subseteq\bF_{q^n}$, and so $\Ga_m\subseteq\Ga_n$. One
would like to study all the sets $\gh_n$ simultaneously, to
establish connections between them for various $n$, and to relate
them to the algebraic group $G$. If $G$ is reductive, the theory
of character sheaves developed by Lusztig (cf.
\cite{char-sheaves}, and also \cite{lusztig} and the references
therein) solves this problem.

\subsection{A note on Frobenius morphisms}\label{ss:frobenii} In what
follows we will need to work with Frobenius morphisms for schemes of
characteristic $p>0$. Since there are several variants of these, we will use
this subsection to fix our definitions and notation, in order to avoid all
possible confusion. If $S$ is a scheme over $\bF_q$, we denote by
$\Phi_{S,q}:S\to S$ the $\bF_q$-morphism defined as the identity on the
underlying topological space, and the map $f\mapsto f^q$ on local sections of
the structure sheaf $\cO_S$. If the scheme $S$ and/or the prime power $q$ are
fixed throughout a particular discussion, we may omit them from the notation.
We say that $S$ is {\em perfect} if $\Phi_{S,q}$ is an automorphism of $S$.
This notion depends only on the characteristic of $S$ and not on $q$, since if
$q=p^s$ for some prime $p$, then $S$ can be viewed as a scheme over $\bF_p$,
and $\Phi_{S,q}=\Phi_{S,p}^s$, whence $\Phi_{S,q}$ is invertible if and only if
$\Phi_{S,p}$ is.

\mbr

If $S,T$ are two schemes over $\bF_q$, the two endomorphisms of
$\operatorname{Mor}_{\bF_q}(S,T)$ induced by $\Phi_{S,q}$ and $\Phi_{T,q}$
clearly coincide\footnote{In other words, if $f:S\to T$ is an $\bF_q$-morphism,
then $\Phi_{T,q}\circ f=f\circ\Phi_{S,q}$.}. Thus we obtain an endomorphism of
$\operatorname{Mor}_{\bF_q}(S,T)$ which we denote by $\Fr_q$. Note that $\Fr_q$
is an automorphism provided either of the two schemes $S$ and $T$ is perfect.
In particular, we obtain an automorphism $\Fr_q$ of the set
$S(\bF)=\operatorname{Mor}_{\bF_q}(\Spec\bF,S)$ of geometric points of $S$. For
any $n\in\bN$, it restricts to an automorphism of the set $S(\bF_{q^n})$, which
we again denote by $\Fr_q$. If $S$ is an algebraic group over $\bF_q$, then
$\Fr_q$ is a group automorphism. As a rule, it is this automorphism that we
refer to below as {\em the Frobenius}.

\mbr

On a few occasions (see \S\ref{ss:perfect-duality} and \S\ref{aa:perfect}) we
will need to use the relative Frobenius morphism. If $X$ is a scheme over an
arbitrary field $k$ of characteristic $p$, we will write $X^{(p)}$ for the
scheme over $k$ obtained as the fiber product of the structure morphism
$X\to\Spec k$ and the morphism $\Phi_{k,p}:\Spec k\to \Spec k$. By the
universal property of the fiber product, the morphism $\Phi_{X,p}:X\to X$ and
the structure morphism $X\to\Spec k$ induce a morphism $\Phi_{X/k}:X\to
X^{(p)}$ of schemes over $\Spec k$; it is called the {\em relative Frobenius
morphism}.

\subsection{Elementary questions}\label{ss:main-geom-questions} In the situation of
\S\ref{ss:alg-aspects}, we would like to pose the following questions, even
though they are imprecisely stated and probably too naive.

\begin{enumerate}[(1)]
\item Do there exist natural ``base change maps''
\begin{equation}\label{e:map}
T_m^n : \gh_m \rar{} \bigl(\gh_n\bigr)^{\Gal(\bF_{q^n}/\bF_{q^m})}
\end{equation}
for all pairs of positive integers $m|n$? Note that we cannot
state this question precisely because at this point we do not know
what ``natural'' means. However, naturality should at least
include the equivariance of $T_m^n$ with respect to
$\widehat{\Fr_q}$, and the compatibility condition
\[
T_m^k = T_n^k \circ T_m^n \qquad \text{whenever} \quad m|n|k.
\]
We expect that the answer to this question may be negative in
general, but positive for many interesting classes of examples.
For instance, if $G$ is connected and unipotent of nilpotence
class $<p$, Kirillov's orbit method provides a positive answer:
see \S\ref{ss:orbmet-unip}. On the other hand, if $G$ is a
unipotent algebra group as defined in Section
\ref{s:algebra-groups}, base change maps for $G$ were constructed
in \cite{base-change}; this case is substantially different. \sbr
\item Assuming that the answer to question (1) is positive, form
the direct limit\footnote{We remark that $\gh$ is {\em not} the set of
irreducible characters of any group, but should be thought of as just one
symbol, defined by \eqref{e:gamma-hat}.}
\begin{equation}\label{e:gamma-hat}
\gh:=\indlim{n} \gh_n
\end{equation}
with respect to the base change maps, and consider the induced
$\widehat{\Fr_q}$-equivariant maps
\begin{equation}\label{e:maps}
T_n : \gh_n\rar{} \bigl(\gh\bigr)^{\Gal(\bF/\bF_{q^n})}.
\end{equation}
Are the base change maps $T_n$ injective (resp., surjective)? We know that the
answer to this question is negative in general; however, the investigation of
the conditions under which the answer is positive leads to interesting
geometric questions. For example, $T_n$ is always surjective when $G$ is
connected and unipotent of nilpotence class $<p$ (see \S\ref{ss:l-indistinguishability}),
and $T_n$ is always injective when $G$ is a unipotent algebra group (see
\cite{base-change}). \sbr
\item Assuming again that the answer to question (1) is positive, does
there exist a ``geometric object'' $\Gh$ defined over $\bF_q$ together with a
$\widehat{\Fr_q}$-equivariant bijection $\Gh(\bF)\rar{\simeq}\gh$? This
question is also not well posed, since we do not define the class of
``geometric objects'' that we will consider. Nevertheless, we expect that the
answer is positive for many interesting classes of examples of {\em unipotent}
groups over $\bF_q$. For example, this is the case when $G$ is unipotent of
nilpotence class $<p$. For unipotent algebra groups the answer is unknown to us
at present. (On the other hand, the answer is {\em negative} for all groups
that are not unipotent; see \S\ref{ss:basic-examples} for the case $G=\bG_m$.)
\end{enumerate}

\begin{defin}\label{d:l-indistinguishability}
Assume that the answer to question (1) is positive. The fibers of the maps
\eqref{e:maps} will be called {\em $\bL$-packets}, and two irreducible
representations of $\Ga_n$ that have the same image in $\gh$ will be said to be
{\em $\bL$-indistinguishable}. The $\bL$-indistinguishability phenomenon is
discussed in more detail in \S\ref{ss:l-indistinguishability} and
\S\ref{ss:fake-heisenberg}.
\end{defin}

\begin{rem}\label{r:langlands}
In the representation theory of reductive groups over \emph{local} fields there is a conjectural notion of an $L$-packet introduced by R.P.~Langlands \cite{langlands}. It is hard to compare it with our notion of $\bL$-packet because technically the two definitions are given in
quite different terms. But philosophically the two notions are closely related. Namely,
as explained to us by R.~Bezrukavnikov, $\bL$-packets are philosophically similar to {\em geometric} $L$-packets, which are, in general, bigger than the $L$-packets defined by Langlands.\footnote{Conjecturally, $L$-packets bijectively correspond to ``Langlands parameters".  Geometric $L$-packets should correspond to Frobenius-invariant ``geometric Langlands parameters" (one gets the geometric Langlands parameters from the usual ones by replacing the Weil-Deligne group $W'_K$ with $\Ker (W'_K\twoheadrightarrow\bZ)\,$). Thus each geometric $L$-packet is a union of several usual $L$-packets.}
\end{rem}

\subsection{Character sheaves}\label{ss:charsheaves-introduction}
We still consider an algebraic group $G$ over a finite field
$\bF_q$. We also fix a prime $\ell\neq p$ and take
$\overline{E}=\ql$ as the field of coefficients for our
representations. According to the conventions of
\S\ref{ss:frobenii}, we have the absolute Frobenius morphism
$\Phi_q:G\rar{}G$, and we may extend scalars and form the morphism
$\Fr_q:=\Phi_q\tens\id:G\tens_{\bF_q}\bF\rar{}G\tens_{\bF_q}\bF$
of algebraic groups over $\bF$ (not to be confused with the
absolute Frobenius endomorphism of $G\tens_{\bF_q}\bF$). Note that
the symbol $\Fr_q$ has also appeared in \S\ref{ss:frobenii}, but
the notation introduced there is consistent with the one we use
here. More precisely, there is a natural identification of
$G(\bF)$ with the set of closed points of $G\tens_{\bF_q}\bF$, and
under this identification the action of $\Fr_q$ introduced in
\S\ref{ss:frobenii} corresponds to the action induced by
$\Phi_q\tens\id$.

\mbr

From the point of view of representation theory of finite groups,
the study of character sheaves is motivated by the following
question, which is at least as important as the questions posed in
\S\ref{ss:main-geom-questions}, and should be viewed as a part of
that list.

\mbr

(4) Does there exist a collection $\cs(G)$ of $\ell$-adic
complexes $\cF\in D^b_c(G\tens_{\bF_q}\bF,\ql)$ (the notation is
explained in \S\ref{aa:derived-constructible}) which enjoy
properties (4-i)--(4-iii) below?

\begin{enumerate}[$($4-i$)$]
\item The construction of $\cs(G)$ only depends
on the algebraic group $\overline{G}:=G\tens_{\bF_q}\bF$ and not
on $G$, i.e., the set of isomorphism classes of the complexes in
$\cs(G)$ is invariant under all automorphisms of $\overline{G}$ as
an algebraic group over $\bF$.
\item The complexes in $\cs(G)$ are irreducible perverse sheaves
on $\overline{G}$, cf. \S\ref{aa:perverse}.
\item For each $n\in\bN$, consider the
subset $\cs(G)^{\Fr_q^n}\subseteq\cs(G)$ consisting of those
complexes $\cF\in\cs(G)$ such that $(\Fr_q^n)^*(\cF)\cong\cF$.
Then it is possible to choose an isomorphism
$\psi_{n,\cF}:(\Fr_q^n)^*(\cF)\rar{\simeq}\cF$ for every
$\cF\in\cs(G)^{\Fr_q^n}$, such that $\psi_{k,\cF}$ is induced by
$\psi_{n,\cF}$ whenever $n\lvert k$ in the sense that
\[
\psi_{k,\cL} = \psi_{n,\cL}\circ (\Fr_q^n)^*(\psi_{n,\cL}) \circ
(\Fr_q^{2n})^*(\psi_{n,\cL}) \circ \dotsb \circ
(\Fr_q^{k-n})^*(\psi_{n,\cL}),
\]
and such that the trace functions
\[
t_{n,\cF} : G(\bF_{q^n}) \rar{} \ql
\]
defined by
\[
t_{n,\cF}(g) = \sum_{i\in\bZ} (-1)^i \tr\bigl( \psi_{n,\cF};
\cH^i(\cF)_g \bigr)
\]
are precisely the irreducible characters of $G(\bF_{q^n})$ over
$\ql$.
\end{enumerate}

\begin{rems}\label{r:charsheaves}
\begin{enumerate}[(a)]
\item In the formula above, $\cH^i(\cF)$ denotes the $i$-th
cohomology sheaf of the complex $\cF$, and $\cH^i(\cF)_g$ is its
stalk at the point $g\in G(\bF_{q^n})$, which is viewed as a
$\Fr_q^n$-stable closed point of the group $\overline{G}$. Thus
$\psi_{n,\cF}$ does indeed act on the vector space $\cH^i(\cF)_g$,
so the formula makes sense.
\item The passage from the pair consisting of the complex $\cF$ on $\overline{G}$
and the isomorphism $\psi_{n,\cF}:(\Fr_q^n)^*(\cF)\rar{\simeq}\cF$
to the function $t_{n,\cF}$ on $G(\bF_{q^n})$ is known as the
\emph{functions-sheaves correspondence} (see, for instance,
\cite{sga4.5}, \emph{Sommes trig.}).
\item In view of the requirement (4-ii), the isomorphisms
$\psi_{n,\cF}$ are unique up to scaling.
\item A positive answer to question (4) would yield a positive
answer to question (1). Indeed, if there exists a collection
$\cs(G)$ satisfying the properties listed above, then, in
particular, property (4-iii) implies that there are natural
bijections between $\widehat{G(\bF_{q^n})}$ and $\cs(G)^{\Fr_q^n}$
for all $n\in\bN$, and one can use these bijections, together with
the obvious inclusions $\cs(G)^{\Fr_q^m}\subseteq\cs(G)^{\Fr_q^n}$
for $m\lvert n$, to define the base change maps.
\end{enumerate}
\end{rems}

If the answer to question (4) is positive, the elements of the collection
$\cs(G)$ will be called the \emph{character sheaves} on $\overline{G}$ (which
explains the notation). For example, when $G=GL_n$, Lusztig proved
\cite{char-sheaves} that the answer to question (4) is indeed positive.
However, Lusztig showed that \emph{the answer is negative for many other
connected reductive groups}. More recently, Lusztig also observed
\cite{lusztig} that \emph{the answer is negative for some classes of connected
unipotent groups} as well. This observation is explained in
\S\ref{ss:dim-irreps}.

\mbr

In view of these comments, one may wish to relax condition (4-iii)
as follows: instead of requiring that the functions $t_{n,\cF}$
are precisely the irreducible characters of $G(\bF_{q^n})$, one
can merely ask these functions to form a basis for the space of
class functions $G(\bF_{q^n})\rar{}\ql$. With this formulation,
question (4) has a chance of having a positive answer for
\emph{connected} algebraic groups $G$. Indeed, Lusztig proved
\cite{char-sheaves} that this is so for connected reductive $G$,
and we conjecture that this is also true for connected
\emph{unipotent} $G$; see Section \ref{s:charsheaves-general}
where the definition of character sheaves for such groups is
given, and where the conjecture is stated in a precise way.
Moreover, we prove (see \S\ref{ss:charsheaves-orbmethod}) that
this conjecture holds whenever the nilpotence class of $G$ is
$<p$. In this section (\S\ref{ss:charsheaves-commutative}) we will
show that (the strong form of) question (4) has a positive answer
for connected \emph{commutative} algebraic groups.

\subsection{Characters of abelian groups}\label{ss:characters-abelian-groups}
Until the end of the section we study commutative algebraic
groups. Let us fix a power $q$ of a prime number $p$. Thus we will
simply write $\Fr$ in place of $\Fr_q$. Also, by abuse of
notation, we will write $\Fr$ in place of $\widehat{\Fr}$ whenever
this cannot cause confusion (an exception is Remark
\ref{r:frobenius}). We begin by recalling a famous theorem of
Serge Lang \cite{lang}: if $G$ is an algebraic group over $\bF_q$,
and $\Ga=G(\bF)$ with the natural action of the absolute Galois
group $\sG=\Gal(\bF/\bF_q)\cong\widehat{\bZ}$, then
\[
H^1(\sG,\Ga) = \begin{cases} \text{trivial if } G \text{ is
connected}, \\
\text{finite in general}.
\end{cases}
\]

\mbr

Let $G$ be a commutative algebraic group over $\bF_q$, and put
$\Ga_n=G(\bF_{q^n})$, as before. We will write the group operation in $G$ {\em
additively}, which differs from Lusztig's convention \cite{lusztig}, since he
writes it multiplicatively. Accordingly, what Lusztig calls the norm maps will
be called the trace maps here:
\[
\tr_m^n : \Ga_n \rar{} \Ga_m \quad \text{for  } m|n,
\]
\[
\tr_m^n(\ga)=\ga+\Fr^m(\ga)+\Fr^{2m}(\ga)+\dotsb+\Fr^{n-m}(\ga).
\]
Observe that $\tr^n_m$ is the restriction of the group homomorphism
\[
\tau_m^n  = \id + \Fr^m + \dotsb + \Fr^{n-m} : \Ga \rar{} \Ga.
\]

\mbr

Until further notice, we assume that $\overline{E}=\bC$, i.e., we will work
with complex representations. The map $\tr_m^n$ is a group homomorphism, so via
pullback it induces a map
\[
\widehat{\tr}_m^n : \gh_m=\Hom(\Ga_m,\cst) \rar{} \gh_n=\Hom(\Ga_n,\cst).
\]
Note that the maps $\widehat{\tr}_m^n$ are natural in the sense explained in
question (1) of \S\ref{ss:main-geom-questions}.
\begin{prop}[cf.~\cite{lusztig}, \S4]\label{p:base-change-abelian}
If $G$ is connected\footnote{Recall that an algebraic group $G$ over a field
$F$ is connected if and only if $G\tens_F \overline{F}$ is connected. Indeed,
we may assume that $F$ is perfect, and the result follows by noting that the
neutral component of $G\tens_F \overline{F}$ is fixed by
$\Gal(\overline{F}/F)$, and hence is defined over $F$. See also
\cite{waterhouse}, Theorem 6.6.}, this map gives an isomorphism
\[
\gh_m\rar{} \bigl(\gh_n\bigr)^{\Fr^m}.
\]
\end{prop}

\begin{rem}\label{r:disconnected}
The proposition may fail if $G$ is not connected: for example, consider the
case where $G$ is a discrete finite abelian group over $\bF_q$.
\end{rem}

\begin{proof}
Let us write $L_n:\Ga\to\Ga$ for the $n$-th Lang isogeny,
$L_n(\ga)=\Fr^n(\ga)-\ga$. We have to prove the exactness of the sequence
\[
0 \rar{} \gh_m \rar{\widehat{\tr}^n_m} \gh_n \rar{\widehat{L}_m} \gh_n,
\]
which is equivalent to the exactness of
\begin{equation}\label{e:Dag}
\Ga_n \rar{L_m} \Ga_n \rar{\tr_m^n} \Ga_m \rar{} 0.
\end{equation}
Since $G$ is connected, Lang's theorem implies that $L_n:\Ga\rar{}\Ga$ is
surjective. On the other hand, we have $L_n=\tau^n_m\circ L_m$, so the map
$\tau_m^n$ is surjective. Now given $\ga\in\Ga_m$, let $x\in\Ga$ be such that
$\ga=\tau_m^n(x)$. Then $x\in\Ga_n$, because
$L_n(x)=L_m(\tau_m^n(x))=L_m(\ga)=0$. Thus $\tr_m^n$ is also surjective, which
means that \eqref{e:Dag} is exact on the right.

\mbr

To prove its exactness in the middle, let $\ga\in\Ga_n$ be such that
$\tr_m^n(\ga)=0$; in any case, we have $\ga=L_m(x)$ for some $x\in\Ga$, and
then $L_n(x)=\tr_m^n(\ga)=0$ shows that $x\in\Ga_n$.
\end{proof}

\mbr

In the situation of the proposition, let us consider $\gh$ as defined above:
$\gh=\indlim{n}\gh_n$, the transition homomorphisms being the
$\widehat{\tr}_m^n$. We see that the natural map
\[
\gh_n \rar{} \bigl(\gh\bigr)^{\Gal(\bF/\bF_{q^n})}
\]
is an isomorphism for all $n$, when $G$ is connected. So in this
case the abelian group $\gh$, together with the action of $\sG$,
``captures'' the representation theory of all the finite groups
$\Ga_n$ in a compatible way. Thus we have solved the desired
problem at the algebraic level, i.e., we have (affirmatively)
answered questions (1) and (2) of \S\ref{ss:main-geom-questions}.
(The preceding discussion corresponds more or less to \S4 and part
of \S5 of \cite{lusztig}.)

\subsection{Basic examples}\label{ss:basic-examples}
From the geometric point of view, one should ask whether $\gh$ can be naturally
identified with the set of $\bF$-points of some commutative group scheme $\Gh$
defined over $\bF_q$. (Again, ``naturally'' means in particular that the
isomorphism $\Gh(\bF)\cong\gh$ should commute with the action of the
Frobenius.) The answer is not always positive.

\begin{example}\label{ex:additive}
For $G=\bG_a$, we claim that $\gh$ can be naturally identified with
$\bG_a(\bF)$, so that $\widehat{\bG}_a\cong\bG_a$. Namely, let us fix a
nontrivial additive character $\psi:\bF_q\to\cst$. For each $n\in\bN$, we can
use the {\em ring scheme} structure on $\bG_a$ to define natural maps
\[
\bG_a(\bF_{q^n}) \xrar{\ \ \simeq\ \ } \widehat{\bG_a(\bF_{q^n})},
\]
\[
x\mapsto \psi\bigl(\tr_{\bF_{q^n}/\bF_q}(x\cdot ?)\bigr).
\]
These maps are clearly $\Fr$-equivariant and compatible with each other for
different $n$ in the obvious sense. So, in this case, $\Gh$ exists as a
commutative algebraic group.
\end{example}

\begin{rem}\label{r:frobenius}
Let $\widehat{\Fr}:\gh\to\gh$ be the group homomorphism corresponding to
$\Fr:\Ga\to\Ga$. Using the above identification of $\gh$ with $\bG_a(\bF)$, we
can view $\widehat{\Fr}$ as an endomorphism of $\bG_a(\bF)$. Then
\begin{equation}\label{e:(star)}
\widehat{\Fr}(u)=u^{1/q}, \qquad u\in \bG_a(\bF)=\bF.
\end{equation}
Indeed, checking this for $u\in\bF_{q^n}$ amounts to the orthogonality of
$\Fr:\bF_{q^n}\to\bF_{q^n}$ with respect to the scalar product
$\bF_{q^n}\times\bF_{q^n}\to\cst$ defined by
$(x,y)\mapsto\psi\bigl(\tr_{\bF_{q^n}/\bF_q}(xy)\bigr)$. In particular,
$\widehat{\Fr}:\bG_a(\bF)\to\bG_a(\bF)$ {\em is not induced by a regular map of
algebraic varieties}.
\end{rem}

\begin{example}\label{ex:multiplicative}
For $G=\bG_m$, the answer to our question is negative. For if $\Gh$ exists,
then it is easy to show that for a prime $l$ different from
$\operatorname{char}\bF_q$, the Tate module $T_l(\widehat{\bG}_m)$ should be
dual to the Tate module $T_l(\bG_m)$. So $\widehat{\Fr}:T_l(\widehat{\bG}_m)\to
T_l(\widehat{\bG}_m)$ should have $q^{-1}$ as an eigenvalue, which is
impossible (the eigenvalues of the Frobenius acting on the Tate module of a
commutative algebraic group are known to be algebraic integers).
\end{example}

\subsection{Commutative unipotent groups}\label{ss:perfect-duality}
Let $G$ be an arbitrary connected commutative algebraic group over $\bF_q$.
From now on, to avoid confusion, we switch notation and replace $\Gh$ with
$G^*$. (Thus, the notation $G^*$ will always be used in the commutative case,
while $\Gh$ will be used in the general case.) The argument used in Example
\ref{ex:multiplicative} shows that $G^*$ cannot exist as a commutative
algebraic group over $\bF_q$ unless $G$ is unipotent. On the other hand, when
$G$ is unipotent, it turns out that $G^*$ always exists. It is known under the
name ``Serre dual'' of $G$, and the idea of the of duality theory for
commutative unipotent groups goes back to Serre's article \cite{grps-proalg}.
However, in the form needed for our purposes, the duality appears to be due to
L.~Begueri (see \cite{begueri} and \cite{saibi}).

\mbr

To be more precise, let us fix an arbitrary field $k$ of characteristic $p$.
The duality functor $G\mapsto G^*$ is defined as a functor
$\cC_k^{opp}\to\cC_k$, where $\cC_k$ is a certain localization of the category
$\widetilde{\cC}_k$ of connected commutative unipotent algebraic groups over
$k$. (Remark \ref{r:frobenius} shows that one cannot expect to have a duality
functor $\widetilde{\cC}_k\to\widetilde{\cC}_k$.) The definition of the
localization $\cC_k$ is as follows: one inverts all
$\widetilde{\cC}_k$-morphisms $f:G_1\to G_2$ such that the corresponding map
$G_1(\overline{k})\to G_2(\overline{k})$ is bijective. (Of course, the same
category $\cC_k$ is obtained if one only inverts the relative Frobenius
morphism $G\to G^{(p)}$ for each $G\in\widetilde{\cC}_k$.) One has $G^{**}=G$.

\mbr

A brief outline of this duality theory is given in Appendix
\ref{a:duality}. The reader may prefer to skip it, since for our
purposes it is interesting enough to consider only those
$G\in\cC_k$ which are isomorphic\footnote{It is well known and
easy to prove that if $G\in\cC_k$ and $p\cdot G=0$, then
$G\cong\bG_a^n$ for some $n$. If $k$ is perfect, this is already
true in $\widetilde{\cC}_k$.} to $\bG_a^n$, and on this full
subcategory the duality functor has the following explicit
description, which can be used as an {\em ad hoc} definition. (Of
course, this description agrees with Example \ref{ex:additive} and
Remark \ref{r:frobenius}.)

\mbr

Since the duality functor is additive, it suffices to know how it acts on
$\bG_a$. One has $(\bG_a)^*=\bG_a$, and it remains to describe the
anti-automorphism $f\mapsto f^*$ of the ring $\End_{\cC_k}(\bG_a)$. In fact, it
is enough to describe $f^*$ in the cases $f(x)=x^p$, $f(x)=x^{p^{-1}}$,  and
$f(x)=cx$ (this is a system of generators of $\End_{\cC_k}(\bG_a)$). In the
first two cases, $f^*=f^{-1}$, and in the third case, $f^*=f$.

\subsection{Character sheaves for connected commutative
groups}\label{ss:charsheaves-commutative} Let $G$ be a connected
commutative algebraic group over $\bF_q$, and let
$\overline{G}=G\tens_{\bF_q}\bF$, as before. Character sheaves for
$\overline{G}$ were constructed by Lusztig in \S5 of
\cite{lusztig}. We will explain this result from a somewhat
different point of view, which is closer in spirit to Section
\ref{s:charsheaves-general} below.

\mbr

Let $\al:G\times G\to G$ denote the group operation. By abuse of
notation, the group operation on $\Gb$ will also be denoted by
$\al$. We define a \emph{character sheaf} on $\Gb$ to be a local
system (i.e., a lisse $\ql$-sheaf) $\cL$ of rank $1$ on $\Gb$
satisfying $\al^*(\cL)\cong\cL\boxtimes\cL$. The collection
$\cs(G)$ of character sheaves on $\Gb$ clearly satisfies property
(4-i) of \S\ref{ss:charsheaves-introduction}. It satisfies
property (4-ii) up to cohomological shift, in the sense that if
$\cL\in\cs(G)$, then $\cL[\dim G]$ is an irreducible perverse
sheaf on $G$ (here we use the assumption that $G$ is connected).
In the present setting we prefer to work with local systems rather
than perverse sheaves in order to avoid having to keep track of an
extra factor of $(-1)^{\dim G}$ in our computations. This issue is
absolutely unimportant.

\mbr

Let us check that $\cs(G)$ satisfies property (4-iii). First of
all, for each $n\in\bN$ and each $\cL\in\cs(G)^{\Fr^n}$, we choose
the unique isomorphism $\psi_{n,\cL}:(\Fr^n)^*\cL\rar{\simeq}\cL$
which acts as the identity on the stalk of $\cL$ at $1\in
G(\bF_q)$. In view of Remark \ref{r:charsheaves}(c), it is clear
that
\[
\psi_{k,\cL} = \psi_{n,\cL}\circ (\Fr^n)^*(\psi_{n,\cL}) \circ
(\Fr^{2n})^*(\psi_{n,\cL}) \circ \dotsb \circ
(\Fr^{k-n})^*(\psi_{n,\cL}) : (\Fr^k)^*\cL\rar{\simeq}\cL
\]
whenever $n\lvert k$. To complete the verification of property
(4-iii), we prove three lemmas.
\begin{lem}
The trace function $t_{n,\cL}:G(\bF_{q^n})\to\ql$ determined by
$\psi_{n,\cL}$ $($as in \S\ref{ss:charsheaves-introduction}$)$
takes values in $\ql^\times$, and is a $1$-dimensional character
of $G(\bF_{q^n})$.
\end{lem}
\begin{proof}
By Remark \ref{r:charsheaves}(c), our normalization of
$\psi_{n,\cL}$ implies that for any isomorphism
$\nu:\al^*\cL\rar{\simeq}\cL\boxtimes\cL$, the following diagram
commutes:
\[
\xymatrix{
  \al^*(\Fr^n)^*\cL = (\Fr^n)^*\al^*\cL \ar[rr]^{(\Fr^n)^*\nu\ \ \ \ \ \ \ \ }
  \ar[d]_{\al^*\psi_{n,\cL}}  &
  &
  (\Fr^n)^*(\cL\boxtimes\cL) = (\Fr^n)^*\cL\boxtimes(\Fr^n)^*\cL
  \ar[d]^{\psi_{n,\cL}\boxtimes\psi_{n,\cL}}
  \\
  \al^*\cL \ar[rr]^{\nu} &   &  \cL\boxtimes\cL
   }
\]
(By abuse of notation, the equality signs in this diagram stand
for the canonical isomorphisms.) Now let
\[
t_{n,\al^*\cL} : G(\bF_{q^n})\times G(\bF_{q^n}) \rar{} \ql
\]
denote the trace function associated to the isomorphism
$\al^*\psi_{n,\cL}:(\Fr^n)^*\al^*\cL\rar{\simeq}\al^*\cL$, as in
\S\ref{ss:charsheaves-introduction}, and let
\[
t_{n,\cL\boxtimes\cL} : G(\bF_{q^n})\times G(\bF_{q^n}) \rar{} \ql
\]
denote the trace function associated to the isomorphism
$\psi_{n,\cL}\boxtimes\psi_{n,\cL}:(\Fr^n)^*(\cL\boxtimes\cL)\rar{\simeq}\cL\boxtimes\cL$.
Then it is clear that $t_{n,\al^*\cL}=\al^*t_{n,\cL}$, the
pullback of the function $t_{n,\cL}$ by the map
$\al:G(\bF_{q^n})\times G(\bF_{q^n})\to G(\bF_{q^n})$, whereas
$t_{n,\cL\boxtimes\cL}=t_{n,\cL}\boxtimes t_{n,\cL}$, the function
defined by
\[
(t_{n,\cL}\boxtimes t_{n,\cL})(g_1,g_2) = t_{n,\cL}(g_1)\cdot
t_{n,\cL}(g_2) \qquad \forall\, g_1,g_2\in G(\bF_{q^n}).
\]
The commutativity of the diagram above implies that
$t_{n,\cL}\boxtimes t_{n,\cL}=\al^* t_{n,\cL}$, which is
equivalent to the statement that $t_{n,\cL}:G(\bF_{q^n})\to\ql$ is
multiplicative. Since $t_{n,\cL}(1)=1$ by construction, it follows
that $t_{n,\cL}:G(\bF_{q^n})\to\ql^\times$ is in fact a
($1$-dimensional) character.
\end{proof}

The attentive reader may observe that so far we have not used the
assumption that $G$ is commutative. However, it will be used in
the next lemma.

\begin{lem}
Every $1$-dimensional character $\chi:G(\bF_{q^n})\to\ql^\times$
has the form $t_{n,\cL}$ for some $\cL\in\cs(G)^{\Fr^n}$.
\end{lem}
\begin{proof}
Let $L_n:\Gb\to\Gb$ denote the $n$-th Lang isogeny, defined by
$L_n(\ga)=\Fr^n(\ga)-\ga$ (we denote the group operation on $\Gb$ additively,
as in \S\ref{ss:characters-abelian-groups}). By Lang's theorem \cite{lang},
this map makes $\Gb$ a $G(\bF_{q^n})$-torsor over $\Gb$. In other words, it
induces an isomorphism between $\Gb/G(\bF_{q^n})$ and $\Gb$. Standard \'etale
descent theory implies that pullback via $L_n$ induces an equivalence between
the category of local systems on $\Gb$ and the category of
$G(\bF_{q^n})$-equivariant local systems on $\Gb$, where $G(\bF_{q^n})$ acts on
$\Gb$ by translations. In particular, there exists a (unique up to isomorphism)
rank $1$ local system $\cL$ on $\Gb$ such that $L_n^*\cL$ is the constant local
system $\ql$ equipped with the $G(\bF_{q^n})$-equivariant structure defined by
the homomorphism $\chi:G(\bF_{q^n})\to\ql^\times$.

\mbr

We claim that $\cL\in\cs(G)^{\Fr^n}$, and that $t_{n,\cL}=\chi^{-1}$.
First of all, the commutativity of $G$ implies that
$\al\circ(L_n\times L_n)=L_n\circ\al:\Gb\times\Gb\to\Gb$. Thus
$\al^*\cL$ is the unique local system on $\Gb\times\Gb$ such that
$(L_n\times L_n)^*\al^*\cL$ is the trivial rank $1$ local system
on $\Gb\times\Gb$ equipped with the $G(\bF_{q^n})\times
G(\bF_{q^n})$-equivariant structure given by the character
$\al^*\chi$. Similarly, $(L_n\times L_n)^*(\cL\boxtimes\cL)$ is
the trivial rank $1$ local system on $\Gb\times\Gb$ equipped with
the $G(\bF_{q^n})\times G(\bF_{q^n})$-equivariant structure given
by the character $\chi\boxtimes\chi$. But
$\chi\boxtimes\chi=\al^*\chi$, and it follows from the fact
mentioned in the previous paragraph that
$\al^*\cL\cong\cL\boxtimes\cL$.

\mbr

Similarly, it is easy to construct an isomorphism
$\psi_{n,\cL}:(\Fr^n)^*\cL\rar{\simeq}\cL$. Indeed, $L_n$ commutes with
$\Fr^n$, so
$L_n^*(\Fr^n)^*\cL\cong(\Fr^n)^*L_n^*\cL\cong(\Fr^n)^*\ql\cong\ql\cong
L_n^*\cL$ (where $\ql$ denotes the constant rank $1$ local system on $\Gb$ and
all isomorphisms are canonical). Moreover, multiplication by $G(\bF_{q^n})$
commutes with $\Fr^n$, which implies that the isomorphisms above are in fact
isomorphisms of $G(\bF_{q^n})$-equivariant local systems, where $\ql$ is
equipped with the equivariant structure given by the character $\chi$. Thus the
isomorphisms above descend to the desired isomorphism
$\psi_{n,\cL}:(\Fr^n)^*\cL\rar{\simeq}\cL$. The last observation is that if
$g\in G(\bF_{q^n})$ and $\ga\in G(\bF)$ is such that $g=L_n(\ga)$, then, by
definition, $\ga=\Fr^n(\ga)-g$, whence $\psi_{n,\cL}$ acts as multiplication by
$\chi(g)^{-1}$ on the stalk $\cL_g$. Thus $t_{n,\cL}=\chi^{-1}$. Of course, we
can replace $\chi$ with $\chi^{-1}$, and the proof of the lemma is complete.
\end{proof}

\begin{lem}
If $\cL,\cK\in\cs(G)^{\Fr^n}$ and $t_{n,\cL}=t_{n,\cK}$, then
$\cL\cong\cK$.
\end{lem}
\begin{proof}
In view of the argument used in the proof of the previous lemma,
it suffices to show that if $\cL\in\cs(G)^{\Fr^n}$, then
$L_n^*\cL$ is a trivial rank $1$ local system on $\Gb$. (Indeed,
it then follows that $\cL$ is determined up to isomorphism by the
$G(\bF_{q^n})$-equivariant structure on $L_n^*\cL$, which is given
by the character $t_{n,\cL}^{-1}$.) Let $\iota:\Gb\to\Gb$ denote the
inversion map, and let $\De:\Gb\to\Gb\times\Gb$ denote the
diagonal. By definition,
$L_n=\al\circ(\Fr^n\times\iota)\circ\De:\Gb\to\Gb$. If
$\cL\in\cs(G)$, then $\al^*\cL\cong\cL\boxtimes\cL$, which implies
that $L_n^*\cL\cong\bigl((\Fr^n)^*\cL\bigr)\tens\cL^{-1}$. If, in
addition, $(\Fr^n)^*\cL\cong\cL$, then $L_n^*\cL$ is trivial, as
desired.
\end{proof}

\mbr

It is clear that the three lemmas above imply that if $G$ is a
connected commutative algebraic group over $\bF_q$, then the
collection $CS(G)$ of local systems on $\Gb$ that we have defined
at the beginning of this subsection satisfies property (4-iii) of
\S\ref{ss:charsheaves-introduction}. Thus question (4) has an
affirmative answer in this case.


\section{The orbit method for finite nilpotent groups}\label{s:orbmethod}

The orbit method was originally discovered by A.A.~Kirillov \cite{kirillov} for
connected and simply connected real nilpotent Lie groups. In this section we
discuss the orbit method in a different context of {\em finite} nilpotent
groups of ``small'' nilpotence class\footnote{The fact that the orbit method
works in this context was probably clear as soon as the original orbit method
over $\bR$ was discovered. The usefulness of the orbit method for finite groups
was demonstrated by D.~Kazhdan \cite{kazhdan}.}. In Section
\ref{s:orbmethod-unipotent} we will apply the orbit method to the groups of the
form $\Ga_n=G(\bF_{q^n})$, where $G$ is a unipotent algebraic group over
$\bF_q$ whose {\em nilpotence class} (defined in \S\ref{ss:alg-aspects}) is less than $p$.

\mbr

For more information about the orbit method for nilpotent Lie
groups we refer to \cite{kirillov} and \cite{solvable} (the latter
also contains results about more general solvable Lie groups). A
detailed exposition of the orbit method for finite nilpotent
groups appears in \cite{orbmethod}.

\subsection{Lazard's construction}\label{ss:lazard}
Fix $c\in\bN$. Let $\Nilp_c$ be the category of (possibly infinite) nilpotent
groups $\Ga$ of nilpotence class $\leq c$ such that the map $\Ga\to\Ga$,
$g\mapsto g^k$, is invertible for all $k\leq c$. Let $\nilp_c$ be the category
of nilpotent Lie algebras over $\bZ\bigl[\frac{1}{c!}\bigr]$ of nilpotence
class $\leq c$. Michel Lazard \cite{lazard} constructed a canonical isomorphism
of categories $\Exp:\nilp_c\rar{\simeq}\Nilp_c$. The functor $\Exp$ is defined
as follows. The group $\Exp\fg$ has the same underlying set as the Lie algebra
$\fg\in\nilp_c$, and the group operation on the set $\Exp\fg=\fg$ is defined by
the formula
\begin{equation}\label{e:(**)}
x*y := \sum_{i\leq c} CH_i(x,y), \qquad x,y\in\fg,
\end{equation}
where $CH_i$ is the homogeneous component of degree $i$ of the
Campbell-Hausdorff series $CH(X,Y):=\log\bigl(\exp(X)\exp(Y)\bigr)$. (For
instance, $CH_2(X,Y)=[X,Y]/2$.) It is known that $CH_i$ is a Lie polynomial
with coefficients in $\bZ\bigl[\frac{1}{i!}\bigr]$ (by definition, this means
that $CH_i$ belongs to the free Lie algebra over $\bZ\bigl[\frac{1}{i!}\bigr]$
with generators $X,Y$). So the RHS of \eqref{e:(**)} is well defined. It is
easy to see that $\Exp\fg$ is a group of nilpotence class $\leq c$, and that if
$f:\fg_1\to\fg_2$ is a morphism in $\nilp_c$, then $f$, viewed as a map
$\Exp\fg_1\to\Exp\fg_2$, is a group homomorphism. So $\Exp$ is a functor
$\nilp_c\to\Nilp_c$.
\begin{thm}[M.~Lazard]\label{t:lazard}
$\Exp:\nilp_c\to\Nilp_c$ is an isomorphism of categories.
\end{thm}
For the proof, see \cite{lazard} or \cite{kh}, \S10.2. The inverse functor to
$\Exp$ will be denoted by $\Log$. If $\fg\in\nilp_c$ and
$\Ga=\Exp\fg\in\Nilp_c$, we prefer to denote the identity maps between $\fg$
and $\Ga$ by $\exp:\fg\to\Ga$ and $\log:\Ga\to\fg$, respectively. Thus our
notation is compatible with the one used in classical Lie theory.

\subsection{Orbit method: formulation} Let $\Ga$ be a finite group of
nilpotence class $\leq c$ such that all the prime divisors of the
order $\abs{\Ga}$ are greater than $c$. Applying Lazard's
construction (see \S\ref{ss:lazard}) we get a finite Lie ring
$\fg:=\Log\Ga\in\nilp_c$ and a bijection
$\log:\Ga\rar{\simeq}\fg$. The abelian group
\[
\fg^*:=\Hom(\fg,\cst)
\]
carries a natural action of $\Ga$, called the {\em coadjoint action}. The {\em
orbit method} in its simplest incarnation provides a natural bijection between
$\gh$ and the set $\fg^*/\Ga$ of $\Ga$-orbits on $\fg^*$.

\begin{thm}\label{t:orbmethod}
For every $\Ga$-orbit $\Om\subset\fg^*$ there exists $\rho_\Om\in\gh$ whose
character is given by the formula
\begin{equation}\label{e:char-formula}
\tr\bigl(\rho_\Om(g)\bigr) = \frac{1}{(\operatorname{card}\Om)^{1/2}}
\sum_{f\in\Om} f(\log g) \qquad \forall\,g\in\Ga.
\end{equation}
The map $\Om\mapsto\rho_\Om$ is a bijection $\fg^*/\Ga\rar{\simeq}\gh$.
\end{thm}

\begin{rems}\label{r:char-formula}
\begin{enumerate}[(i)]
\item Formula \eqref{e:char-formula} means that the character of $\rho_\Om$
equals the (inverse) Fourier transform of the characteristic function of $\Om$
up to a constant factor.
\item It follows from \eqref{e:char-formula} that
\begin{equation}\label{e:dim-rho-Om}
\dim\rho_\Om = (\operatorname{card}\Om)^{1/2}.
\end{equation}
\item By \eqref{e:dim-rho-Om}, $\operatorname{card}\Om$ is a full square.
There is also a direct proof of this fact (see \cite{orbmethod}),
very much in the spirit of the proof of the corresponding well
known fact from Lie theory: the dimension of any coadjoint orbit
for every Lie group is even \cite{kirillov}.
\item As a representation, $\rho_\Om$ is defined only up to a
{\em non-unique}\footnote{If $V$ is the space of the representation $\rho_\Om$
and $\Ga_\Om=\bigl\{g\in\Aut(\Ga) \,\st\, g\Om=\Om \bigr\}$, then it is not
always possible to make the projective representation of $\Ga_\Om$ in $V$ a
true representation. In fact, $\Ga_\Om\supset\Ga /Z$, where $Z$ is the center
of $\Ga$, and the  projective representation of $\Ga /Z$ in $V$ can be made a
true representation only if $Z\cap [\Ga ,\Ga ]\subset\Ker\rho_{\Om}$. This
condition does not hold if $\Ga$ is the Heisenberg group and
$\operatorname{card}\Om>1$.} isomorphism. However, it follows from
\cite{lion-perrin} that, given $f\in\Om$, there is a way to construct a
representation $\rho_f$ with class $\rho_\Om\in\gh$ which is defined up to a
{\em unique} isomorphism.
\end{enumerate}
\end{rems}
The proof of Theorem \ref{t:orbmethod} given in \S\S\ref{ss:fourier},
\ref{ss:char-formula} below is well known, but we were unable to find it in the
literature. Kirillov's original proof \cite{kirillov} was different. It is
longer but allows one to construct the irreducible representations of $\Ga$
(not only their characters), and it works for nilpotent groups over $\bR$ or
$\bQ_p$. His proof is based on the notion of a {\em polarization}, which will
be recalled in \S\ref{ss:Polarizations}. Some of the references in the
literature that adapt Kirillov's method to the situation of Theorem
\ref{t:orbmethod} are \cite{kazhdan,basmanova,matveev,orbmethod}.

\subsection{Construction of the orbit method bijection}\label{ss:fourier}
In this subsection we will construct a bijection
\begin{equation}\label{e:*}
\fg^*/\Ga \rar{\simeq} \gh;
\end{equation}
in \S\ref{ss:char-formula} we will prove \eqref{e:dim-rho-Om}, and then
\eqref{e:char-formula}. One has
\begin{equation}\label{e:specs}
\fg^*/\Ga = \Spec \Fun(\fg^*)^\Ga, \qquad \gh=\Spec\Z\Ga,
\end{equation}
where $\Fun(\fg^*)$ is the algebra of functions $\fg^*\to\bC$ with pointwise
multiplication, and $\Z\Ga$ is the center of the group algebra $\bC\Ga$ of
$\Ga$. Note that $\bC\Ga$ is naturally isomorphic to the algebra $\Meas(\Ga)$
of complex measures on $\Ga$ with convolution as the multiplication operation,
and under this isomorphism $\Z\Ga$ corresponds to the subalgebra
$\Meas(\Ga)^\Ga\subseteq\Meas(\Ga)$ of measures invariant under conjugation.
The Fourier transform provides a $\Ga$-equivariant isomorphism of algebras
\[
\cF:\Meas(\fg)\rar{\simeq}\Fun(\fg^*),
\]
\[
\cF(\mu)(f) = \int_\fg f(x)^{-1}\,d\mu(x),
\]
so constructing
\eqref{e:*} amounts to constructing an algebra isomorphism
\begin{equation}\label{e:**}
\Meas(\fg)^\Ga \rar{\simeq} \Meas(\Ga)^\Ga.
\end{equation}
The bijection $\exp:\fg\to\Ga$ induces an isomorphism {\em of vector spaces}
\begin{equation}\label{e:***}
\exp_* : \Meas(\fg) \rar{\simeq} \Meas(\Ga).
\end{equation}
We define \eqref{e:**} to be the restriction of \eqref{e:***}. Notice that if
$\Ga$ is not commutative, then $\exp:\fg\to\Ga$ is not a group homomorphism,
and therefore \eqref{e:***} is not an algebra morphism.
\begin{prop}\label{p:alg-isom}
The map \eqref{e:**} we have defined is an isomorphism of algebras.
\end{prop}

A complete proof of this proposition is contained in \cite{orbmethod}. It rests
on Lemma \ref{l:phi-psi} below. We leave the proof of Proposition
\ref{p:alg-isom} as an exercise for the reader (it can be easily done by
imitating the proof of Proposition \ref{p:convol}).

\begin{lem}\label{l:phi-psi}
For every $c\in\bN$, there exist Lie polynomials $\phi(X,Y)$ and $\psi(X,Y)$
with coefficients in $\bZ\left[\frac{1}{c!}\right]$ such that
\begin{equation}\label{e:phipsi}
\log\bigl( \exp(X)\cdot\exp(Y) \bigr) \equiv \exp\bigl(\ad\phi(X,Y)\bigr)(X) +
\exp\bigl(\ad\psi(X,Y)\bigr)(Y)
\end{equation}
modulo terms of degree $\geq c+1$.
\end{lem}

\sbr

\noindent
(The statement means that $\phi(X,Y)$ and $\psi(X,Y)$ are elements of the free
Lie algebra over $\bZ\left[\frac{1}{c!}\right]$ generated by $X$ and $Y$, and
that the difference between the two sides of \eqref{e:phipsi} is a formal Lie
series in $X$ and $Y$ all of whose terms have degree $\geq c+1$. Note that if
$\fg\in\nilp_c$, then $\phi(x,y)$ and $\psi(x,y)$ can be evaluated in $\fg$ for
all $x,y\in\fg$.)

\sbr

\begin{proof}[Sketch of the proof]
One can construct $\phi$ and $\psi$ inductively. Let us write
$\phi(X,Y)=\sum_{n=1}^{c-1} \phi_n(X,Y)$ and $\psi(X,Y)=\sum_{n=1}^{c-1}
\psi_n(X,Y)$, where $\phi_n$ and $\psi_n$ are homogeneous Lie polynomials of
degree $n$. Then at each step of the induction, we will have to solve an
equation of the form $[X,\phi_n]+[Y,\psi_n]\equiv\eta_{n+1}(X,Y)$ modulo
commutators of order $\geq n+2$, where $\eta_{n+1}(X,Y)$ is already known. It
is easy to check that this is always possible. As the proof shows, we use
nothing special about the Campbell-Hausdorff series: as a matter of fact, the
latter could be replaced by any formal Lie series whose linear part is $X+Y$
and whose terms of degree $\leq c$ have coefficients in
$\bZ\left[\frac{1}{c!}\right]$.
\end{proof}

Proposition \ref{p:alg-isom} means that the composition
\begin{equation}\label{e:composition}
\Phi:\bC \Ga \cong \Meas(\Ga) \xrar{\log_*} \Meas(\fg) \rar{\cF} \Fun(\fg^*)
\end{equation}
restricts to an algebra isomorphism $\Z\Ga\rar{\simeq}\Fun(\fg^*)^\Ga$. In view
of \eqref{e:specs}, this yields a bijection $\fg^*/\Ga\rar{\simeq}\gh$, as
desired. This bijection will be denoted by $\Om\mapsto\rho_\Om$. We warn the
reader that in general, the composition \eqref{e:composition} {\em is not an
isomorphism of $\Z\Ga$-modules}, where the $\Z\Ga$-module structure on
$\Fun(\fg^*)$ is induced by the isomorphism $\Z\Ga\rar{\simeq}\Fun(\fg^*)^\Ga$
(see Appendix \ref{a:counterexamples} for a counterexample).
\begin{rem}\label{r:idemp}
The composition \eqref{e:composition} takes the central idempotent
$e_\Om\in\Z \Ga$ corresponding to the representation $\rho_\Om$ to
the characteristic function $1_\Om$ of the subset
$\Om\subset\fg^*$.
\end{rem}

\subsection{Proof of the character formula}\label{ss:char-formula}
Let $\chi_{reg}:\bC\Ga\to\bC$ denote the character of the regular
representation of $\Ga$, and $\Phi:\bC\Ga\to\Fun(\fg^*)$ the composition
\eqref{e:composition}.
\begin{lem}\label{l:2.7}
We have $\chi_{reg}=\Sg\circ\Phi$, where $\Sg:\Fun(\fg^*)\to\bC$
takes $f\in\Fun(\fg^*)$  to $\sum_{\la\in\fg^*}f(\la)$.
\end{lem}
\begin{proof}
We have $\chi_{reg}(1)=\operatorname{card}\Ga=(\Sg\circ\Phi)(1)$. If
$\ga\in\Ga\setminus\{1\}$, then $\chi_{reg}(\ga)=0=(\Sg\circ\Phi)(\ga)$. Since
$\chi_{reg}$ and $\Sg\circ\Phi$ are linear, the proof is complete.
\end{proof}

\mbr

Let $e_\Om$ and $\rho_\Om$ be as in Remark \ref{r:idemp}. Then
$\chi_{reg}(e_\Om)=(\dim\rho_\Om)^2$. Combining this with Lemma \ref{l:2.7} and
Remark \ref{r:idemp}, we obtain
$(\dim\rho_\Om)^2=\Sg(1_\Om)=\operatorname{card}\Om$, which is equivalent to
\eqref{e:dim-rho-Om}. To prove \eqref{e:char-formula}, notice that the
character of $\rho_\Om$, call it $\chi_\Om:\bC\Ga\to\bC$, is the linear
functional invariant under $\Ga$-conjugation whose value at $e_{\Om'}$ equals
$\tr\rho_\Om(e_{\Om'})$, i.e., $\dim\rho_\Om$ if $\Om'=\Om$ and $0$ if
$\Om'\neq\Om$. By Remark \ref{r:idemp}, $\chi_\Om=\ell_\Om\circ\Phi$, where
$\ell_\Om:\Fun(\fg^*)\to\bC$ is the $\Ga$-invariant linear functional such that
$\ell_\Om(1_{\Om'})=0$ for $\Om'\neq\Om$ and
$\ell_\Om(1_\Om)=\dim\rho_\Om=(\operatorname{card}\Om)^{1/2}$. This is
equivalent to \eqref{e:char-formula}.

\subsection{Polarizations}\label{ss:Polarizations} To end the
section we briefly explain a different approach to the orbit
method which has the advantage of explicitly constructing
irreducible representations, not only their characters.
Furthermore, this approach works not only in the setting of this
section, but also in the setting of nilpotent Lie groups
\cite{kirillov}.

\begin{defin}\label{d:polarization}
\begin{enumerate}[(i)]
\item Let $\fg$ be a nilpotent Lie algebra over a field $k$ and $f\in\fg^*=\Hom_k(\fg,k)$.
A {\em polarization} of $\fg$ at $f$ is a Lie subalgebra $\fh\subseteq\fg$ such
that $f([\fh,\fh])=\{0\}$, and $\fh$ is maximal among all {\em subspaces} of
$\fg$ with this property.
\item If $\fg$ is a finite nilpotent Lie ring and $f\in\fg^*=\Hom_\bZ(\fg,\cst)$, a
{\em polarization} of $\fg$ at $f$ is a Lie subring $\fh\subseteq\fg$ such that
$f([\fh,\fh])=\{1\}$, and $\fh$ is maximal among all {\em additive subgroups}
of $\fg$ with this property.
\end{enumerate}
\end{defin}

\begin{thm}
In each of the situations $($i$)$ and $($ii$)$, $\fg$ has a polarization at
every $f\in\fg^*$.
\end{thm}
This theorem is due to A.A.~Kirillov \cite{kirillov}, and a somewhat different
proof was given by M.~Vergne \cite{vergne}. We explain her proof in Appendix
\ref{a:vergne}.

\begin{rem}
The theorem implies that one can alternately define a polarization to be a
subalgebra $\fh\subseteq\fg$ of maximal possible dimension satisfying
$f([\fh,\fh])=0$; this is the approach of \cite{kirillov}. On the other hand,
it is {\em not true} that a maximal (with respect to inclusion) subalgebra
$\fh\subseteq\fg$ satisfying $f([\fh,\fh])=0$ is necessarily a polarization:
see \S\ref{aa:1-and-2}.
\end{rem}

The following theorem, due to A.A.~Kirillov \cite{kirillov}, explains how to
construct an irreducible representation corresponding to a given coadjoint
orbit.
\begin{thm}
In the situation of Theorem \ref{t:orbmethod}, let $f\in\Om$, and let
$\fh\subseteq\fg$ be a polarization at $f$. Put $H=\Exp\fh$, and let
$\chi_f:H\to\cst$ be the character corresponding to $f$. Then the
representation $\rho_{f,\fh}:=\Ind_H^\Ga \chi_f$ is irreducible, and its
isomorphism class equals $\rho_\Om$.
\end{thm}
For a proof of this result and for more details we refer the
reader to \cite{orbmethod}.


\section{Unipotent groups over $\bF_q$ and $\bL$-indistinguishability}\label{s:orbmethod-unipotent}

In this section we apply the orbit method discussed above to study the
elementary questions posed in \S\ref{ss:main-geom-questions} for unipotent
groups of small nilpotence class.

\subsection{Lazard's construction for unipotent groups}\label{ss:lazard-unipotent}
Let $k$ be a field of characteristic $p>0$, and let $G$ be a
unipotent algebraic group over $k$ which has nilpotence class
$\leq p-1$. For every scheme $S$ over $k$, we get the (possibly
infinite) group $G(S)$ that satisfies the assumptions of Theorem
\ref{t:lazard} with $c=p-1$, and hence Lazard's construction
produces a Lie ring $\fg(S)$ (in fact, a Lie algebra over
$\bZ\bigl[\frac{1}{(p-1)!}\bigr]$). Thus we get a functor $\fg$
from schemes over $k$ to Lie rings, and this functor is obviously
representable since the functors \{$k$-schemes\} $\rar{}$ \{Sets\}
corresponding to $G$ and $\fg$ are the same by definition. Hence
$\fg$ is a Lie ring scheme over $k$. Moreover, the underlying
additive group of $\fg$ is unipotent, because $G$ has a filtration
by normal subgroups with commutative unipotent successive
quotients, which induces a similar filtration on the additive
group of $\fg$. We will write $\fg=\Log G$. Conversely, if we
start with a nilpotent Lie ring scheme $\fg$ over $k$ such that
the underlying additive group of $\fg$ is unipotent and the
nilpotence class of $\fg$ is less than $p$, we can form the
associated unipotent group scheme $G=\Exp\fg$. These two
constructions will also be referred to as Lazard's constructions.

\begin{rem}
All of the above goes through equally well in the case when
$\operatorname{char}k=0$. In this situation no assumption on the
nilpotence class of $G$ or $\fg$ is needed. However, Lazard's
construction is not very useful in this case because the Lie ring
scheme $\fg$ it produces always arises from an honest nilpotent
Lie algebra over $k$, which coincides with the Lie algebra of the
algebraic group $G$ in the classical sense. On the other hand, it
is very important to notice that if $\operatorname{char}k>0$, then
$\Log G$ may be very different from the Lie algebra of $G$, even
if $G$ is annihilated by $p$. An example is provided by the fake
Heisenberg groups defined in \S\ref{ss:fake-heisenberg}, for which
$\Log G$ does not correspond to \emph{any} Lie algebra over $k$.
\end{rem}

\subsection{The dual of $\fg$}\label{ss:equivar-duality}
Now let us assume that $k$ is perfect and $G$ is connected (this implies that
$G$ is geometrically connected, see the footnote to Proposition
\ref{p:base-change-abelian}). Then $\fg$ is connected, so one has its Serre
dual $\fg^*$ as an object of the category $\cC_k$ from
\S\ref{ss:perfect-duality}.\footnote{If $\fg$ is disconnected, then $\fg^*$ is
not a group scheme, but rather a Picard stack.} The general definition of
$\fg^*$ is given in Appendix \ref{a:duality}, but the reader may prefer to
assume that $p\cdot\fg=0$, and to use the {\em ad hoc} definition of $\fg^*$
from \S\ref{ss:perfect-duality} (all phenomena discussed in this article
already appear in this special case).

\mbr

In some sense, one also has the coadjoint action of $G$ on $\fg^*$. More
precisely, the situation is as follows. Recall that the category $\cC_k$ in
which $\fg^*$ lives is a certain localization of the category
$\widetilde{\cC}_k$ of connected commutative unipotent algebraic groups over
$k$, and there is no canonical way of choosing a representative of $\fg^*$ in
$\widetilde{\cC}_k$. Thus, strictly speaking, one should denote such a
representative by a different symbol, such as $\fg'$, as we have done in
Appendices \ref{a:duality} and \ref{a:fourier-deligne}. Furthermore, in
Corollary \ref{c:representative-dual} we show that there exists a
representative $\fg'\in\widetilde{\cC}_k$ of $\fg^*$ and an action of $G$ on
$\fg'$ by group automorphisms which induces the canonical action on $\fg^*$.
However, in the main body of the paper we will denote by $\fg^*$ a chosen
representative of the Serre dual of $\fg$ on which $G$ acts.

\mbr

A simple and instructive example of a coadjoint action will be
discussed in \S\ref{ss:fake-heisenberg}.

\subsection{The orbit method for unipotent groups}\label{ss:orbmet-unip}
From now on assume that $k=\bF_q$, so $G$ is a connected unipotent algebraic
group over $\bF_q$ of nilpotence class $<p$.  Fix a representative
$\fg^*\in\widetilde\cC_k$ of the Serre dual of $\fg=\Log G$ so that $G$ acts on
$\fg^*$ (see \S\ref{ss:equivar-duality}).

\mbr

Let us consider questions (1)-(3) from
\S\ref{ss:main-geom-questions} in this situation. The answer to
question (1) is positive. Indeed, the base change maps
\begin{equation}\label{e:bsmaps}
 \gh_m \rar{} \bigl(\gh_n\bigr)^{\Gal(\bF_{q^n}/\bF_{q^m})}, \qquad m\lvert n,
\end{equation}
can be defined as follows. By Proposition \ref{p:lang}, we have
natural identifications $\fg^*(\bF_{q^n})\cong\fg(\bF_{q^n})^*$
for all $n\in\bN$. We define \eqref{e:bsmaps} to be the
composition
\[
\gh_m = \fg^*(\bF_{q^m})/\Ga_m \rar{} \fg^*(\bF_{q^n})/\Ga_n = \gh_n.
\]
Here the arrow is induced by the inclusion
$\fg^*(\bF_{q^m})\hookrightarrow\fg^*(\bF_{q^n})$ and the
equalities stand for the natural bijections provided by the orbit
method. The map \eqref{e:bsmaps} is clearly $\Fr$-equivariant. The
set $\gh:=\indlim{n} \gh_n$ is naturally identified with
$\fg^*(\bF)/G(\bF)$.

\mbr

Now consider question (3) from \S\ref{ss:main-geom-questions}. It
amounts to whether $\fg^* (\bF )/G(\bF )$ can be interpreted as
the set of $\bF$-points of some kind of geometric object over
$\bF_q$. We do not know a good answer\footnote{Here is an answer
that has {\it a chance\,} to be good. If $X$ is a $G$-scheme of
finite type over a field $k$ one can define $X/G$ as an
\emph{$h$-sheaf}, i.e., as a sheaf on the category of $k$-schemes
of finite type equipped with Voevodsky's h-topology (see the
appendix to \cite{suslin-voevodsky}). Namely, $X/G$ is the
$h$-sheaf associated to the presheaf $S\mapsto X(S)/G(S)$. This
definition of $X/G$ seems to be reasonable at least if the
stabilizers of all points of $X$ have the same dimension.}. But at
least, there exists a stratification of $\fg^*$ by $G$-invariant
locally closed subschemes $S_r$ such that for each $r$ the
subscheme $\{ (g,s)\in G\times S_r | gs=s\}$ is flat over $S_r$;
then the quotient $S_r/G$ exists as an algebraic space and
$(S_r/G)(\bF )=S_r(\bF )/G(\bF )$.

\subsection{The $\bL$-indistinguishability
phenomenon}\label{ss:l-indistinguishability}
Let $G$ be as in \S\ref{ss:orbmet-unip}. We will show
that in general, the answer to question (2) from \S\ref{ss:main-geom-questions} is
negative: in terms of Definition \ref{d:l-indistinguishability}, there may exist
$\bL$-packets with more than one element. We will also describe  the $\bL$-packets
(see Proposition \ref{p:l-indis2}(iii) and Remark \ref{r:l-packets}).

\begin{prop}\label{p:l-indist}
The maps  $\gh_n\rar{} \bigl(\gh\bigr)^{\Gal(\bF/\bF_{q^n})}$ are surjective.
They are injective for all $n$ if and only if the $G$-stabilizer of each point
of $\fg^*$ is connected.
\end{prop}

\sbr

In \S\ref{ss:fake-heisenberg} we will give an example where
the stabilizers of some points of $\fg^*$ are disconnected.
Proposition \ref{p:l-indist} is a consequence of the discussion in \S\ref{ss:orbmet-unip}
and of statements (i)-(ii) of the following proposition applied to $X=\fg^*$.
Statement (iii) for $X=\fg^*$ gives a description of $\bL$-packets.

\sbr

\begin{prop}\label{p:l-indis2}
Let $X$ be an algebraic variety over $\bF_q$ equipped with an
action of a connected algebraic group $G$.
\begin{enumerate}[$($i$)$]
\item The natural map
\begin{equation}\label{e:map-orbits}
f: X(\bF_q)/G(\bF_q) \rar{} \bigl(X(\bF)/G(\bF)\bigr)^{\Fr}
\end{equation}
is surjective.
\item The map $X(\bF_{q^n})/G(\bF_{q^n})\to X(\bF)/G(\bF)$ is injective for each $n$ if and only if
the $G$-stabilizer of each point of $X$ is connected.
\item For each $x\in X(\bF_q)$ there is a canonical bijection $f^{-1}(f(\overline{x}))\isom H^1(\bF_q,\pi_0(G_x))$.
Here $G_x\subset G$ is the stabilizer of $x$, $\pi_0(G_{x}):=G_{x}/G_x^\circ$,
$\overline{x}$ is the image of $x$ in $X(\bF_q)/G(\bF_q)$, and $G_x^\circ$ is
the connected component of the identity in $G_x$.
\end{enumerate}
\end{prop}

This proposition is standard. Its proof will be recalled in
\S\ref{ss:proof-l-indis2}. In the next subsection we describe a class of
connected unipotent groups for which all $\bL$-packets are trivial.

\begin{rem}[A gerby description of the $\bL$-packets]\label{r:l-packets}
A drawback of the description of the fiber $f^{-1}(f(\overline{x}))$ given in
Proposition \ref{p:l-indis2}(iii) is that it depends on the choice of $x$, not
only on $y=f(\overline{x})$. However, there exists a more canonical
description. Namely, let us first consider the groupoid $\sG_y$ whose set of
objects is $\{x_0\in X(\bF_q)\bigl\lvert f(\overline{x_0})=y\}$ and where a
morphism $x_0\to x_1$ is an $\bF_q$-point of $\pi_0\bigl(\{g\in
G\,\bigl\lvert\, gx_0=x_1\}\bigr)$. By Lang's theorem (applied to the neutral
component $G_{x_0}^\circ$ of $G_{x_0}$) such an $\bF_q$-point can be lifted to
an $\bF_q$-point of $\bigl\{g\in G\,\bigl\lvert\,gx_0=x_1\bigr\}$. Thus the
isomorphism classes of objects of $\sG_y$ are in one-to-one correspondence with
the elements of $f^{-1}(y)$.

\mbr

The groupoid $\sG_y$ can in turn be described in terms of another groupoid
$\sA_y$ whose set of objects is $\{x_0\in X(\bF_q)\bigl\lvert
f(\overline{x_0})=y\}$ and where a morphism $x_0\to x_1$ is an $\bF$-point of
$\pi_0\bigl(\{g\in G\,\bigl\lvert\, gx_0=x_1\}\bigr)$. Note that $\sA_y$ is a
gerbe over a point (i.e., a gerbe over $\bF$). That is, $\sA_y$ has precisely
one isomorphism class of objects. Moreover, $\sA_y$ is equipped with a (strict)
action of $\Fr$, and $\sG_y$ can be identified with the groupoid of pairs
$(c,\Fr(c)\rar{\simeq}c)$ consisting of an object $c$ of $\sA_y$ and as
isomorphism between $\Fr(c)$ and $c$.
\end{rem}

\subsection{Groups of exponential type}\label{ss:exp-type}
If $\fg$ is an honest Lie algebra (as opposed to a Lie ring scheme) over
$\bF_q$ of nilpotence class $<p$, then the unipotent algebraic group
$G=\Exp\fg$ produced by Lazard's construction is said to be of {\em exponential
type}. In this case one can choose the group scheme $\fg^*$ from
\S\S\ref{ss:equivar-duality}-\ref{ss:l-indistinguishability} to be equal to the
vector space dual to $\fg$.

\begin{prop}\label{p:exptype}
Let $G$ be an algebraic group over $\bF_q$ of exponential type. Then
\begin{enumerate}[$($i$)$]
\item the $G$-stabilizer of every $f\in\fg^*$ is connected,
\item all $G$-orbits in $\fg^*$ have even dimension.
\end{enumerate}
\end{prop}

\begin{proof}
The stabilizer of $f\in\fg^*$ equals $\Exp\fg^f$, where $\fg^f$ is the set of
all $x\in\fg$ such that $f([x,y])=0$ for all $y\in\fg$. Clearly $\fg^f$
is a linear subspace of $\fg$, so $G^f$ is connected. As $\fg^f$ is
the kernel of the alternating form $f([x,y])$, its
codimension is even.
\end{proof}

\subsection{Functional dimension}\label{ss:fun-dim}
Let $\rho$ be an irreducible representation of $G(\bF_q )$ and
$\Omega\subset\fg^*(\bF_q)$ the corresponding $G(\bF_q)$-orbit.

\begin{defin}
The {\it functional dimension\,} of  $\rho$ is
$$\operatorname{fdim}\rho = \frac{1}{2}\cdot\dim (G\Omega),$$
where $G\Omega$ is the orbit of the algebraic group $G$
containing $\Omega$.
\end{defin}

\begin{rems}
\begin{enumerate}[(i)]
\item If $G$ has exponential type then by Proposition \ref{p:exptype}(ii),
the functional dimension of any $\rho$ is an integer. But in
general, as was first observed by Lusztig \cite{lusztig}, {\it the
functional dimension may fail to be an integer\,} (see
\S\ref{ss:dim-irreps}).
\item The name ``functional dimension" is traditional in the representation
theory of real nilpotent  Lie groups. In this setting it is an integer, and
an irreducible representation of functional dimension $n$ can be realized
in the space of sections of a vector bundle on an $n$-dimensional variety.
\end{enumerate}
\end{rems}


\subsection{The fake Heisenberg groups}\label{ss:fake-heisenberg} In this section
we study the simplest example of a noncommutative connected unipotent group in
positive characteristic. We define a {\em fake Heisenberg group} over a perfect
field $k$ to be any unipotent algebraic group $G$ over $k$ of exponent
$p=\operatorname{char} k$ that can be represented as a {\em noncommutative}
central extension
\[
0 \rar{} \bG_a \rar{} G \rar{} \bG_a \rar{} 0.
\]
The existence of such groups is purely a characteristic $p$
phenomenon. Indeed, in characteristic zero, every extension of
$\bG_a$ by $\bG_a$ splits (this follows from the corresponding
statement for Lie algebras). In characteristic $p$, however, there
are plenty of examples. The name ``fake Heisenberg group'' is
motivated by the fact that in characteristic $0$ the smallest
noncommutative unipotent group is the Heisenberg group (it has
dimension $3$).

\mbr

It is a little easier to understand the corresponding Lie ring schemes. Since a
fake Heisenberg group always has nilpotence class $2$, we assume in this
section that $p>2$. Let $G$ be a fake Heisenberg group as above, and let
$\fg=\Log G$, defined via Lazard's constructions. Then, as an additive group,
$\fg$ is an extension of $\bG_a$ by $\bG_a$. This extension is split because by assumption, $G$ (and hence $\fg$) has exponent $p$, and $k$ is perfect. Thus we only need to specify the Lie bracket on $\fg$, which corresponds to a choice of an alternating bi-additive morphism
\[
B : \bG_a \times \bG_a \rar{} \bG_a.
\]
Such morphisms are exactly the ones of the form
\[
B(x,y) = \sum_{i,j} a_{ij} x^{p^i} y^{p^j},
\]
where $a_{ij}\in k$ satisfy $a_{ij}=-a_{ji}$.

\mbr

Let us now assume that $k=\bF_q$ and concentrate on the simplest example:
$B(x,y)=x^p y - x y^p$. Let $G$ denote the corresponding fake Heisenberg group.
Using the Campbell-Hausdorff formula, which is very simple for groups of class
$2$, we see that $G=\bA^2_k$ as a variety over $k$, with the group law
\[ (x,z)\cdot(y,w) = \Bigl( x+y,
z+w+\frac{1}{2}(x^p y - x y^p) \Bigr).
\]
Let us also fix, once and for all, a nontrivial additive character
$\psi:\bF_p\to\cst$. As explained in Example \ref{ex:additive} (we
remind the reader that in Section \ref{s:overview} we have used
the notation $\widehat{\fg}$ in place of $\fg^*$), this choice
allows us to identify $\fg^*$ with $\bG_a\oplus\bG_a$, using the
pairings
\[
\bigl\langle\cdot,\cdot\bigr\rangle : \bigl(\bG_a(\bF_{q^n})\oplus
\bG_a(\bF_{q^n})\bigr) \times \bigl(\bG_a(\bF_{q^n})\oplus
\bG_a(\bF_{q^n})\bigr) \rar{} \cst,
\]
\[
\bigl\langle (u,v), (y,w) \bigr\rangle =
\psi\bigl(\tr_{\bF_{q^n}/\bF_p}(uy+vw)\bigr).
\]
The adjoint action of $G$ on $\fg$ is easy to compute: given
$(x,z)\in G$, we have
\[
\Ad(x,z) : \fg\to\fg,  \qquad (y,w) \longmapsto \bigl( y, w+x^p y
- x y^p \bigr).
\]
We leave it as an exercise to the reader to compute the coadjoint
action (use Remark \ref{r:frobenius}):
\[
\Ad^*(x,z) : \fg^*\to\fg^*,  \qquad  (u,v) \longmapsto \bigl( u +
x^{1/p} v^{1/p} - x^p v, v \bigr).
\]

\mbr

This provides a concrete explanation of why the functor defining
the Serre dual $\fg^*$ of the underlying commutative unipotent
group of $\fg$ is not representable on the category of all affine
$k$-schemes. Indeed, if it were representable, then the coadjoint
action $G\times\fg^*\to\fg^*$ would be defined as a morphism of
$k$-schemes, while the presence of $x^{1/p}$ in the equation above
shows that it cannot possibly be defined unless one replaces $G$
with its \emph{perfectization} (this notion is defined in \S\ref{aa:perfect}).

\begin{rem}\label{r:topology}
On the other hand, the operation of passing from $G$ to its
perfectization changes neither the group of points over $\bF$, nor
the (Zariski or \'etale) topology of the underlying scheme. In
particular, it is irrelevant for the construction of irreducible
characters of the groups $G(\bF_{q^n})$, and also for the theory
of character sheaves, since perverse sheaves are topological
objects.
\end{rem}

We also obtain an example of a situation with nontrivial
$\bL$-packets. Namely, it is easy to check that the stabilizer in
$G$ of a point $(u,v)\in\fg^*$ is defined by the equation $x^{1/p}
v^{1/p} - x^p v=0$ (where $(x,z)$ are the coordinates on $G$), or,
equivalently, $x^{p^2}v^p=xv$. In particular, if $v\neq 0$, we see
that the stabilizer of $(u,v)$ is not connected. Moreover, it is
easy to see that the Lang isogeny may fail to be surjective for
the stabilizer in this case. Proposition \ref{p:l-indis2}(iii) now
implies that there exist different $G(\bF_q)$-orbits in
$\fg^*(\bF_q)$ that become the same orbit over $\bF$.

\mbr

Finally, observe that the nontrivial $G$-orbits in $\fg^*$ have
dimension $1$. Thus the example of fake Heisenberg groups is very
important, because it is elementary, and yet illustrates all three
phenomena mentioned at the beginning of the introduction.

\subsection{Proof of Proposition \ref{p:l-indis2}}\label{ss:proof-l-indis2}
The proof will appear in the final version of the text. (It is completely straightforward.)


\section{Character sheaves in the orbit method setting}
\label{s:char-sheaves-orbmethod}

\subsection{Construction of character sheaves} Let $k$ be an
algebraically closed field of characteristic $p>0$, and let $G$ be
a connected unipotent algebraic group over $k$ whose nilpotence
class is $\leq p-1$. We use the notation and constructions of
\S\S\ref{ss:lazard-unipotent}, \ref{ss:equivar-duality}. Thus we
have the Lie ring scheme $\fg=\Log G$, and we can choose a
representative $\fg^*\in\widetilde{\cC}_k$ of the Serre dual of
$\fg$ such that $G$ acts on $\fg^*$. Since $G$ is unipotent and
the underlying variety of $\fg^*$ is affine, it is well known that
all orbits of the $G$-action on $\fg^*$ are closed. If
$\Om\subset\fg^*$ is such an orbit and $\cL$ is an irreducible
local system on $\Om$, it follows that the $\ell$-adic complex
(see \S\ref{aa:derived-constructible}) on $\fg^*$ given by
$i_{\Om*}(\cL)[\dim\Om]$, where $i_\Om:\Om\into\fg^*$ denotes the
inclusion, is an irreducible perverse sheaf on $\fg^*$ (see
\S\ref{aa:perverse}).

\begin{defin}\label{d:charsheaves-ad-hoc}
The \emph{character sheaves} for the group $G$ are the $\ell$-adic
complexes on $G$ of the form
$\exp_*\cF^{-1}\bigl(i_{\Om*}(\cL)[\dim\Om]\bigr)$ for an
irreducible \emph{$G$-equivariant} local system $\cL$ on a
$G$-orbit $\Om\subset\fg^*$, where $\exp:\fg\rar{}G$ is the
exponential morphism and
$\cF:D^b_c(\fg,\ql)\rar{}D^b_c(\fg^*,\ql)$ is the Fourier-Deligne
transform. If $\Om$ is fixed, the corresponding character sheaves
are said to lie in the \emph{$\bL$-packet defined by $\Om$}. All
$\bL$-packets are finite (see below).
\end{defin}

\begin{rems}
\begin{enumerate}[(a)]
\item The Fourier-Deligne transform\footnote{Strictly speaking, in order to
define the Fourier-Deligne transform one has to choose an isomorphism between
$\bQ_p/\bZ_p$ and the subgroup $\mu_{p^\infty}\subset\ql^\times$ of roots of
unity whose order is a power of $p$. However, this choice is not important for
us, so throughout this paper we assume that one such choice has been made once
and for all, and we speak of ``the'' Fourier-Deligne transform.} for
$\ell$-adic complexes, which is an analogue of the classical Fourier transform
for functions, is discussed in Appendix \ref{a:fourier-deligne}. It is an
equivalence of triangulated categories which takes perverse sheaves to perverse
sheaves. Thus the character sheaves for $G$ are irreducible perverse sheaves.
\item The general definition of $G$-equivariant local systems or
$\ell$-adic complexes is recalled in \S\ref{aa:equiv}. We note
that since $G$ is connected, being equivariant is a property of an
irreducible local system, rather than extra structure. More
precisely, any two $G$-equivariant structures on an irreducible
local system $\cL$ on $\Om$ are proportional to each other.
Furthermore, the notion of irreducibility is unambiguous, in the
sense that if $\cL$ is a $G$-equivariant local system which has no
nontrivial $G$-equivariant sub-local systems, then it is also
irreducible as an ordinary local system.
\item The definition above is motivated by the classical orbit
method. Indeed, irreducible perverse sheaves on $\fg^*$ of the
form $i_{\Om*}(\cL)[\dim\Om]$ are natural analogues of the
characteristic functions of coadjoint orbits in the setting of
Section \ref{s:orbmethod}. The main difference with the classical
setting is that there may exist nontrivial irreducible
$G$-equivariant local systems on coadjoint orbits
$\Om\subset\fg^*$. More precisely, if $x\in\Om$ and
$\Pi_x=G_x/G_x^\circ$ is the group of connected components of the
stabilizer of $x$ in $G$, then the category of $G$-equivariant
local systems on $\Om$ is equivalent to the category of finite
dimensional representations of $\Pi_x$ over $\ql$. As we have
already observed in the pervious section, $\Pi_x$ may be
nontrivial. We also see that the $\bL$-packet defined by $\Om$ can
be identified with the set $\widehat{\Pi}_x$; in particular, it is
finite.
\end{enumerate}
\end{rems}

\subsection{Relation to irreducible characters} The main result of this section is
\begin{thm}
Let $G_0$ be a connected unipotent group over $\bF_q$ whose nilpotence class is
$\leq p-1$, let $G=G_0\tens_{\bF_q}\bF$, let $\Fr=\Phi_q\tens 1:G\to G$ denote
the Frobenius morphism, and let $\cs(G)$ denote the collection of character
sheaves for $G$ constructed in Definition \ref{d:charsheaves-ad-hoc}.
\begin{enumerate}[$($a$)$]
\item Each $\bL$-packet for $G$ consists of $1$ element if and only
if the stabilizer of every point of $\fg^*$ in $G$ is connected.
\item If all stabilizers of points of $\fg^*$ in $G$ are connected,
then $\cs(G)$ satisfies properties $($4-i$)$--$($4-iii$)$ in
\S\ref{ss:charsheaves-introduction}; thus the answer to question $($4$)$ is
affirmative in this case.
\item In general, $\cs(G)$ satisfies the weaker property stated at
the end of \S\ref{ss:charsheaves-introduction}.
\end{enumerate}
\end{thm}

The proof will appear in the final version of the text. It is not difficult to deduce this theorem from the orbit method for finite nilpotent groups explained in Section \ref{s:orbmethod} and the results of Appendix \ref{a:characters}.


\section{Character sheaves for general unipotent groups}\label{s:charsheaves-general}

In this section we formulate a general definition of character
sheaves for a unipotent group $G$ over an algebraically closed
field $k$ of characteristic $p>0$, and explain why it agrees with
the \emph{ad hoc} definition given in the previous section when
$G$ is connected and has nilpotence class $<p$. We will use some
results on idempotents in monoidal categories; for their proofs
and many more details we refer the reader to \cite{idempotents}.

\subsection{Idempotents in monoidal categories}\label{ss:idempotents}
If $\cM=(\cM,\tens,\e)$ is a monoidal category, where
$\tens:\cM\times\cM\rar{}\cM$ is the monoidal bifunctor and $\e$ is the unit
object, we define an arrow $\e\to e$ in $\cM$ to be an {\em idempotent arrow}
if it becomes an isomorphism after tensoring with $e$ either on the left or on
the right. An object $e$ of $\cM$ will be called an {\em idempotent} if there
exists an idempotent arrow $\e\to e$. This notion is rather rigid: for example,
if $\pi_1,\pi_2:\e\to e$ are two idempotent arrows, then there exists a unique
morphism $f:e\to e$ satisfying $f\circ\pi_1=\pi_2$. By symmetry, $f$ is then
necessarily an isomorphism.
\begin{rem}\label{r:idempotents}
For the sake of brevity, our terminology here differs from that of
\cite{idempotents}. Namely, what we call an idempotent here is
referred to as a \emph{closed idempotent} in \emph{op.~cit.}. The
motivation for the terminology is also explained in
\emph{op.~cit.}
\end{rem}
We define a partial order on the set of isomorphism classes of
idempotents in $\cM$ as follows: $e\leq e'$ whenever $e\tens
e'\cong e$. One can show that this condition is equivalent to
$e'\tens e\cong e$; this order relation can also be characterized
in many other ways, see \cite{idempotents}.

\subsection{Important remark}\label{ss:important-remark}
The notion of an idempotent in fact depends only on the bifunctor
$\tens$ and not on any other ingredients of a monoidal category,
i.e., the associativity constraint and the unit object. Indeed,
the associativity constraint does not appear in the definition,
and since any two unit objects are isomorphic, it is easy to see
that $e\in\cM$ is an idempotent with respect to one unit object if
and only if it is an idempotent with respect to every unit object.
Similarly, the partial order relation on the set of isomorphism
classes of idempotents in $\cM$ introduced above also depends only
on $\tens$, which is obvious from the definition.

\subsection{Hecke subcategories}\label{ss:hecke}
Let $\cM$ be a monoidal category. If $e$ is an idempotent in
$\cM$, we define $e\cM$ to be the full subcategory of
$\cM$ consisting of objects isomorphic to those of the form
$e\tens X$, where $X\in\cM$. It is easy to see that $e\cM$
consists precisely of those $Y\in\cM$ for which $e\tens Y\cong Y$.
The full subcategories $\cM e$ and $e\cM e$ of $\cM$ are
defined similarly. We call $e\cM e$ the {\em Hecke subcategory}
associated to the idempotent $e$. It is closed under $\tens$,
which makes it a monoidal category with unit object $e$. If $\cM$
is a braided monoidal category (BMC) or, more generally, if
$X\tens Y\cong Y\tens X$ for any pair of objects $X$, $Y$ of $\cM$
(we do not even need these isomorphisms to be functorial), it is
clear that the three subcategories $e\cM$, $\cM e$ and $e\cM e$
all coincide.

\subsection{Minimal idempotents}\label{ss:minimal-idempotents}
Let $\cM$ be an additive monoidal category (this means that $\cM$
is an additive category equipped with a monoidal structure such
that the bifunctor $\tens$ is biadditive). We define a {\em
minimal idempotent} in $\cM$ to be an idempotent which is minimal
in the set of nonzero idempotents in $\cM$ with respect to the
partial order relation introduced in \S\ref{ss:idempotents}. If
$\cM$ is a braided monoidal category and $e$, $e'$ are idempotents
in $\cM$, one can prove \cite{idempotents} that $e\tens e'$ is
again an idempotent. This has the consequence that if $\cM$ is an
additive BMC and $e_1$, $e_2$ are minimal non-isomorphic
idempotents in $\cM$, then $e_1\tens e_2=0$. Consequently,
$e_1\cM\cap e_2\cM=0=(e_1\cM)\tens(e_2\cM)$.

\subsection{The categories $\sD(G)$ and $\sD_G(G)$}\label{ss:derived-categories-of-G}
Let $k$ be an algebraically closed field of characteristic $p>0$,
and $G$ a unipotent algebraic group over $k$. We also fix a prime
$\ell\neq p$. The category $\sD(G)=D^b_c(G,\ql)$ of
\emph{constructible $\ell$-adic complexes on $G$} is discussed in
\S\ref{aa:derived-constructible}. It is a monoidal category with
respect to the \emph{convolution} operation defined below (see also
\S\ref{aa:convolution}). The category $\sD_G(G)$ of $G$-equivariant
objects in $\sD(G)$, where $G$ acts on itself by conjugation, is
defined in \S\ref{aa:categories-equivariant}. It is also a
monoidal category with respect to convolution, and, moreover, it
is naturally braided (see below). We have the forgetful functor $\sD_G(G)\rar{}\sD(G)$ by means of which every object of $\sD_G(G)$ can be thought of as an $\ell$-adic complex on $G$. If $G$ is connected, this functor is fully faithful (see \S\ref{aa:ber-lun}).

\mbr

The convolution bifunctor on $\sD(G)$ and $\sD_G(G)$ is defined by
\[
M*N = \mu_! (M\boxtimes N),
\]
where $\mu:G\times G\to G$ is the multiplication morphism. The unit
object  $\e$ in $\sD(G)$ or in $\sD_G(G)$ equals the delta-sheaf
$1_*\ql=1_!\ql$, where $1:\Spec k\to G$ is the unit of~$G$.

\mbr

The monoidal category $\sD_G(G)$ has a natural structure of braided
category. The braiding $\be_{M,N}:M*N\isom N*M$ is defined as
follows. Consider the commutative diagram
\[
\xymatrix{
  G\times G \ar[d]_{\tau} \ar[r]^{\xi} &  G\times G \ar[d]^{\mu} \\
  G\times G \ar[r]^{\mu} & G,
   }
\]
where $\tau(g,h):=(h,g)$ and $\xi(g,h):=(g,g^{-1}hg)$. We have
$M*N=\mu_!(M\boxtimes N)$, and the above diagram shows that
$N*M=(\mu\tau )_!(M\boxtimes N)=\mu_!\xi_!(M\boxtimes N)$. We define
$\be_{M,N}:\mu_!(M\boxtimes N)\isom\mu_!\xi_!(M\boxtimes N)$ by
$\be_{M,N}:=\mu_!(f)$, where $f:M\boxtimes N\isom\xi_!(M\boxtimes N)$
comes from the $G$-equivariant structure on $N$.

\mbr

Another important structure on $\sD_G(G)$ is the canonical automorphism of the identity functor, defined as follows. Let $p_2:G\times G\to G$ denote the second projection, $c:G\times G\to G$ the conjugation action morphism $c(g,h)=ghg^{-1}$, and $\De:G\to G\times G$ the diagonal. Then $c\circ\De=\id_G=p_2\circ\De$. For each $M\in\sD_G(G)$, the $G$-equivariant structure on $M$ yields an isomorphism $c^*M\to p_2^*M$; pulling it back by $\De$, we obtain an isomorphism $\te_M:M\isom M$, called the \emph{canonical automorphism of $M$}. It is clear that the collection of all the $\te_M$'s defines an automorphism of the identity functor $\Id:\sD_G(G)\rar{}\sD_G(G)$. It is related to the braiding via the easily verified formula
\[
\te_{M*N} = \be_{N,M}\circ\be_{M,N}\circ (\te_M*\te_N) \qquad \forall\, M,N,\in\sD_G(G).
\]

\subsection{General definition of character sheaves}\label{ss:charsheaves-general}
For every minimal idempotent $e\in\sD_G(G)$, we let
$\cM_e^{perv}\subset e\sD_G(G)$ denote the full
subcategory consisting of those complexes in $e\sD_G(G)$ that are
perverse sheaves on $G$. Let us also recall that
$e\sD_G(G)=\sD_G(G)e=e\sD_G(G)e$, because $\sD_G(G)$ is braided.

\begin{defin}\label{d:charsheaves-general}
A {\em character sheaf} on $G$ is an indecomposable object of the
category $\cM_e^{perv}$ for some minimal idempotent
$e\in\sD_G(G)$. The collection of isomorphism classes of all
character sheaves on $G$ will be denoted by $\cs(G)$. If
$e\in\sD_G(G)$ is a minimal idempotent, the set of isomorphism
classes of character sheaves in $e\sD_G(G)$ will be called the
{\em $\bL$-packet} associated to $e$. An $\bL$-packet is said to be
{\em trivial} if it consist of only one element. Character sheaves
that belong to the same $\bL$-packet are said to be {\em
$\bL$-indistinguishable}.
\end{defin}

\begin{rems}
\begin{enumerate}[(1)]
\item The subcategory $\cM_e^{perv}$ of $\sD_G(G)$ is
closed under taking direct summands. Thus the notion of an
indecomposable object of $\cM_e^{perv}$ is unambiguous.
\item A character sheaf determines the corresponding minimal
idempotent $e$ uniquely, because of the remark in
\S\ref{ss:minimal-idempotents}. Thus different $\bL$-packets are
disjoint.
\end{enumerate}
\end{rems}

\subsection{Conjectural properties}\label{ss:conjectures} We now
state several conjectures about character sheaves which we expect
to hold for general unipotent groups (regardless of their
nilpotence class). We also introduce the notion of functional
dimension for a character sheaf that agrees with the geometric
notion of functional dimension whenever the orbit method is
applicable.
\begin{conj}\label{conj:1}
\begin{enumerate}[$($a$)$]
\item For every minimal idempotent $e\in\sD_G(G)$, there exists
a number $n_e\in\{0,1,\dotsc,\dim G\}$ such that $\bD^-_G e\cong e[-2n_e]$,
where we set $\bD^-_G=\bD_G\circ\iota^*=\iota^*\circ\bD_G$, $\iota:G\rar{}G$
being the inversion map.
\item The complex $e[-n_e]$ $($where we forget the $G$-equivariant
structure on $e${}$)$ is a perverse sheaf on $G$ $($see
\S\ref{aa:perverse} for the definition of perverse sheaves$)$.
\end{enumerate}
\end{conj}

\begin{defin}\label{d:fun-dim}
Assuming that Conjecture \ref{conj:1} holds, the number
$d_e=\frac{\dim G-n_e}{2}$ is called the \emph{functional
dimension} of $e$. Note that it may fail to be an integer. We also
call $d_e$ the functional dimension of every character sheaf in
the $\bL$-packet defined by $e$.
\end{defin}

\begin{conj}\label{conj:2}
\begin{enumerate}[$($a$)$]
\item Character sheaves are irreducible perverse sheaves.
\item The triangulated subcategory $e\sD_G(G)\subset\sD_G(G)$ is
generated by $\cM_e^{perv}$.
\item We have $\Ext^i_{\sD_G(G)}(M_1,M_2)=0$ for $i>0$ and any two
character sheaves $M_1,M_2$.
\item All $\bL$-packets for $G$ are finite.
\end{enumerate}
\end{conj}

\begin{cor}
Assuming that Conjecture \ref{conj:2} holds, $\cM_e^{perv}$ is a
semisimple abelian category with finitely many isomorphism classes
of irreducible objects, and its bounded derived category is
naturally equivalent to $e\sD_G(G)$.
\end{cor}

\begin{conj}\label{conj:3}
\begin{enumerate}[$($a$)$]
\item The subcategory $\cM_e:=\cM_e^{perv}[n_e]\subset e\sD_G(G)$
is closed under convolution, and is thus a monoidal category with
unit object $e$ $($if Conjecture \ref{conj:1}$($b$)$ holds$)$.
\item The monoidal category $\cM_e$ is rigid, and
$M^*\cong(\bD_G^-M)[2n_e]$ for $M\in\cM_e$.
\item $\cM_e$ is a modular category, where the ribbon structure, or, equivalently, the ``twist'', is given by the restriction of the canonical automorphism $\te$ to $\cM_e$.
\end{enumerate}
\end{conj}

\begin{conj}\label{conj:4}
If an $\bL$-packet for $G$ is trivial, the functional dimension of
the corresponding minimal idempotent $e\in\sD_G(G)$ is an integer.
\end{conj}

\begin{rem}
If Conjecture \ref{conj:1}(b) holds, and $e\in\sD_G(G)$ is a
minimal idempotent corresponding to a trivial $\bL$-packet, then
this $\bL$-packet consists of the complex $e[-n_e]$.
\end{rem}

\begin{defin}\label{d:easy-groups}
An algebraic group $H$ over $k$ is \emph{easy} if each $h\in H(k)$
lies in the neutral connected component of its centralizer.
\end{defin}

Note that this definition makes sense for \emph{any} algebraic group over \emph{any} field. In particular, it makes sense for reductive groups.

\begin{conj}\label{conj:5}
The group $G$ as above is easy if and only if all $\bL$-packets for
$G$ are trivial.
\end{conj}

The last two conjectures are concerned with the relationship
between character sheaves and irreducible characters. Thus we now
assume that our ground field is $k=\bF$, an algebraic closure of a
field with $p$ elements, and that $G=G_0\otimes_{\bF_q}\bF$ for a
unipotent group $G_0$ over $\bF_q$. We have the corresponding
Frobenius morphism $\Fr:G\rar{}G$ (see
\S\ref{ss:charsheaves-introduction}).

\begin{conj}\label{conj:6}
Suppose $G$ is connected and all $\bL$-packets for $G$ are trivial. Let
$e\in\sD_G(G)$ be a minimal idempotent such that $\Fr^*e\cong e$, let
$\psi_e:\Fr^*e\rar{\simeq}e$ be the unique isomorphism such that
$\psi_e\circ\Fr^*(\pi)=\pi$ for some, and hence for every, idempotent arrow
$\pi:\e\rar{}e$, and let $t_e:G_0(\bF_q)\rar{}\ql$ denote the corresponding
trace function, as in \S\ref{ss:charsheaves-introduction}. Then the function
$\chi_e=q^{d_e+n_e}t_e$ is an irreducible character of $G_0(\bF_q)$, and every
irreducible character of $G_0(\bF_q)$ has this form. Moreover,
$\chi_e(1)=q^{d_e}$.
\end{conj}

\begin{rems}
\begin{enumerate}[(a)]
\item We remind the reader that $\e$ denotes the delta-sheaf at
the identity element of $G$, which is a unit object in the
monoidal category $\sD_G(G)$.
\item One can prove the existence of an isomorphism $\psi_e$ with
the property stated above using the results of \cite{idempotents};
the details will appear elsewhere.
\item In the situation of Conjecture \ref{conj:6}, we obtain an
affirmative answer to question (4) in
\S\ref{ss:charsheaves-introduction}, at least up to cohomological
shift, as in \S\ref{ss:charsheaves-commutative}. Namely, for every
minimal idempotent $e\in\sD_G(G)$ such that $(\Fr^n)^*e\cong e$,
let $\psi_{n,e}:(\Fr^n)^*e\rar{\simeq}e$ denote the isomorphism as
in the statement of the conjecture, let
$\psi'_{n,e}=q^{n(d_e+n_e)}\cdot\psi_{n,e}$, and let
$t'_{n,e}:G_0(\bF_{q^n})\rar{}\ql$ denote the trace function
corresponding to $\psi'_{n,e}$. It is easy to check that
$\psi'_{k,e}$ is induced by $\psi'_{n,e}$ whenever $n\lvert k$,
and, moreover, it follows from Conjecture \ref{conj:6} that the
functions $t'_{n,e}$ are precisely the irreducible characters of
$G_0(\bF_{q^n})$.
\end{enumerate}
\end{rems}

\begin{cor}\label{c:dims-irreps}
In the situation above, assume that $G$ is easy. If Conjectures
\ref{conj:4}, \ref{conj:5} and \ref{conj:6} hold, then the
dimension of every irreducible representation of $G_0(\bF_q)$ is a
power of $q$.
\end{cor}

This statement is obvious. However, as we have already shown in
\S\ref{ss:dim-irreps}, the answer to question (4) in
\S\ref{ss:charsheaves-introduction} is negative in general.
However, we still expect that the following weaker statement is
true.

\begin{conj}\label{conj:7}
Let $G_0$ be any connected unipotent group over $\bF_q$, and let
$G=G_0\tens_{\bF_q}\bF$ and $\Fr:G\rar{}G$ be as above. Let
$\cs(G)^{\Fr}$ denote the set of those $M\in\cs(G)$ for which
$\Fr^*M\cong M$, and for each such $M$, choose an isomorphism
$\psi_M:\Fr^*M\rar{\simeq}M$, and let $t_M:G_0(\bF_q)\rar{}\ql$
denote the corresponding trace function. Then the functions
$\bigl\{t_M\bigr\}_{M\in\cs(G)^{\Fr}}$ form a basis for the space
of class functions on $G_0(\bF_q)$.
\end{conj}


\subsection{Character sheaves and the orbit method}\label{ss:charsheaves-orbmethod}
In this subsection we prove

\begin{thm}\label{t:charsheaves-orbmethod}
Let $G$ be a connected unipotent group over an algebraically
closed field $k$ of characteristic $p>0$, whose nilpotence class
is $\leq p-1$.
\begin{enumerate}[$($a$)$]
\item There is a natural bijection between $($isomorphism classes of$)$ minimal
idempotents in $\sD_G(G)$ and $G$-orbits in $\fg^*$. If $e\in\sD_G(G)$
corresponds to $\Om\subset\fg^*$ under this bijection, then $n_e$ is the
codimension of $\Om$ in $\fg^*$ and $d_e=\frac{1}{2}\dim\Om$.
\item The notions of character sheaves and $\bL$-packets for $G$
provided by Definitions \ref{d:charsheaves-ad-hoc} and
\ref{d:charsheaves-general} coincide.
\item Conjectures 1, 2, 3$($a$)$ and 3$($b$)$ hold in this situation.
\end{enumerate}
\end{thm}

The proof of this theorem is based on the following
\begin{prop}\label{p:convol}
Every choice of Lie polynomials $\phi(X,Y)$ and $\psi(X,Y)$ satisfying the
conditions of Lemma \ref{l:phi-psi} with $c=p-1$ determines a collection of
isomorphisms
\begin{equation}\label{e:isoms}
\exp^*(M*N)\cong\exp^*(M)*\exp^*(N) \qquad\forall\,M,N\in\sD_G(G)
\end{equation}
that are bifunctorial with respect to $M$ and $N$.
\end{prop}
\begin{proof}
We first reformulate the result entirely in terms of the category $\sD_G(\fg)$.
Let us write $p_1,p_2:\fg\times\fg\rar{}\fg$ for the two projections,
$\al:\fg\times\fg\rar{}\fg$ for the addition morphism, and
$CH:\fg\times\fg\to\fg$ for the morphism defined by
$(x,y)\mapsto\sum_{i=1}^{p-1}CH_i(x,y)$, where $CH(X,Y)$ is the
Campbell-Hausdorff series. Then $G$ is exactly $\fg$ equipped with the
operation $CH$. Thus we now obtain two bifunctors on $\sD_G(\fg)$, defined by
\[
M* N = \al_!(p_1^*M\tens p_2^* N) \qquad\text{and}\qquad M\bigstar N =
CH_!(p_1^*M\tens p_2^* N).
\]
The first is the usual convolution of complexes on $\fg$, the second one
corresponds to the convolution of complexes on $G$. We must prove that every
choice of $\phi$ and $\psi$ as in Lemma \ref{l:phi-psi} yields an isomorphism
of bifunctors between $*$ and $\bigstar$.

\begin{lem}\label{l:h-isomorphism}
The morphism
\[
h : \fg\times\fg\rar{}\fg\times\fg, \qquad (x,y)\longmapsto\bigl(
e^{\ad\phi(x,y)}(x), e^{\ad\psi(x,y)}(y) \bigr),
\]
is an isomorphism of schemes.
\end{lem}
\begin{proof}
It suffices to check that in the situation of Lemma \ref{l:phi-psi}, the map
$h:\fg'\times\fg'\rar{}\fg'\times\fg'$ defined by the formula above is
bijective for every $\fg'\in\nilp_c$. The functor $\nilp_c\rar{}\cS{}ets$
defined by $\fg'\mapsto\fg'\times\fg'$ is (co)representable by the free
nilpotent algebra $\fF$ over $\bZ\bigl[\frac{1}{c!}\bigr]$ of class $c$ with
generators $x,y$, i.e., $\Hom(\fF,\fg')=\fg'\times\fg'$ for any
$\fg'\in\nilp_c$. The map $h$ comes from an endomorphism of $\fF$. This
endomorphism induces the identity on $\fF/[\fF,\fF]$, so it is invertible (this
is a Lie algebra version of the inverse function theorem, with $\fF$ playing
the role of the maximal ideal in the ring of commutative formal power series
and $\fF/[\fF,\fF]$ playing the role of $\fm/\fm^2=\text{cotangent space}$).
\end{proof}

Continuing with the proof of Proposition \ref{p:convol}, note that we have
$CH=\al\circ h$ by construction, so in view of Lemma \ref{l:h-isomorphism}, it
suffices to prove that there are isomorphisms $h^*(p_1^*M\tens p_2^* N)\cong
p_1^*M\tens p_2^* N$, bifunctorial with respect to $M, N\in \sD_G(\fg)$. We
will show that, in fact, there are functorial isomorphisms $h^*p_1^*M\cong
p_1^*M$ and $h^*p_2^* N\cong p_2^* N$.

\mbr

Observe that there is a commutative diagram
\[
\xymatrix{
\fg\times\fg \ar[r]^f \ar[rd]_{p_1\circ h}   &   \fg\times\fg \ar[d]^{\ga}    \\
  &  \fg
   }
\]
where $f$ is given by
\[
f(x,y) = \bigl( x, \phi(x,y) \bigr),
\]
and where $\ga$ is the conjugation action map\footnote{Note that it is the
second copy of $\fg$ that acts on the first one, and not the other way
around.}, given by
\[
\ga(s,t) = e^{\ad t}(s).
\]
Note also that $p_1\circ f=p_1$. Now, by the definition of equivariance, we
have functorial isomorphisms $\ga^*M\cong p_1^*M$ for all $M\in \sD_G(\fg)$.
Thus we obtain a sequence of functorial isomorphisms:
\[
h^*p_1^*M \cong (p_1\circ h)^*M = (\ga\circ f)^*M \cong f^*\ga^*M \cong f^*
p_1^*M \cong (p_1\circ f)^*M = p_1^*M
\]
for every $M\in \sD_G(\fg)$. An analogous construction yields functorial
isomorphisms $h^* p_2^* N\cong p_2^* N$ for all $N\in \sD_G(\fg)$. This proves
the proposition.
\end{proof}

We warn the reader that the result we have proved \emph{does not} imply that
$\sD_G(G)$ and $\sD_G(\fg)$ are equivalent as monoidal categories. In fact, in
general there does not exist a way of choosing the isomorphisms \eqref{e:isoms}
so that they would be compatible with the natural associativity constraints on
$\sD_G(G)$ and $\sD_G(\fg)$.

\begin{proof}[Proof of Theorem \ref{t:charsheaves-orbmethod}] A detailed proof will appear in the final version of the paper. Let us indicate the main ideas, from which the full argument can be easily recovered.

\mbr

In view of
Proposition \ref{p:convol} and the remarks in \S\ref{ss:important-remark}, the
functor $\exp_*:\sD_G(\fg)\rar{}\sD_G(G)$ induces a bijection between the set
of (isomorphism classes of) minimal idempotents in $\sD_G(\fg)$ and that in
$\sD_G(G)$. On the other hand, by Proposition \ref{p:fourier-equivariant}, the
Fourier-Deligne transform induces an equivalence of monoidal categories
$\cF:\sD_G(\fg)\rar{\sim}\sD_G(\fg^*)$, where $\sD_G(\fg^*)$ is equipped with
the monoidal structure $M\odot N=(M\tens N)[-\dim G]$. Hence minimal
idempotents in $\sD_G(G)$ correspond to minimal idempotents in $\sD_G(\fg^*)$
with respect to $\odot$. Let $i_\Om:\Om\into\fg^*$ be the inclusion of a
$G$-orbit, let $(\ql)_\Om$ denote the constant sheaf on $\Om$ with stalk $\ql$,
and let $(\ql)_{\fg^*}$ denote the constant sheaf on $\fg^*$ with stalk $\ql$.
Put $d=\dim G$. Then $(\ql)_{\fg^*}[d]$ is a unit object in $\sD_G(\fg^*)$ with
respect to $\odot$. Let $e_\Om=(i_\Om)_* (\ql)_\Om[d]$. By adjunction, we have a natural morphism $(\ql)_{\fg^*}[d]\to e_\Om$, which is clearly an idempotent arrow. The Hecke subcategory of $\sD_G(\fg^*)$ defined by $e_\Om$ identifies with $\sD_G(\Om)$. Moreover, it is easy to check that $e_\Om$ is a minimal idempotent of $\sD_G(\fg^*)$, and every minimal idempotent of $\sD_G(\fg^*)$ has this form. The rest is rather computational.
\end{proof}


\section{Unipotent groups arising from associative algebras}\label{s:algebra-groups}

In this section we introduce a large collection of examples that
are quite different from the ones we have considered so far, in
the sense that the groups we will discuss may have ``large''
nilpotence class (relative to the characteristic of the ground
field), and therefore the orbit method cannot be used to study
their representations.

\subsection{Unipotent linear groups}\label{ss:gutkin-main}
Let $n\in\bN$. We define a group scheme $UL_n$ over $\Spec\bZ$ as
follows. For every commutative ring $R$, we let $UL_n(R)$ be the
group of unipotent upper-triangular matrices over $R$ of size $n$.
It is clear that the functor $\Spec R\longmapsto UL_n(R)$ is
representable by an affine group scheme whose underlying scheme is
the affine space $\bA_{\bZ}^{n(n-1)/2}$ of dimension $n(n-1)/2$
over $\Spec\bZ$. Following A.A.~Kirillov, we call $UL_n$ the
\emph{unipotent linear group} (of size $n$). Note that the
nilpotence class of $UL_n$ is equal to $n-1$.

\mbr

In particular, if $q$ is a power of a prime $p$, we have the finite group
$UL(n,q):=UL_n(\bF_q)$ of order $q^{n(n-1)/2}$. Complex irreducible
representations of $UL(n,q)$ have been studied by several authors: see, for
example, \cite{kazhdan,kir-triangular,andre-1,ning-yan,andre-3}. There are some
interesting results in this theory, but also many unanswered questions. The
main difficulty lies in the cases where $n$ is large compared to $p$; for
example, if $n>p$, then the orbit method explained in Section \ref{s:orbmethod}
cannot be applied to the group $UL(n,q)$. One of the first nontrivial facts
about representations of $UL(n,q)$ is the following
\begin{thm}[Isaacs]\label{t:isaacs} The dimension of every complex irreducible
representation of $UL(n,q)$ is a power of $q$.
\end{thm}

This result was conjectured by J.~Thompson and proved by I.M.~Isaacs in
\cite{isaacs}\footnote{We thank Jon Alperin and George Glauberman for providing
us with this reference.} (for $n\leq p$, it was proved in \cite{kazhdan}).
E.~Gutkin claimed to have proved this fact in \cite{gutkin}; in fact, he stated
a much stronger result (see Theorem \ref{t:gutkin} below). However, his proof
has a gap, and to the best of our knowledge Isaacs' proof is the first correct
one.

\subsection{Unipotent algebra groups: character theory}\label{ss:unip-alg-grps}
We now present a wide generalization of the unipotent linear
groups.
\begin{defin}\label{d:algebra-groups}
Let $k$ be a commutative ring, and let $A$ be an associative
$k$-algebra in which every element is nilpotent (in particular,
$A$ is then non-unital). We define $A^\times$ to be the group that
has $A$ as the underlying set and $x\circ y=x+y+xy$ as the group
operation. $($Informally speaking, $A^\times=1+A.)$ Note that the
assumption that every element of $A$ is nilpotent is used to
ensure that elements of $A^\times$ have inverses. A group of the
form $A^\times$ is called an \emph{algebra group} over $k$.
\end{defin}

\begin{defin}\label{d:unipotent-algebra-groups}
In the situation of Definition \ref{d:algebra-groups}, assume that
$A$ is free of finite rank as a $k$-module. For every commutative
unital $k$-algebra $R$, we have the group $(A\otimes_k R)^\times$.
It is clear that the functor $\Spec R\longmapsto(A\otimes_k
R)^\times$ is representable by a unipotent group scheme over $k$.
By abuse of notation, we denote it by $A^\times$, and call it the
\emph{unipotent algebra group} defined by $A$.
\end{defin}

\begin{example}\label{ex:ul-n-q}
Let $A$ be the algebra of $n\times n$ strictly upper triangular
matrices over $\bZ$ (with respect to matrix multiplication). Then
$A^\times\cong UL_n$ as group schemes over $\bZ$.
\end{example}

\begin{rem}
Note that with the standard definition of an algebra group
\cite{isaacs,halasi} one starts with a finite dimensional
associative \emph{unital} algebra $A$ over a field $k$, and forms
the group $1+J$, which is a subgroup of the group of units of $A$,
where $J$ is the Jacobson radical of $A$. However, this is a
special case of Definition \ref{d:algebra-groups}, because $J$ is
an associative algebra over $k$ in its own right, and every
element of $J$ is nilpotent. Conversely, if $A$ is as in
Definition \ref{d:algebra-groups}, then we can form an algebra
$\widetilde{A}=k\cdot 1\oplus A$ over $k$ by formally adjoining
$1$ to $A$, and it is clear that $A$ is then the Jacobson radical
of $\widetilde{A}$.
\end{rem}

\begin{thm}[Gutkin-Halasi]\label{t:gutkin}
Let $A$ be a finite dimensional associative nilpotent algebra over
$\bF_q$. Then every complex irreducible representation of
$A^\times$ is induced from a $1$-dimensional representation of a
subgroup of the form $B^\times$, where $B\subseteq A$ is an
associative subalgebra.
\end{thm}

\begin{rems}
\begin{enumerate}[(1)]
\item In the situation of Theorem \ref{t:gutkin}, Isaacs proved
(\cite{isaacs}, Theorem A) that the dimension of every irreducible
representation of $A^\times$ is a power of $q$. In view of Example
\ref{ex:ul-n-q}, Thompson's conjecture (Theorem \ref{t:isaacs}) is
a special case of his result.
\item In turn, it is clear that Isaacs' result would follow
immediately from Corollary \ref{c:dims-irreps} (which depends on
Conjectures \ref{conj:4}, \ref{conj:5} and \ref{conj:6} in
\S\ref{ss:conjectures}) applied to the group
$G=A^\times\tens_{\bF_q}\bF$.
\item On the other hand, it is also clear that Theorem \ref{t:gutkin}
is substantially stronger than Isaacs' result. The theorem was
stated by Gutkin in \cite{gutkin} and proved by Halasi in
\cite{halasi}. However, Halasi's proof itself relies on Isaacs'
result; in other words, one cannot use \cite{halasi} to reprove
Theorem A of \cite{isaacs}.
\item The situation was somewhat improved in \cite{base-change}
where a direct proof of Theorem \ref{t:gutkin} which is based on
Halasi's methods but avoids using \cite{isaacs} was given.
\item In the special case where $A^p=(0)$, Theorem \ref{t:gutkin}
was proved by C.A.M.~Andr\'e \cite{andre}. In fact, Andr\'e
pointed out that for every functional $f\in
A^*=\Hom_{\bF_q}(A,\bF_q)$, there exists a polarization of $A$ at
$f$ which is multiplicatively closed. Here we consider $A$ as a
Lie algebra in the usual way: $[a,b]=ab-ba$. (The proof of this
claim is given in Theorem \ref{t:vergne2}.) In the case where
$A^p=(0)$, the group $A^\times$ is of exponential type, namely,
the exponential map $\exp:A\to A^\times$ can be defined by the
usual series, $\exp(x)=x+x^2/2+x^3/6+\dotsb$. If $B\subseteq A$ is
a multiplicatively closed polarization at $f$, it follows that
$B^\times=\exp(B)$, which together with the usual orbit method
(cf.~\cite{kirillov,orbmethod}) implies Theorem \ref{t:gutkin}.
However, when $A^p\neq (0)$, the proof of Gutkin's claim is
substantially more difficult.
\end{enumerate}
\end{rems}

\mbr

To end this discussion, we would like to mention that many
interesting results about characters of the unipotent linear
groups obtained in \cite{ning-yan} and \cite{andre-1,andre-3} have
recently been generalized to algebra groups over finite fields in
\cite{diaconis-isaacs}. Unfortunately, this theory lies beyond the
scope of our article.

\subsection{Unipotent algebra groups: geometric aspects}
\label{ss:unip-alg-grps-geom} We briefly consider unipotent
algebra groups from the point of view of the questions posed in
Section \ref{s:overview}. Let $A$ be a finite dimensional
associative nilpotent algebra over $\bF_q$, and $A^\times$ the
corresponding unipotent group over $\bF_q$.
\begin{thm}[see \cite{base-change}]
There exist injective $\Fr_q$-equivariant maps
\[
T_m^n : \widehat{A^\times(\bF_{q^m})} \rar{}
\bigl(\widehat{A^\times(\bF_{q^n})}\bigr)^{\Fr_q^m}
\]
for all pairs of positive integers $m\lvert n$, which satisfy
$T_m^k=T_n^k\circ T_m^n$ whenever $m\lvert n\lvert k$ and commute
with all automorphisms induced by algebra automorphisms of $A$
over $\bF_q$.
\end{thm}
Thus the answer to question (1) is positive. One can also prove
\cite{base-change} that the base change maps are surjective in some special
cases, but in general question (2) remains open, as do questions (3) and (4).
(However, if Conjectures \ref{conj:4}, \ref{conj:5} and \ref{conj:6} in
\S\ref{ss:conjectures} hold, then the answer to question (4) is positive for
unipotent algebra groups.)

\subsection{Other examples} We have seen one way of generalizing
the unipotent linear groups to produce interesting examples of
unipotent algebraic groups. Another way of producing infinite
families of examples is to replace the groups $UL_n$ with maximal
unipotent subgroups of classical simple groups other than $SL_n$.
It was already observed in \cite{isaacs} (see also
\cite{previtali}) that these groups exhibit behavior different
from that of the algebra groups. Lusztig pointed out
\cite{lusztig} that these groups also illustrate some nontrivial
aspects of the theory of character sheaves for unipotent groups.
In the remainder of this section we will concentrate on maximal
unipotent subgroups of the symplectic groups $Sp_{2n}(k)$ (i.e.,
simple groups of type $C$) and their generalizations, and in the
next section we will present Lusztig's example, which is the case
$n=2$ and $k=\bF_{2^s}$.

\subsection{Unipotent symplectic groups}\label{ss:symplectic}
These groups are defined as maximal unipotent subgroups
$USp(2n,q)=USp_{2n}(\bF_q)$ of the symplectic groups
$Sp(2n,q)=Sp_{2n}(\bF_q)$. This is a special case of the following
\begin{defin}\label{d:unip-symp}
Let $A$ be an associative algebra over a commutative ring $k$ in which every
element is nilpotent, and let $\sg:A\rar{} A$ be an {\em anti-involution},
i.e., $\sg$ is $k$-linear and satisfies $\sg(xy)=\sg(y)\sg(x)$ for all $x,y\in
A$, and $\sg^2=1$. It is clear that $\sg$ induces an anti-automorphism of the
group $A^\times$. The {\em generalized unipotent symplectic group} associated
to this data is the subgroup
\[
Sp(A,\sg)=\bigl\{x\in A^\times \st x+\sg(x)+x\sg(x)=0 \bigr\}.
\]
In other words, since we (formally) have $1+x\circ y=1+x+y+xy=(1+x)(1+y)$, so
that $A^\times$ should be thought of as the group $1+A$, we should likewise
think of $Sp(A,\sg)$ as the group of elements $g\in 1+A$ such that
$\sg(g)=g^{-1}$.
\end{defin}

\begin{rem}
Let $V$ be a finite dimensional vector space over a field $k$,
equipped with a symplectic form $\om$, and construct a complete
flag $0=V_0\subset V_1\subset\dotsb\subset V_{2n}=V$ of subspaces
of $V$ in the following way. Pick a Lagrangian subspace $L\subset
V$ with respect to $\om$, let $0=V_0\subset
V_1\subset\dotsb\subset V_n=L$ be an arbitrary complete flag of
subspaces of $L$, and define $V_{j}$ to be the orthogonal
complement of $V_{2n-j}$ in $V$ with respect to $\om$ for all
$n<j\leq 2n$. If $A$ is the algebra of endomorphisms $\phi$ of $V$
satisfying $\phi(V_i)\subseteq V_{i-1}$ for $1\leq i\leq 2n$, then
the equation $\om(\sg(\phi)(x),y)=\om(x,\phi(y))$ defines an
involution $\sg$ of $A$, and it is easy to check that $Sp(A,\sg)$
is a maximal unipotent subgroup of $Sp(V)$. In particular,
$USp(2n,q)$ is indeed a special case of Definition
\ref{d:unip-symp}.
\end{rem}

When $p=\operatorname{char}(\bF_q)=2$, the groups $USp(2n,q)$ behave
differently from $UL(2n,q)$; in particular, in Section \ref{s:usp4} we explain,
following \cite{lusztig}, the classification of irreducible representations of
the group $USp(4,2^s)$ ($s\in\bN$), and show that it has irreducible
representations of dimension $2^{s-1}$ if $s\geq 2$. However, when $p>2$, this
phenomenon does not occur. Andrea Previtali proved in \cite{previtali} that if
$p>2$, then the dimension of every complex irreducible representation of
$USp(2n,q)$ is a power of $q$. The proof of the next proposition is very
similar to Previtali's proof, and is only included here for the sake of
completeness. Before giving it, we should point out that the key step in the
proof (apart from using Isaacs' results) is the fact that the coefficients of
the power series expansion for the function $(1+t)^{1/2}$ only have powers of
$2$ in the denominators.
\begin{prop}\label{p:unip-symp}
Let $\Ga=Sp(A,\sg)$ be a generalized unipotent symplectic group, as in
Definition \ref{d:unip-symp}, where $k=\bF_q$ and $A$ is a finite dimensional
associative nilpotent algebra over $\bF_q$. If
$p=\operatorname{char}(\bF_q)>2$, then the dimension of every complex
irreducible representation of $\Ga$ is a power of $q$.
\end{prop}
\begin{proof}
We use Theorem D of \cite{isaacs}. Following \loccit, we say that a subgroup
$H\subseteq A^\times$ is {\em strong} if for every $\bF_q$-subalgebra
$B\subseteq A$, the order of the intersection $H\cap B^\times$ is a power of
$q$ (in particular, $\abs{H}$ itself must be a power of $q$). By Isaacs'
Theorem D, if $H\subseteq A^\times$ is a strong subgroup, then the dimension of
every $\rho\in\widehat{H}$ is a power of $q$. Thus it suffices to prove that
$\Ga$ is a strong subgroup of $A^\times$.

\mbr

To avoid confusion, we will think of $A^\times$ as the group of elements of the
form $g=1+x$ in what follows. Let $B\subseteq A$ be any $\bF_q$-subalgebra, and
let $C=B\cap\sg(B)$; this is a $\sg$-stable $\bF_q$-subalgebra of $A$, and
hence $\sg$ restricts to an anti-involution of $C$. Now if $g=1+x\in\Ga\cap
B^\times$, then $x\in B$, and since $1+\sg(x)=\sg(g)=g^{-1}$ must also lie in
$\Ga\cap B^\times$, we see that $x\in C$. It follows that
\[
\Ga\cap B^\times = \bigl\{ g\in C^\times \st g^{-1}=\sg(g) \bigr\},
\]
and hence, replacing $A$ with $C$, we are reduced to the following
\begin{lem}
With the assumptions of Proposition \ref{p:unip-symp}, the order of $\Ga$ is a
power of $q$.
\end{lem}
\textit{Proof.} Since $p>2$, we can write $A=A_+\oplus A_-$ (direct sum of
$\bF_q$-subspaces), where $A_{\pm}=\Ker(\sg\mp\id:A\to A)$. We will define an
explicit bijection between $\Ga$ and $A_-$, which will prove the lemma.

\mbr

Pick $x\in A$ and write $x=x_++x_-$, with $x_\pm\in A_\pm$. If $1+x\in\Ga$,
then in particular $1+\sg(x)$ must commute with $1+x$, whence $[x_+,x_-]=0$.
Thus
\[
0 = x+\sg(x)+x\sg(x)=2x_+ +x_+^2 - x_-^2,
\]
which we rewrite as
\[
(1+x_+)^2=1+x_-^2.
\]
We now show that, given $x_-\in A_-$, there exists a unique $x_+\in A_+$
satisfying this equation, and, moreover, this $x_+$ commutes with $x_-$. This
will imply that the map
\[
\Ga\rar{} A_-, \qquad 1+x\longmapsto x_-
\]
is a bijection, and will complete the proof. Given $x_-\in A_-$, let us define
$x_+^\circ\in A$ using Newton's binomial formula:
\[
x_+^\circ = -1+(1+x_-^2)^{1/2} = -1 + \sum_{j\geq 0} {1/2\choose j} \cdot x_-^{2j} =  \frac{x_-^2}{2} - \frac{x_-^4}{8} + \dotsb;
\]
there is no question of convergence because $A$ is nilpotent. All coefficients
have the form $k/2^m$, where $k\in\bZ$ and $m\in\bN$, whence the expression above is well defined.\footnote{In fact, a more general statement holds. If $p$ is a prime, $a\in\bZ_p$, and $j\in\bZ$, $j\geq 0$, then the binomial coefficient $a\choose j$ lies in $\bZ_p$. This is well known for $a\in\bZ$, and follows for all $a\in\bZ_p$ by continuity.}

\mbr

By construction, $x_+^\circ\in A_+$, $(1+x_+^\circ)^2=1+x_-^2$, and $x_+^\circ$
commutes with every element of $A$ that commutes with $x_-^2$. Now if $x_+\in
A_+$ satisfies $(1+x_+)^2=1+x_-^2$, then $x_+$ commutes with $x_-^2$ and
therefore with $x_+^\circ$. Since $(1+x_+)^2=(1+x_+^\circ)^2$, and since the
order of $A^\times$ is odd, it is immediate that $x_+=x_+^\circ$, which
completes the proof of the lemma, and hence of the proposition.
\end{proof}


\section{Maximal unipotent subgroup of $Sp_4$ (after G. Lusztig)}\label{s:usp4}

G.~Lusztig \cite{lusztig} classified the irreducible representations of a
maximal unipotent subgroup of $Sp_4(k)$, where $k$ is a finite field. He
observed that if $k$ has characteristic $2$ and order greater than $2$, then
the dimensions of some of these representations are not powers of $q$. The
significance of this observation is explained in \S\ref{ss:dim-irreps}.

\mbr

We reproduce Lusztig's classification in \S\ref{ss:lusztig-answer}. We describe
two slightly different ways to derive it (see \S\ref{ss:littlegroups} and
\S\S\ref{ss:approach2}--\ref{ss:proof-ex2}). Both methods rely on standard
material from \S\S\ref{ss:abstract1}--\ref{ss:unip-symplectic}. In these
subsections the ground field $k$ can be arbitrary.

\subsection{The Siegel parabolic}\label{ss:abstract1} Let
$(V,\om)$ be a finite dimensional symplectic vector space over $k$ (later we
will assume that $\dim V=4$, but for now $V$ is arbitrary). Let $Sp(V)$ denote
the group of all linear automorphisms of $V$ that preserve $\om$. Fix a
Lagrangian subspace $L\subset V$ and a complementary Lagrangian $L^*\subset V$,
i.e., $V=L\oplus L^*$. This notation is justified by the fact that $\om$
induces a perfect pairing $L\times L^*\to k$, and hence identifies $L^*$
canonically with the dual space of $L$.

\mbr

Define
\[
P=\bigl\{g\in Sp(V)\st g(L^*)\subseteq L^* \bigr\},
\]
a parabolic subgroup of $Sp(V)$. Then
\[
L_P=\bigl\{g\in Sp(V) \;\st\; g(L^*)\subseteq L^*,\ \ g(L)\subseteq L\bigr\}
\]
is a Levi subgroup of $P$, and
\[
U_P = \bigl\{ g\in Sp(V) \st g\vert_{L^*}=\id_{L^*} \bigr\}
\]
is the unipotent radical of $P$. It is possible to describe $U_P$ more
explicitly. Namely, if $g\in U_P$, write $g=1+C$ for some $C\in\End(V)$. Then
$C(L^*)=0$, and it is easy to check that the condition that $g$ preserves $\om$
is equivalent to the condition that $C(L)\subseteq L^*$ and that the bilinear
form $B(u,w)=\om(Cu,w)$ on $L$ is {\em symmetric}. This gives the
identification
\[
U_P \cong (Sym^2 L)^*;
\]
more precisely, $U_P$ is the underlying additive group of the vector space
$(Sym^2 L)^*$. Since $L_P\cong\Aut(L)$, we obtain
\[
P \cong \Aut(L) \ltimes (Sym^2 L)^*
\]
(with $(Sym^2 L)^*$ being the normal subgroup).

\subsection{The maximal unipotent subgroup}\label{ss:abstract2}
If we choose a maximal unipotent subgroup $H\subset \Aut(L)$, then $U=H\ltimes
(Sym^2 L)^*$ is a maximal unipotent subgroup of $P$, and hence of $Sp(V)$. Note
also that the dual unipotent group $U_P^*$ can obviously be identified with the
underlying additive group of $Sym^2 L$. The latter can be naturally viewed as
the space of homogeneous quadratic polynomials on $L^*$, and the action of $U$
(and of $P$) on $Sym^2 L$ is induced by its natural action on $L^*$.

\mbr

Since $U_P$ is abelian, the action of $U$ on $U_P$ by conjugation factors
through the action of $H$, which will be denoted by $B\mapsto B^h$ (for $B\in
U_P,h\in H$). Explicitly, if $B\in U_P$ is thought of as a symmetric bilinear
form on $L$, and $h\in H\subset\Aut(L)$, then $B^h(u,w)=B(h^{-1}u,h^{-1}w)$.

\subsection{Matrix realization}\label{ss:unip-symplectic} Let us recall from
\S\ref{ss:gutkin-main} that $UL_n(k)$ denotes the group of unipotent
upper-triangular matrices of size $n$ over $k$. In order to relate the
construction of the maximal unipotent subgroup $U\subset Sp(V)$ given above to
the description of $USp_4(k)$ used in \cite{lusztig}, we recall one of the
standard realizations of the group $Sp_{2n}(k)$. If $M$ is a square matrix, we
will denote by $M^T$ the matrix obtained by reflecting $M$ with respect to its
{\em antidiagonal} (i.e., the diagonal going from the lower left corner to the
upper right corner). This is not to be confused with the usual transpose of
$M$, which is denoted by $M^t$. Now let $S\in GL_{2n}(k)$ be the diagonal
matrix with diagonal entries $(1,-1,1,-1,\dotsc,1,-1)$. Then $Sp_{2n}(k)$ can
be defined as the subgroup of $GL_{2n}(k)$ consisting of matrices $M$
satisfying $M^{-1}=SM^T S^{-1}$. With this realization, the intersection
$USp_{2n}(k):=Sp_{2n}(k)\cap UL_{2n}(k)$ is a maximal unipotent subgroup of
$Sp_{2n}(k)$.

\mbr

Further, let $J$ denote the $2n\times 2n$ antidiagonal matrix with all
antidiagonal entries equal to $1$. It is then easy to check that $M^T=JM^t
J^{-1}$, which implies that the condition $M^{-1} = SM^T S^{-1}$ is equivalent
to $M^t (JS) M = JS$. Thus $Sp_{2n}(k)$ can alternatively be described as the
group of linear automorphisms of $k^{2n}$ which preserve the symplectic form
$\om$ given by the matrix $JS$. This matrix is also antidiagonal, with
antidiagonal entries $-1,1,-1,1,\dotsc,-1,1$ (from the top right corner to the
bottom left corner). As a Lagrangian subspace $L\subset k^{2n}$ we choose the
span of the first $n$ standard basis vectors, and as a maximal unipotent
subgroup of $\Aut(L)\cong GL_n(k)$ we choose $H=UL_n(k)$. We now see that the
subgroup $U=H\ltimes (Sym^2 L)^*$ coincides with $USp_{2n}(k)$.

\subsection{Little groups method}\label{ss:littlegroups}
Observe now that if $n=2$, then $H\cong\bG_a$ is abelian. Hence, if $k$ is
finite, it is easy to classify all irreducible representations of $USp_4(k)$
and find their dimensions using the decomposition $USp_4(k)\cong H\ltimes
(Sym^2 L)^*$ and the ``little groups method'' of Wigner and Mackey. We leave it
as an exercise\footnote{Hint: if $\chr k=2$ the $H$-stabilizers of some
elements of $(Sym^2 L)\setminus\{0\}$ have order $2$, while if $\chr k\neq 2$
all the stabilizers are trivial. To see this, identify $Sym^2 L$ with the space
of polynomials $f\in k[x]$ of degree not greater than $2$ and $H$ with the
group of translations $x\mapsto x+a$, $a\in k$. If $\chr k=2$ and $f$ has two
{\em distinct} roots in $\overline{k}$, then the group of translations
$x\mapsto x+a$ preserving $f$ has order $2$.}; having done it, the reader will
see that if $q$ is a power of $2$, then some irreducible representations have
dimension $q/2$ (which is not a power of $q$ unless $q=2$). The particular case
of the little groups method that suffices to do the exercise can be summarized
as follows.

\begin{prop}\label{p:little-groups}
Let $H$, $A$ be finite abelian groups, let $H$ act on $A$ by group
automorphisms, and form the semidirect product $G=H\ltimes A$. Let $A^*$ denote
the Pontryagin dual of $A$ and consider the induced action of $H$ on $A^*$.
There is a natural bijection
\[
\Gh \overset{\simeq}{\longleftrightarrow} \bigl\{ \text{pairs } (\Om,\psi)
\bigr\},
\]
where $\Om\subset A^*$ is an $H$-orbit and $\psi$ is a $1$-dimensional
character of the stabilizer $H^\chi\subseteq H$ for some $\chi\in\Om$. $($Note
that $H^\chi$ does not depend on the choice of $\chi\in\Om$ because $H$ is
abelian.$)$ Given a pair $(\Om,\psi)$ of this form, the corresponding
irreducible representation of $G$ is constructed as
\[
\rho(\Om,\psi) = \Ind_{H^\chi\ltimes A}^G
(\widetilde{\psi}\otimes\widetilde{\chi}),
\]
where $\widetilde{\psi}$ and $\widetilde{\chi}$ are the natural extensions of
$\psi$ and $\chi$ to $H^\chi\ltimes A$, defined by
$\widetilde{\psi}(ha)=\psi(h)$ and $\widetilde{\chi}(ha)=\chi(a)$ for $h\in
H^\chi$ and $a\in A$, respectively.
\end{prop}

\begin{rem}\label{r:little-groups}
The little groups method can be formulated in much greater generality; for
instance, it can be used to classify the irreducible representations of a
finite group $G$ with a nontrivial normal subgroup $N$ in terms of irreducible
representations $\rho$ of $N$ and irreducible representations of the
corresponding ``little groups'' (the little group $G^\rho\subset G$ is the
stabilizer of the isomorphism class of $\rho$). In this context the method is
usually known under the name ``Clifford theory''. We refer the reader to
\cite{n-t}, \S3.3 for a nice exposition of this topic.
\end{rem}

\begin{rem}
The method described in this subsection has the following drawback. In
characteristic $2$ there exists a purely inseparable isogeny $\xi:Sp_4\to Sp_4$
such that the parabolic $\xi(P)$ is not conjugate to $P$ (i.e., $\xi(P)$ is the
stabilizer of a line). We can choose $\xi$ so that $USp_4$ is stable under
$\xi$. Then $\xi:USp_4\to USp_4$ induces an automorphism of $USp_4(\bF_{2^k})$.
Unfortunately, the approach to the classification of the irreducible
representations of $USp_4(\bF_{2^k})$ presented above breaks the
$\xi$-symmetry. However, this symmetry is preserved by the approach explained
below.
\end{rem}

\subsection{Second approach}\label{ss:approach2}
We now explain a different approach to the classification of irreducible
representations of $USp_4(\bF_q)$, which yields an answer that is closer in
spirit to the one described in \S7 of \cite{lusztig}. Consider first $U=USp_4$
as an algebraic group over an arbitrary field $k$. The quotient $U/[U,U]$ is
isomorphic to a direct sum of two copies of $\bG_a$ (more precisely, these two
copies correspond to the two simple roots for $Sp_4(k)$). It is easy to see
that $[U,U]$ is also isomorphic to a direct sum of two copies of $\bG_a$; in
particular $[U,U]$ is commutative.

\begin{lem}\label{l:ex1} If $\operatorname{char} k=2$, then the center
$Z$ of $U$ coincides with $[U,U]$. Otherwise, it is a subgroup of $[U,U]$
isomorphic to $\bG_a$ $($more precisely, it is the subgroup of $[U,U]$
corresponding to the maximal root of $Sp_4(k)${}$)$.
\end{lem}
The proof is given in \S\ref{ss:proof-ex1}.

\mbr

Now suppose that $k=\bF_q$ is a finite field of characteristic $2$. Fix a
character $\chi:Z(k)\to\{\pm 1\}$, and let us classify the irreducible
representations of $U(k)$ on which $Z(k)$ acts via $\chi$. This is the same as
classifying the irreducible representations of the amalgamated sum (=coproduct)
of $U(k)$ and $\{\pm 1\}$ over $Z(k)$ on which the element $-1\in\{\pm 1\}$
acts as multiplication by $-1$. Let us denote this amalgamated sum by
$H(\chi)$. By Proposition \ref{p:heisenberg}, such representations correspond
bijectively to the $1$-dimensional characters of the center of $H(\chi)$ which
take value $-1$ on the element $-1\in\{\pm 1\}$. It therefore remains to
describe the center of $H(\chi)$.

\begin{lem}\label{l:ex2} \begin{enumerate}[$($i$)$]
\item If $\chi$ is trivial, then $H(\chi)$ is commutative.
\item There are $2(q-1)$ characters $\chi$ for which the center of
$H(\chi)$ equals $\{\pm 1\}$.
\item For all other characters $\chi$, the center of $H(\chi)$ is
a direct sum of $\{\pm 1\}$ and two other subgroups of order $2$.
\end{enumerate}
\end{lem}

The proof is given in \S\ref{ss:proof-ex2}. At the end of the section we
present Lusztig's formulas for the irreducible characters of $USp_4(\bF_q)$,
where $q$ is a power of $2$. They can be obtained using the approach we have
just outlined.

\subsection{Proof of Lemma \ref{l:ex1}}\label{ss:proof-ex1} We return to the
coordinate-free framework described in \S\ref{ss:abstract1} and
\S\ref{ss:abstract2}. We assume that $\dim V=4$, so that $\dim L=2$. We fix a
line $\ell\subset L$ such that $H$ consists of all automorphisms $h$ of $L$
that act trivially on $\ell$ and on the quotient line $L/\ell$. As remarked
before, $H\cong\bG_a$, and in particular is abelian. Moreover,
$U_P\cong\bG_a^3$ in this case.

\mbr

Let us describe the center of $U$. We will write the elements of $U$ as pairs
$(h,B)$, where $h\in H$ and $B$ is a symmetric bilinear form on $L$. In this
notation, the multiplication on $U$ looks like
\[
(h,B)\cdot (h',B') = (hh',B^{(h')^{-1}}+B').
\]
Since $H$ is commutative, the condition that $(h,B)$ commutes with all of $U$
breaks up into two conditions: $(h,0)$ must commute with $(1,B')$ for every
$B'$, and $(1,B)$ must commute with $(h',0)$ for every $h'$. In other words, it
is equivalent to the following two conditions:
\begin{enumerate}[(a)]
\item $B^{(h')^{-1}}=B$ for all $h'\in H$, and
\item $(B')^{h^{-1}}=B'$ for every symmetric bilinear form $B'$ on
$L$.
\end{enumerate}

Since $H$ obviously acts faithfully on $(Sym^2 L)^*$, the second condition is
equivalent to $h=1$. The first condition is equivalent to
\begin{equation}\label{e:1}
B(h'u,h'w)=B(u,w) \qquad  \text{for all } h'\in H \text{ and all } u,w\in L.
\end{equation}
Fix $v\in\ell$, $v\neq 0$. Taking $w=v$ in the equation above and letting $h'$
be arbitrary shows that $B(v,v)=0$. In characteristic $2$ this condition is
also sufficient for \eqref{e:1} to hold, since for any $u\in L$ we have
$h'u=u+\al v$ for some $\al\in k$, whence
\[
B(h'u,h'u)=B(u,u)+\al^2B(v,v).
\]
However, in characteristic different from $2$, we have
\[
B(h'u,h'u) = B(u,u) + \al^2 B(v,v) + 2\al B(u,v),
\]
so the condition \eqref{e:1} is equivalent to $B(u,v)=0$ for all $u\in L$. We
conclude that the center of $U$ consists of all $(0,B)$ such that
$B(\ell,\ell)=0$ (resp., $B(L,\ell)=0$) if $\operatorname{char} k=2$ (resp.,
$\operatorname{char} k\neq 2$). Thus the center is isomorphic to $\bG_a^2$ in
characteristic $2$ and to $\bG_a$ in characteristic different from $2$.

\mbr

In any characteristic the commutant of $U$ equals the set of pairs $(0,B)$ such
that $B(\ell,\ell)=0$ (this group is isomorphic to $\bG_a^2$). So we see that
the center of $U$ is contained in $[U,U]$, and it equals $[U,U]$ if and only if
the characteristic of $k$ equals $2$. This proves Lemma \ref{l:ex1}.

\subsection{Proof of Lemma \ref{l:ex2}}\label{ss:proof-ex2} We now consider the
case where $k=\bF_q$, with $q$ a power of $2$. Note that specifying a character
$\chi:Z(k)\to\{\pm 1\}$ is the same as specifying a linear map $f:Z(k)\to k$.
Namely, given such an $f$, we obtain a character $\chi$ defined by
\[
\chi(B)=\psi_0(f(B)), \qquad B\in Z(k),
\]
where $\psi_0:k\to\cst$ is a fixed nontrivial additive character. As explained
above, $Z(k)$ is identified with the space of linear functionals $B:Sym^2 L\to
k$ which vanish on $Sym^2\ell\subset Sym^2 L$. Hence the dual space
$Z(k)^*=\Hom(Z(k),k)$ is naturally identified with the quotient of $Sym^2 L$ by
the line $Sym^2\ell$. Now modulo this line, every element of $Sym^2 L$ is
easily seen to be decomposable, and in fact, we have $4$ essentially different
possibilities for an element of $(Sym^2 L)/(Sym^2\ell)$:
\begin{enumerate}[(1)]
\item $0$;
\item $w\tens v$, where $v\in\ell$ and $w\in L\setminus\ell$;
\item $w\tens w$, where $w\in L\setminus\ell$;
\item $(w+v)\tens w$, where $v\in\ell$ and $w\in L\setminus\ell$.
\end{enumerate}
It is a trivial exercise to show that these in fact exhaust all the
possibilities. Moreover, since in the cases (2) and (3) $w$ is only determined
up to adding an element of $\ell$, it is easy to see that there are $(q-1)$
different elements of $(Sym^2 L)/(Sym^2\ell)$ of the form (2) and $(q-1)$
different elements of the form (3).

\mbr

On the other hand, the same argument that was used to find the center of $U(k)$
also shows that an element $(h,B)\in U(k)$ maps to a central element of
$H(\chi)$ if and only if the following two conditions hold:
\begin{enumerate}[(a$'$)]
\item $\chi(B^{(h')^{-1}}-B)=1$ for all $h'\in H$, and
\item $\chi((B')^{h^{-1}}-B')=1$ for every symmetric bilinear form $B'$ on
$L$.
\end{enumerate}

\mbr

Let $f\in Z(k)^*$ be a linear functional represented by an element of one of
the types (1)--(4) above, and let $\chi:Z(k)\to\{\pm 1\}$ denote the
corresponding character. Let us describe the center of $H(\chi)$. If $f$ is of
type (1), then $H(\chi)$ is obviously abelian. In all other cases the natural
projection $U(k)\to H(\chi)$ is surjective, so it suffices to describe the set
of all elements $(h,B)\in U(k)$ satisfying conditions (a$'$) and (b$'$) above.

\mbr

Assume that $f$ is of type (2), so that $f(B)=B(w,v)$ for every $B\in Z(k)$.
Then condition (a$'$) is equivalent to $\psi_0(B(w,v))=\psi_0(B(h'w,v))$ for
every $h'\in H$. But since $w\not\in\ell$, this obviously forces $B(v,v)=0$.
Similarly, condition (b$'$) is equivalent to $\psi_0(B'(hw,v))=\psi_0(B'(w,v))$
for every $B'$, which obviously forces $h=1$. Thus, in this case, the center of
$H(\chi)$ is the image of the center of $U(k)$, and hence coincides with $\{\pm
1\}$.

\mbr

Assume next that $f$ is of type (3). Then condition (a$'$) is equivalent to
$\psi_0(B(w,w))=\psi_0(B(h'w,h'w))$ for every $h'\in H$. Writing $h'w=w+\al v$,
the condition becomes equivalent to $\psi_0(\al^2 B(v,v))=1$, which again
forces $B(v,v)=0$ since the map $k\to k$, $\al\mapsto\al^2$, is bijective.
Condition (b$'$) is equivalent to $\psi_0(B'(hw,hw))=\psi_0(B'(w,w))$ for every
$B'$, which similarly forces $h=1$. Thus the center of $H(\chi)$ is again
$\{\pm 1\}$ in this case.

\mbr

Finally, let $f$ be of type (4). Then condition (a$'$) is equivalent to
$\psi_0(B(w+v,w))=\psi_0(B(w+\al v + v,w+\al v))$ for all $\al\in k$, which
reduces to $\psi_0((\al^2+\al)B(v,v))=1$. This no longer forces $B(v,v)=0$.
Indeed, the image of the map $k\to k$, $\al\mapsto\al^2+\al$, is an additive
subgroup of $k$ of index $2$, so there are two possibilities for $B(v,v)$.
Moreover, condition (b$'$) is equivalent to
$\psi_0(B'(hw+v,hw))=\psi_0(B'(w+v,w))$ for all $B'$. Write $hw=w+\be v$ for
some $\be\in k$; then the condition becomes $B'(w+\be v + v,w+\be v)=B'(w+v,w)$
(since $B'$ is arbitrary, we can ignore $\psi$). It is clear that again we have
two choices: $\be=0$ and $\be=1$.

\mbr

Combining these results, we see that when $f$ is of type (4), the center of
$H(\chi)$ has order $8$, and is the homomorphic image of the set of pairs
$(h,B)\in U(k)$ satisfying the conditions in the previous paragraph. This
completes the proof of Lemma \ref{l:ex2}.

\mbr

To summarize, we have explicitly described the center of $H(\chi)$ in each of
the cases (1)--(4) listed above, and this yields a classification of the
irreducible representations of $U(k)$ that act as $\chi$ on $Z(k)$. In
particular, one obtains their characters.

\subsection{Irreducible characters}\label{ss:lusztig-answer}
As in \S\ref{ss:unip-symplectic}, we identify $USp_4(k)$ with
$U(k)$, and it is easy to see that $U(k)$ can in turn be
identified with the group consisting of all matrices of the form
$$\left(\begin{matrix} 1&a&b&c\\
          0&1&d&b+ad\\
          0&0&1&a\\
          0&0&0&1                               \end{matrix}\right)$$
where $a,b,c,d\in k$. We will denote this matrix by $[a,b,c,d]$. Lusztig's
description (\cite{lusztig}, \S7) of the irreducible characters of $U(k)$ is as
follows.

\begin{enumerate}[(i)]
\item There are $q^2$ one dimensional characters $U(k)\rar{}\cst$ of the form
$$[a,b,c,d]\longmapsto\psi_0(xa+yd)$$ (one for each pair $x,y\in\bF_q$), where
$\psi_0:\bF_q\to\cst$ is a fixed nontrivial additive character. \sbr
\item There are $q-1$ irreducible characters of degree $q$ of the form
$[0,b,c,0]\mapsto q\cdot\psi_0(xb)$ (all other elements are mapped to $0$), one
for each $x\in\bF_q^\times$. \sbr
\item There are $q-1$ irreducible characters of degree $q$ of the form
$[0,b,c,0]\mapsto q\cdot\psi_0(xc)$ (all other elements are mapped to $0$), one
for each $x\in\bF_q^\times$. \sbr
\item There are $4(q-1)^2$ irreducible characters of degree $q/2$, one for each
quadruple \\
$(a_0,d_0,\eps_1,\eps_2)$ where
$a_0\in\bF_q^\times,d_0\in\bF_q^\times$ and
$\eps_1:\{0,a_0\}\rar{}\{\pm 1\}$, $\eps_2:\{0,d_0\}\rar{}\{\pm
1\}$ are group homomorphism. Namely, the character corresponding
to such a quadruple is given by
$[a,b,c,d]\mapsto(q/2)\cdot\eps_1(a)\cdot\eps_2(d)\cdot\psi_0\bigl(a_0^{-2}\cdot
d_0^{-1}\cdot (ba+ba_0+c)\bigr)$ if $a\in\{0,a_0\}$ and
$d\in\{0,d_0\}$; all other elements are sent to $0$.
\end{enumerate}

\subsection{Dimensions of representations}\label{ss:dim-irreps}
Let us fix $q=2^s$, $s\geq 2$. We use Lusztig's computation presented above to show that question (4) in \S\ref{ss:charsheaves-introduction} has a \emph{negative} answer for the group $USp_4$ over the finite field $\bF_q$.
\begin{lem}\label{l:4-cor}
Suppose that question $($4$)$ has a positive answer for a given
algebraic group $G$ over $\bF_q$, and let $\chi$ be an irreducible
character of $G(\bF_q)$ over $\ql$. Then there exist a complex
$\cF\in D^b_c(G\tens_{\bF_q}\bF,\ql)$ and an isomorphism
$\psi:\Fr^*\cF\rar{\simeq}\cF$ satisfying the following property.
For each $n\in\bN$, let $\psi_n:(\Fr^n)^*\cF\rar{\simeq}\cF$
denote the isomorphism given by
\[
\psi_n = \psi \circ \Fr^*\psi \circ (\Fr^2)^*\psi \circ\dotsb\circ
(\Fr^{n-1})^*\psi,
\]
and let $t_n:G(\bF_{q^n})\rar{}\ql$ denote the corresponding trace
function, as in \S\ref{ss:charsheaves-introduction}. Then
$t_1=\chi$, and for each $n\in\bN$, $t_n$ is an irreducible
character of $G(\bF_{q^n})$.
\end{lem}
This statement is obviously a special case of the property (4-iii)
in \S\ref{ss:charsheaves-introduction}. We also have
\begin{prop}\label{p:geometric-progressions}
Suppose $\La\subset\ql^\times$ is a finite nonempty subset,
suppose we are given a function $\La\to\bZ\setminus\{0\}$ denoted
$\la\mapsto k_\la$, and form
\begin{equation}\label{e:F}
F(n) = \sum_{\la\in\La} k_\la \la^n \qquad \text{for }
n=1,2,3,\dotsc
\end{equation}
Then the formula
\begin{equation}\label{e:bad}
F(n) = \begin{cases} q/2 & \text{for } n=1, \\
q^n, \ q^n/2 \text{ or } 1 & \text{for } n\geq 2
\end{cases}
\end{equation}
cannot hold.
\end{prop}
The proof is given in \S\ref{ss:proof-progressions}. Now it is
easy to see that Lusztig's computation, together with the lemma
and the proposition above, imply that question (4) in
\S\ref{ss:charsheaves-introduction} has a negative answer for $G=USp_4$ over $\bF_q$. Indeed, if the answer is positive, let $\chi$ be an irreducible character of $G(\bF_q)$ with $\chi(1)=q/2$, and let $\cF$, $\psi$ be as in
Lemma \ref{l:4-cor}. Put $F(n)=t_n(1)$. Then \eqref{e:bad} holds
by Lusztig's computation. On the other hand, let
$\La\subset\ql^\times$ be the set of nonzero eigenvalues of $\psi$
acting on the stalks $\cH^i(\cF)_1$ of the cohomology sheaves of
$\cF$ at $1\in G(\bF_q)$, and let $k_\la$ be their multiplicities
taken with the sign $(-1)^i$. Then \eqref{e:F} is satisfied by
construction, and we have a contradiction.

\subsection{Proof of Proposition \ref{p:geometric-progressions}}\label{ss:proof-progressions}
This subsection is completely independent from the rest of the
paper. The statement of the proposition we need to prove is purely
algebraic, and hence we can replace $\ql$ with $\bC$. This will
allow us to use analytic arguments below.

\mbr

For each $\la\in\cst$, consider the character $f_\la:\bZ\to\cst$
given by $f_\la(n)=\la^n$. Let us recall that a subset
$S\subseteq\bN$ is said to have \emph{density} $\de$ if
\[
\lim\limits_{N\to\infty} \frac{\# \bigl(S\cap [1,N]\bigr)}{N} =
\de.
\]
\begin{lem}\label{l:1}
Given a subset $S\subseteq\bN$ of density $1$ and a function
$f:S\to\bC$, define
\[
I_S(f) = \lim\limits_{N\to\infty} \frac{1}{N} \sum_{n\in S\cap
[1,N]} f(n),
\]
provided that this limit exists. If $\la\in\cst$, $\abs{\la}\leq
1$, then $I_S(f_\la)$ exists and is given by
\[
I_S(f_\la) = \begin{cases} 1 & \text{if } \la=1, \\
0 & \text{if } \la\neq 1.
\end{cases}
\]
\end{lem}
\begin{proof}
Since $\abs{\la}\leq 1$, we have $\abs{f_\la(n)}\leq 1$ for all
$n\in\bN$. As $S$ has density $1$, this implies that
\[
\lim\limits_{N\to\infty} \frac{1}{N} \sum_{n\in (\bN\setminus
S)\cap [1,N]} f(n) = 0.
\]
Thus $I_S(f_\la)$ exists if and only if $I_\bN(f_\la)$ exists, and
if they do exist, they are equal. The case $S=\bN$ of the lemma is
completely straightforward.
\end{proof}

\mbr

\begin{rem}
It is clear that the operation $f\mapsto I_S(f)$ is a ``partially
defined linear functional'' on the space of functions $f:S\to\bC$.
Namely, if $f,g:S\to\bC$ are such that $I_S(f)$ and $I_S(g)$
exist, then for any $a,b\in\bC$, $I_S(af+bg)$ also exists an
equals $aI_S(f)+bI_S(g)$.
\end{rem}

\mbr

\begin{cor}\label{c:lin-ind}
If $S\subseteq\bN$ is any subset of density $1$, then the
restrictions of the functions $\{f_\la\}_{\la\in\cst}$ to $S$ are
linearly independent.
\end{cor}
\begin{proof}
Otherwise there exists a finite nonempty subset $\La\subset\cst$
and a function $\La\to\cst$, $\la\mapsto c_\la$, such that
$\sum_{\la\in\La} c_\la f_\la(n)=0$ for all $n\in S$. Let
$\mu\in\La$ be an element of the maximal absolute value. From the
previous remark and the lemma, it follows that
$I_S\left(\sum_{\la\in\La} c_\la f_{\mu^{-1}\la}\right) =
c_\mu\neq 0$. On the other hand, $\sum_{\la\in\La} c_\la
f_{\mu^{-1}\la}(n) = \mu^{-n} \sum_{\la\in\La} c_\la f_\la(n)=0$
for all $n\in S$, which is a contradiction.
\end{proof}

\mbr

\begin{proof}[Proof of Proposition \ref{p:geometric-progressions}
with $\ql$ replaced by $\bC$] Assume that \eqref{e:bad} holds. We
will derive a sequence of statements that will eventually lead to
a contradiction. Let $R=\max_{\la\in\La}\abs{\la}$. By Corollary
\ref{c:lin-ind}, one of the values $q^n$ or $q^n/2$ must be
achieved infinitely many times in \eqref{e:bad} (indeed, this must
be so on a subset of $\bN$ whose complement does not have density
$1$). This immediately implies that $R\geq q$. Now let $\mu\in\La$
be an element with $\abs{\mu}=R$. Then on the one hand,
$I_\bN(f_{\mu^{-1}}\cdot F)=k_\mu$, a nonzero integer. On the
other hand, \eqref{e:bad} implies that for every $N\in\bN$, we
have
\[
\frac{1}{N} \sum_{n=1}^N \abs{f_{\mu^{-1}}(n)\cdot F(n)} \leq
\frac{1}{N} \sum_{n=1}^N \abs{q/\mu}^n = \frac{1}{N} \sum_{n=1}^N
(q/R)^n;
\]
note that the RHS is always $\leq 1$, and $\to 0$ as $N\to\infty$
if $R>q$. This means that $R=q$, and, in addition (since
$\abs{k_\mu}\geq 1$), the set of values $n\in\bN$ for which
$F(n)=q^n$ has density $1$. But then $F(n)=q^n$ for all $n\in\bN$
by Corollary \ref{c:lin-ind}, which contradicts \eqref{e:bad}.
\end{proof}


\appendix


\section{Characters of finite groups}\label{a:characters}

\subsection{Classical theory}\label{aa:classical} The main goal of this section is the study
of irreducible characters of a finite group $\Pi$ that are
invariant under a given automorphism $\phi$ of $\Pi$; this study
is motivated by its application to the theory of character sheaves
in the setting of the orbit method (see Section
\ref{s:char-sheaves-orbmethod}). On the other hand, the proofs of
the facts that we need are natural extensions of the proofs of the
classical results about irreducible characters of finite groups
(which is the case when $\phi=\id_\Pi$). Therefore we begin by
reviewing a bit of the classical theory, and then bring the
automorphism $\phi$ into the picture.

\mbr

We take $\bC$ as our field of coefficients to avoid introducing more notation,
even though all of what follows is equally valid over an arbitrary
algebraically closed field of characteristic zero. Let us fix a finite group
$\Pi$ and consider two {\em nonisomorphic} irreducible representations, $V$ and
$W$, of $\Pi$ over $\bC$. Recall that $\bC\Pi$ denotes the group algebra of
$\Pi$.

\begin{lem}
Let $m:\Pi\rar{}\bC$ be a matrix element of $V$. Then the element
\[
a_m = \sum_{\ga\in\Pi} m(\ga^{-1})\cdot\ga\in\bC\Pi
\]
acts by $0$ on $W$. \emph{(This lemma is in fact valid for
representation over any field.)}
\end{lem}
\begin{proof}
Recall that to say that $m$ is a matrix element of $V$ means that
$m(\ga)=f_0(\ga v_0)$ for some $v_0\in V$ and $f_0\in V^*$. Now
fix $w_0\in W$. The map $v\mapsto f_0(v) w_0$ is a linear map
$V\rar{} W$. Hence the map
\[
v\longmapsto \sum_{\ga\in\Pi} f_0(\ga^{-1}v) \ga w_0
\]
is a $\Pi$-equivariant linear map $V\rar{}W$, and is therefore $0$
by (the weak form of) Schur's lemma. But this map takes $v_0$ to
$a_m\cdot w_0$, whence $a_m\cdot w_0=0$. Since $w_0\in W$ is
arbitrary, the proof is complete.
\end{proof}

\begin{cor}\label{c:e-V}
In the same situation, define
\[
e_V = \frac{\dim V}{\card(\Pi)}\cdot \sum_{\ga\in\Pi}
\tr(\ga^{-1};V)\cdot\ga\in\bC\Pi.
\]
Then $e_V$ acts by the identity on $V$ and by $0$ on $W$.
\end{cor}
\begin{proof}
It is easy to see that $e_V(W)=0$. Indeed, let $(v_i)$ be a basis
of $V$, let $(f_i)$ be the dual basis of $V^*$, and let
$m_i(\ga)=f_i(\ga v_i)$ be the corresponding matrix elements of
$V$. Then
\[
e_V = \frac{\dim V}{\card(\Pi)}\cdot \sum_{i\in I} a_{m_i},
\]
whence the lemma applies.

\mbr

Moreover, it is clear that $e_V$ is a central element of $\bC\Pi$,
so by (the strong form of) Schur's lemma it acts by a scalar on
$V$. Thus to complete the proof of the corollary it suffices to
check that $\tr(e_V;V)=\dim V$. However, it is clear that
$\tr(e_V;\bC\Pi)=(\dim V)^2$, where $\bC\Pi$ is viewed as the
regular representation of $\Pi$. Since $\bC\Pi$ contains exactly
$\dim V$ copies of $V$, the corollary follows.
\end{proof}

\mbr

As a consequence, we immediately obtain the well known

\begin{cor}[Orthogonality relations]\label{c:orthogonality}
The irreducible characters of $\Pi$ are orthonormal with respect
to the inner product on the space of $\bC$-valued functions
$\Fun(\Pi)$ on $\Pi$ defined by
\[
\eval{f_1}{f_2} = \frac{1}{\card(\Pi)} \sum_{\ga\in\Pi} f_1(\ga)
f_2(\ga^{-1})
\]
\end{cor}
\begin{proof}
In the situation above, let $\chi_V$ and $\chi_W$ denote the
characters of $V$ and $W$, respectively. It is clear that
$\eval{\chi_W}{\chi_V}=(\dim V)^{-1}\cdot\tr(e_V;W)$, which
implies that $\eval{\chi_W}{\chi_V}=0$ in view of the previous
corollary. Similarly, $\eval{\chi_V}{\chi_V}=(\dim
V)^{-1}\cdot\tr(e_V;V)=1$.
\end{proof}
Later on we will need the following result. If $V$ is an
irreducible representation of $\Pi$, we define $I_V$ to be the set
of elements of $\bC\Pi$ that act by $0$ on every irreducible
representation of $\Pi$ that is not isomorphic to $V$. Then $I_V$
is a minimal (nonzero) two-sided ideal of $\bC\Pi$, and we call it
\emph{the minimal ideal corresponding to $V$}.
\begin{cor}\label{c:min-ideal}
With the same notation, consider the linear map
\[
\psi_V : \End_\bC(V) \rar{}\bC\Pi,
\]
\[
f\longmapsto \frac{\dim V}{\card(\Pi)} \sum_{\ga\in\Pi}
\tr(f\circ\ga^{-1};V)\cdot\ga.
\]
For each $f\in\End_\bC(V)$, the element $\psi_V(f)\in\bC\Pi$ acts as $f$ on
$V$, and if $W$ is an irreducible representation of $\Pi$ which is not
isomorphic to $V$, then $\psi_V(f)$ acts as $0$ on $V$. Consequently, $\psi_V$
is multiplicative, and it maps $\End_\bC(V)$ isomorphically onto $I_V$.
\end{cor}
\begin{proof}
The last statement of the corollary follows from the first. The morphism
$\bC\Pi\rar{}\End_\bC(V)$ induced by the action of $\Pi$ on $V$ is surjective,
thus it suffices to prove the first statement for the case where $f$ is given
by the action of an element $g\in\Pi$. Now
\[
\sum_{\ga\in\Pi}\tr(g\ga^{-1};V)\cdot\ga =
\sum_{\ga\in\Pi}\tr(g\ga^{-1};V)\cdot\ga g^{-1}\cdot g =
\sum_{\ga'\in\Pi}\tr(\ga^{\prime-1};V)\cdot\ga'\cdot g,
\]
whence $\psi_V(g)=e_V\cdot g$. In view of Corollary \ref{c:e-V}, this completes
the proof.
\end{proof}

\subsection{Convolution}\label{aa:convol-finite}
A special case of the following result has implicitly appeared in Section
\ref{s:orbmethod}. Let us recall that under the natural identification of
$\bC\Pi$ with the space $\Fun(\Pi)$ (where each $\ga\in\Pi$ corresponds to the
delta-function at $\ga$), the multiplication in the algebra $\bC\Pi$
corresponds to the convolution of functions, defined by
\[
(f_1*f_2)(\ga) = \sum_{h\in\Pi} f_1(h) f_2(h^{-1}\ga).
\]
\begin{prop}\label{p:idempotents}
As before, let $\chi_V,\chi_W:\Pi\to\bC$ denote the characters of nonisomorphic
irreducible representations $V$ and $W$ of $\Pi$, respectively. Then
\[
\chi_V*\chi_W = 0 \qquad\text{and}\qquad \chi_V*\chi_V = \frac{\card(\Pi)}{\dim
V} \cdot\chi_V.
\]
\end{prop}
Note that this result is a strengthening of Corollary \ref{c:orthogonality},
because if $f_1,f_2\in\Fun(\Pi)$, then
$\eval{f_1}{f_2}=(\card\Pi)^{-1}\cdot(f_1*f_2)(1)$, and on the other hand
$\chi_V(1)=\dim V$.
\begin{proof}
By definition, for any $\ga\in\Pi$, we have
\[
\frac{\dim V}{\card(\Pi)}\cdot (\chi_V*\chi_W)(\ga) = \frac{\dim
V}{\card(\Pi)}\cdot \sum_{h\in\Pi} \tr(h;V)\cdot\tr(h^{-1}\ga;W) =
\tr(e_V\cdot\ga;W).
\]
Since $e_V$ acts by $0$ on $W$ by Corollary \ref{c:e-V}, we see that
$\chi_V*\chi_W=0$. Replacing $W$ by $V$ in the last computation and using
Corollary \ref{c:e-V} again shows that $\frac{\dim V}{\card(\Pi)}\cdot
(\chi_V*\chi_V)=\chi_V$.
\end{proof}

\subsection{Twisted conjugacy classes}\label{aa:twisted} Let $\Pi$ be a finite
group and $\phi$ an automorphism of $\Pi$. We define the
\emph{$\phi$-conjugation action} of $\Pi$ on itself by
$\ga:x\mapsto\phi(\ga)x\ga^{-1}$. (Note that unless
$\phi=\id_\Pi$, this is {\em not} an action by group
automorphisms.) Its orbits will be called the
\emph{$\phi$-conjugacy classes} in $\Pi$. For the purposes of
Section \ref{s:char-sheaves-orbmethod} it is important to study
the relationship between functions $\Pi\rar{}\bC$ that are
constant on $\phi$-conjugacy classes on the one hand, and the
irreducible representations of $\Pi$ whose characters are
$\phi$-invariant on the other hand. To this end, let $\phh$ denote
the automorphism of $\pih$ induced by $\phi$. For every
$\rho\in\bigl(\pih\bigr)^{\phh}$, let us choose a realization
$\rho:\Pi\to\Aut(V_\rho)$ and an automorphism $\phi_\rho$ of
$V_\rho$ such that\footnote{Note that this condition determines
$\phi_\rho$ uniquely up to scaling.}
$\rho(\phi(\ga))=\phi_\rho^{-1}\rho(\ga)\phi_\rho$ for all
$\ga\in\Pi$. With this notation, we have

\begin{prop}\label{p:twisted-conjugacy}
If $\rho\in\bigl(\pih\bigr)^{\phh}$, the function
$\widetilde{\chi}_\rho$ on $\Pi$ defined by
$\widetilde{\chi}_\rho(\ga)=\tr(\phi_\rho\rho(\ga))$ is invariant
under $\phi$-conjugation. Moreover, together the functions
$\widetilde{\chi}_\rho$ form a basis for the space of functions on
$\Pi$ that are invariant under $\phi$-conjugation.
\end{prop}

The set of $\phi$-conjugacy classes in $\Pi$ is naturally
identified with $H^1(\bzh,\Pi)$, where $\bzh$ acts on $\Pi$ via
$\phi$. Thus we obtain the following
\begin{cor}
In the same situation, we have $\card
H^1(\bzh,\Pi)=\card\bigl(\pih\bigr)^{\phh}$.
\end{cor}

\noindent
$\bigl($Note that in general there is no natural bijection between
$H^1(\bzh,\Pi)$ and $\bigl(\pih\bigr)^{\phh}$, as is already
demonstrated by the special case $\phi=\id_\Pi$.$\bigr)$

\subsection{Proof of Proposition \ref{p:twisted-conjugacy}} The proof will be included in the final version of the paper, but the idea is to let $d$ denote the order of $\phi$, consider the semidirect product $\Ga=\Pi\rtimes(\bZ/d\bZ)$, where $\bZ/d\bZ$ acts on $\Pi$ via $\phi$, and apply the standard results recalled above to the irreducible characters of the group $\Ga$.


\section{Heisenberg representations}\label{a:heisenberg}


Representations of finite groups of nilpotence class $\leq 2$ (in
particular, of the fake Heisenberg groups) can be studied and
classified using the following approach. In effect, it amounts to
a special case of the orbit method; however, this case is
technically simpler, and, historically, it was understood before
the general orbit method was developed. Thus we prefer to discuss
it independently of Section \ref{s:orbmethod}.
\begin{defin}\label{d:heisenberg}
We say that an irreducible representation $\rho:\Ga\to\Aut(V)$ of
a finite group $\Ga$ is a {\em Heisenberg representation} if the
quotient of $\Ga$ by $N=\{g\in \Ga \,\big\vert\, \rho(g)\in\cst\}$
is abelian. This amounts to the same as requiring that the image
of $\Ga$ under $\rho$ is a nilpotent group of class at most $2$,
i.e., that $\rho\bigl([\Ga,[\Ga,\Ga]]\bigr)=\{\id_V\}$.
\end{defin}
\begin{lem}\label{l:nondegenerate}
If $\rho$ is a Heisenberg representation and $N$ is as above, the
pairing
\begin{equation}\label{e:pairing}
\Ga/N \times \Ga/N \rar{} \cst,
\end{equation}
obtained by composing $\rho$ with the commutator in $\Ga$, is
nondegenerate.
\end{lem}
\begin{proof}
Suppose that $g\in \Ga$ is such that $\rho([g,g'])=1$ for all
$g'\in \Ga$. This means that $\rho(g)$ commutes with $\rho(\Ga)$,
and so $\rho(g)$ is a scalar by Schur's lemma, i.e., $g\in N$, as
required.
\end{proof}
\begin{cor}\label{c:char-heis}
With the same notation, the character of $\rho$ vanishes outside of $N$.
\end{cor}
\begin{proof}
Let $g\in\Ga$, $g\not\in N$. By the lemma, there exists $\ga\in\Ga$ such that
$\rho([\ga,g])\neq 1$. But $\rho(\ga g\ga^{-1})=\rho([\ga,g])\rho(g)$, and
since $\rho([\ga,g])$ is a scalar $\neq 1$ and $\tr(\rho(\ga
g\ga^{-1}))=\tr(\rho(g))$, this forces $\tr(\rho(g))=0$, as claimed.
\end{proof}
We denote by $\Heis(\Ga)\subseteq\gh$ the set of isomorphism classes of
Heisenberg representations of $\Ga$. It is a simple exercise to show that every
Heisenberg representation of $\Ga$ is induced from a $1$-dimensional
representation of a subgroup of $\Ga$; in fact, that subgroup can be taken to
be the preimage of any subgroup of $\Ga/N$ that is Lagrangian with respect to
the pairing \eqref{e:pairing}. There is also a more precise description of
$\Heis(\Ga)$:
\begin{prop}\label{p:heisenberg}
If $\Ga$ is a finite group, there is a natural bijection between $\Heis(\Ga)$
and the set $S(\Ga)$ of pairs $(\nu,\widetilde{\nu})$ consisting of a
$\Ga$-invariant character $\nu:[\Ga,\Ga]\to\cst$ and an extension
$\widetilde{\nu}$ of $\nu$ to the preimage in $\Ga$ of the center\footnote{Note that $\Ga$ acts on $[\Ga,\Ga]$ by conjugation, and if a character $\nu:[\Ga,\Ga]\to\cst$ is invariant under this action, then $\Ker\nu$ is a normal subgroup of $\Ga$, so the quotient $\Ga/\Ker\nu$ is also a group. Hence the definition of $S(\Ga)$ makes sense.} of
$\Ga/\Ker\nu$.
\end{prop}
\begin{proof}
We define a map
$\al:\Heis(\Ga)\rar{}S(\Ga)$ as follows. Given a Heisenberg representation
$\rho:\Ga\to\Aut(V)$, the subgroup $\rho([\Ga,\Ga])\subset\Aut(V)$ commutes
with $\rho(\Ga)$, so $\rho([\Ga,\Ga])$ consists of scalars
by Schur's lemma. Thus there exists a character $\nu:[\Ga,\Ga]\to\cst$ such
that $\rho(\ga)=\nu(\ga)\id_V$ for all $\ga\in[\Ga,\Ga]$. It is obviously
$\Ga$-invariant. We have $\Ker\nu\subseteq\Ker\rho$, whence the representation
$\rho$ factors through $\Ga/\Ker\nu$. Applying Schur's lemma again, we see that
the center, $Z$, of $\Ga/\Ker\nu$ acts by scalars on $V$, which determines an
extension $\widetilde{\nu}$ of $\nu$ to a character of the preimage in $\Ga$ of
$Z$. We put $\al(\rho)=(\nu,\widetilde{\nu})$.

\mbr

Next we define a map $\be:S(\Ga)\rar{}\Heis(\Ga)$. Given
$(\nu,\widetilde{\nu})\in S(\Ga)$, note that the $\Ga$-invariance of $\nu$
means that $[\Ga,[\Ga,\Ga]]\subseteq\Ker\nu$, so if $C\subseteq\Ga$ denotes the
preimage of the center of $\Ga/\Ker\nu$, then $[\Ga,\Ga]\subseteq C$. In
particular, $\Ga/C$ is an abelian group. Moreover, the character $\nu$
determines a commutator map $c_\nu:\Ga\times\Ga\rar{}\cst$,
$(\ga_1,\ga_2)\mapsto\nu([\ga_1,\ga_2])$, which by the definition of $C$
factors through a nondegenerate pairing
$\overline{c}_\nu:(\Ga/C)\times(\Ga/C)\rar{}\cst$. Let $L\subseteq\Ga/C$ be a
Lagrangian subgroup with respect to $\overline{c}_\nu$, and let
$\widetilde{L}\subseteq\Ga$ denote its preimage. By construction, $\nu$
vanishes on $[\widetilde{L},\widetilde{L}]$, whence the character
$\widetilde{\nu}:C\rar{}\cst$ admits a (possibly non-unique) extension to a
character $f:\widetilde{L}\rar{}\cst$. Consider the induced representation
$\rho_{L,f}=\Ind_{\widetilde{L}}^\Ga f$. (We use this notation because \emph{a
priori} $\rho_{L,f}$ depends both on the choice of $L$ and on the choice of
$f$.) Let us prove three assertions:
\begin{enumerate}[1)]
\item The representation $\rho_{L,f}$ is irreducible. This follows immediately from Mackey's
irreducibility criterion (\cite{serre-reps}, \S7.5) and the fact that
$L\subseteq\Ga/C$ is Lagrangian. Indeed, if $\ga\in\Ga$ is an element such that
$f(\ga x\ga^{-1})=f(x)$ for all $x\in\widetilde{L}$, then $f$ is trivial on
$[\ga,\widetilde{L}]$, which means that the image $\overline{\ga}$ of $\ga$ in
$\Ga/C$ is orthogonal to $L$ with respect to the pairing $\overline{c}_\nu$,
and therefore $\overline{\ga}\in L$, i.e., $\ga\in\widetilde{L}$.
\item The representation $\rho_{L,f}$ is Heisenberg. This is clear because $f$
is trivial on $[\Ga,[\Ga,\Ga]]$ by construction, and $[\Ga,[\Ga,\Ga]]$ is
normal in $\Ga$, so $\rho_{L,f}$ is also trivial on $[\Ga,[\Ga,\Ga]]$.
\item Up to isomorphism, $\rho_{L,f}$ depends only on $\widetilde{\nu}$ (and
$\nu$), but not on the choices of $L$ or $f$. To see this we use the standard
formula for the character of an induced representation (\cite{serre-reps},
\S7.2). Since $C$ is normal in $\Ga$ and $\widetilde{\nu}:C\rar{}\cst$ is
$\Ga$-invariant, this formula implies that the restriction of $\tr(\rho_{L,f})$
to $C$ is a multiple of $\widetilde{\nu}$ which is independent of $L$ and $f$.
On the other hand, it is clear that $C=\{g\in\Ga\, \bigl\lvert\,
\rho_{L,f}(g)\in\cst\}$, whence $\tr(\rho_{L,f})$ vanishes outside of $C$ by
the previous step and by Corollary \ref{c:char-heis}. So $\tr(\rho_{L,f})$
depends only on $(\nu,\widetilde{\nu})$, and hence so does $\rho_{L,f}$.
\end{enumerate}
This allows us to define $\be:S(\Ga)\rar{}\Heis(\Ga)$ by
$\be(\nu,\widetilde{\nu})=\rho_{L,f}$ with the notation above. Finally, it is
clear that $\al\circ\be$ is the identity, and the fact that $\be\circ\al$ is
the identity follows from the formula for the character of an induced
representation mentioned above.
\end{proof}



\section{A ``reduction process'' for finite nilpotent
groups}\label{a:reduction}

In this appendix we prove that every irreducible representation of
a finite nilpotent group can be obtained canonically by inducing a
Heisenberg representation (Appendix \ref{a:heisenberg}) of a subgroup.
The method of proof also yields an almost canonical construction
of polarizations. For definiteness we will work with representations over the field $\bC$.

\subsection{Reminder on induced representations}\label{ss:induced}
Let $V$ be a finite dimensional complex vector space and $\rho: \Ga\to\Aut(V)$
an irreducible representation of a group $\Ga$. We define a {\em realization of
$\rho$ as an induced representation} to be a direct sum decomposition
\begin{equation}\label{e:star}
V = \bigoplus_{i\in I} V_i, \qquad V_i\neq (0),
\end{equation}
such that each $g\in \Ga$ takes each $V_i$ to some $V_{g(i)}$, where $g(i)\in
I$. Since $\rho$ is irreducible, these conditions imply that $\Ga$ acts
transitively on $I$. For the sake of naturality, we will not fix any particular
$i\in I$; but if one chooses $i$, one can easily see that $\rho$ is induced
from the representation of $\Ga_i=\{g\in \Ga\big\vert g(i)=i\}$ in the subspace
$V_i$.

\subsection{Reduction process for finite nilpotent groups}\label{ss:reduction} We claim that
{\em if $\rho:\Ga\to\Aut(V)$ is an irreducible representation of a
finite nilpotent group $\Ga$, then there is a canonical
realization \eqref{e:star} of $V$ as an induced representation,
such that the representation of each $\Ga_i$ in $V_i$ is
Heisenberg in the sense of Definition \ref{d:heisenberg}}. Here,
``canonical'' means ``constructed without making any choices'';
more precisely, given two irreducible representations
$\rho:\Ga\to\Aut(V)$, $\rho':\Ga'\to\Aut(V')$ and compatible
isomorphisms $\phi:\Ga\rar{\simeq}\Ga'$, $\psi:V\rar{\simeq}V'$,
the decomposition \eqref{e:star} for $V$ is mapped to that for
$V'$ under $\psi$.

\mbr

This canonical decomposition of $V$ can be constructed as follows. Put
$N=\{g\in \Ga \bigl\lvert \rho(g)\in\cst\}$, and let $Z$ be the center of
$\Ga/N$. Of course, unless $\Ga/N$ is abelian, the pairing \eqref{e:pairing} is
not defined, but one always gets a pairing
\begin{equation}\label{e:pairing2}
Z\times(\Ga/N)\rar{}\cst
\end{equation}
by composing $\rho$ with the commutator in $\Ga$. Let $Z_0\subseteq Z$ denote
the kernel of the restriction of this pairing to $Z\times Z$, and let $A\subset
\Ga$ be the preimage of $Z_0$. By construction, $A$ is a normal subgroup of
$\Ga$ (because $A/N$ is central in $\Ga/N$). Moreover, $\rho(A)$ is abelian.
Put $S=\{\chi\in A^*\st V_\chi\neq (0)\}$, where $V_\chi$ is the
$\chi$-eigenspace of $A$ in $V$. The decomposition
\[
V = \bigoplus_{\chi\in S} V_\chi
\]
satisfies the conditions of \S\ref{ss:induced}. The following lemma shows that
if $\rho$ is not Heisenberg, then $\operatorname{card}(S)>1$, which allows us
to proceed by induction.
\begin{lem}
If $\rho(A)\subseteq\cst$, then $\rho$ is a Heisenberg representation.
\end{lem}
\begin{proof}
If $\rho(A)\subseteq\cst$, then $A=N$, so $Z_0$ is trivial. This means that the
commutator pairing $Z\times Z\to\cst$ is nondegenerate, i.e., it induces an
isomorphism $Z\rar{\simeq} Z^*=\Hom(Z,\cst)$. By the definition of
\eqref{e:pairing2}, this isomorphism factors through $\Ga/N$, so $Z$ splits off
as a direct factor of $\Ga/N$. Since the center of $\Ga/N$ equals $Z$, this
implies that the center of $(\Ga/N)/Z$ is trivial. As $\Ga/N$ is nilpotent,
this means that $\Ga/N$ is abelian, i.e., $\rho$ is Heisenberg.
\end{proof}

\begin{rems}
\begin{enumerate}[(i)]
\item A classical theorem says that every irreducible representation of a
finite nilpotent group $\Ga$ can be realized as a representation
induced from a $1$-dimensional character of some
subgroup\footnote{This follows, e.g., from the construction
described above, because every Heisenberg representation of a
finite group can be induced from a $1$-dimensional representation
of a subgroup, see Appendix \ref{a:heisenberg}.}. However, there is no
{\em canonical} realization of this type: if $\rho$ is a
Heisenberg representation of $\Ga$ and $N$ is as before, then the
decompositions \eqref{e:star} with $\dim V_i=1$ correspond
bijectively to Lagrangian subgroups of $\Ga/N$, and it may happen
that $\Ga/N$ has no Lagrangian subgroups stable under
$\Aut(\Ga,\rho)$.
\item The previous remark is one of the reasons why we prefer to consider
Heisenberg representations as ``atoms'' rather than realizing them
as representations induced from $1$-dimensional characters. There
is also another reason. Namely, we hope that some version of the
``reduction process'' described above works in the geometric
setting of unipotent {\em algebraic} groups (instead of abstract
finite nilpotent groups). In this setting it is important to avoid
having to choose Lagrangian subgroups: e.g., the quotient of a
fake Heisenberg group (see \S\ref{ss:fake-heisenberg}) by its
commutator has {\em no algebraic Lagrangian subgroups at all},
because its dimension is odd.
\end{enumerate}
\end{rems}

\subsection{Reduction process for Lie rings. A construction of
polarizations}\label{ss:polarizations} We will describe an analogue of the
reduction process of \S\ref{ss:reduction} for finite nilpotent Lie rings. It
leads to an ``almost canonical'' construction of polarizations. After obvious
changes our definitions and constructions also apply to finite dimensional
nilpotent Lie algebras over a field\footnote{These changes are left to the
reader. The only difference is in the definition of $\fg^*$: if $\fg$ is a Lie
algebra over a field $k$, then $\fg^*$ stands for $\Hom_k(\fg,k)$.}.

\mbr

Given a finite Lie ring $\fg$, we put $\fg^*=\Hom_\bZ(\fg,\cst)$. If
$f\in\fg^*$, we define an alternating pairing $B_f:\fg\times\fg\rar{}\cst$ by
\[
B_f(x,y) = f\bigl( [x,y] \bigr), \qquad x,y\in\fg.
\]
Given an additive subgroup $\fa\subseteq\fg$, we write $\fa^{\perp_f}$ for the
orthogonal complement to $\fa$ with respect to $B_f$. We abbreviate
$\fg^{\perp_f}$ as $\fg^f$.

\begin{defin}
An element $f\in\fg^*$ will be called {\em Heisenberg} if its restriction to
$[\fg,[\fg,\fg]]$ is trivial, or, equivalently, if $[\fg,\fg]\subseteq\fg^f$.
\end{defin}

\begin{defin}
A {\em quasi-polarization} for $f\in\fg^*$ is a Lie subring $\fa\subseteq\fg$
which is coisotropic as an additive subgroup (i.e., satisfies
$\fa^{\perp_f}\subseteq\fa$). A quasi-polarization $\fa$ at $f$ is {\em
Heisenberg} if the restriction of $f$ to $[\fa,[\fa,\fa]]$ is trivial.
\end{defin}

Let us also recall (\S\ref{ss:polarizations}) that a \emph{polarization} for $f\in\fg^*$ is a Lie subring $\fh\subseteq\fg$ such that $f([\fh,\fh])=\{1\}$ and $\fh$ is maximal among all additive subgroups of $\fg$ with this property.

\mbr

The following result is obvious:

\begin{lem}\label{l:obvious}\begin{enumerate}[$($a$)$]
\item If $f\in\fg^*$ is Heisenberg, then every maximal isotropic
additive subgroup of $\fg$ with respect to $B_f$ is a Lie subring of $\fg$, and
hence a polarization of $\fg$ at $f$.
\item If $\fa\subseteq\fg^*$ is a quasi-polarization at $f$, then every
polarization $($resp., quasi-polarization$)$ for $\fa$ at $f\bigl\lvert_{\fa}$
is also a polarization $($resp., quasi-polarization$)$ for $\fg$ at $f$.
\end{enumerate}
\end{lem}

According to this lemma, in order to prove the existence of a polarization of
$\fg$ at a given $f\in\fg^*$, it is enough to show the existence of a
Heisenberg quasi-polarization at $f$. We prove a stronger statement:

\begin{prop}\label{p:construction}
Given $f\in\fg^*$, there exists a canonical\footnote{In the sense explained in
\S\ref{ss:reduction}.} way of constructing a Heisenberg quasi-polarization
$\fa_f$ for $\fg$ at $f$. The construction is compatible with the reduction
process of \S\ref{ss:reduction}; more precisely, it satisfies the condition
stated below.
\end{prop}

\sbr

\noindent
\textbf{Compatibility condition.} Suppose that $\fg$ has nilpotence class $c$ such that
$c!$ is invertible on $\fg$. Let $\rho_f$ be the irreducible representation of
$\Ga:=\exp\fg$ corresponding to $f$, and let $\rho_f^{\fa_f}$ be the
irreducible representation of $A_f:=\Exp\fa_f$ corresponding to the restriction
$f\bigl\lvert_{\fa_f}$. Then $\rho_f^{\fa_f}$ is Heisenberg, and since $\fa_f$
is coisotropic, the orbit method \cite{kirillov,orbmethod} implies that
$\Ind_{A_f}^\Ga\rho_f^{\fa_f}\cong\rho_f$. {\em This realization of $\rho_f$ as
a representation induced from a Heisenberg representation coincides with the
realization obtained in \S\ref{ss:reduction}.}

\bbr

To prove Proposition \ref{p:construction}, we will show that the Heisenberg
quasi-polarization $\fa_f$ can be constructed by mimicking the inductive
argument used in \S\ref{ss:reduction}. First, let $\fn$ denote the maximal
ideal of $\fg$ that is contained in $\fg^f$, let $\widetilde{\fz}=\bigl\{
x\in\fg\big\vert [x,\fg]\subseteq\fn \bigr\}$, and let
$\fg_1=\widetilde{\fz}+\widetilde{\fz}^{\perp_f}$.

\begin{lem}\label{l:construction}
\begin{enumerate}[$($i$)$]
\item $\fg_1\subset\fg$ is coisotropic with respect to $B_f$.
\item $\fg_1\supseteq [\fg ,\fg ]$.
\item $\fg_1=\fg$ if and only if $f$ is Heisenberg.
\end{enumerate}
\end{lem}

Before proving the lemma, let us note that it implies Proposition \ref{p:construction}.
Namely, the canonical Heisenberg quasi-polarization $\fa_f$ is constructed by induction
on the order of $\fg$. If $f$ is already Heisenberg take $\fa_f:=\fg$. Otherwise, Lemma
\ref{l:construction} tells us that $\fg_1$ is a proper coisotropic Lie subring of $\fg$,
so  the canonical Heisenberg quasi-polarization $\fa_{f\lvert_{\fg_1}}$ is defined by the
induction assumption, and by
Lemma \ref{l:obvious}(b) we can take $\fa_f=\fa_{f\lvert_{\fg_1}}$.
Now let us prove Lemma \ref{l:construction}.

\begin{proof}
(i) We have $\fg_1^{\perp_f}\subseteq\widetilde{\fz}^{\perp_f}\subseteq\fg_1$.

\mbr

(ii) It suffices to show that
\begin{equation}\label{e:my}
\widetilde{\fz}^{\perp_f}\supseteq[\fg ,\fg ],
\end{equation}
i.e., that
$B_f([x,y],z)=1$ for $z\in\widetilde{\fz}$ and $x,y\in\fg$. Write $B_f([x,y],z)$  as
$f([[x,y],z])$, apply the Jacobi identity, and recall that
$[\fg ,\widetilde{\fz}]\subseteq\fn\subseteq\fg^f$.

\mbr

(iii) If $f$ is Heisenberg, then $\widetilde{\fz}=\fg$, and therefore $\fg_1 =\fg$.
Let us prove that if $\fg_1=\fg$ then
$f$ is Heisenberg. This is equivalent to $\fg/\fn$ being abelian. As $\fg/\fn$ is nilpotent,
to prove that it is abelian it suffices to show that its commutant has
zero intersection with its center. Our $\widetilde{\fz}$ is the preimage of this center
in $\fg$, so we have to show that if $\fg_1=\fg$ then
$\widetilde{\fz}\cap [\fg ,\fg ]\subseteq\fn$.

\mbr

By \eqref{e:my},
$\widetilde{\fz}\cap [\fg ,\fg ]\subseteq\widetilde{\fz}\cap\widetilde{\fz}^{\perp_f}
\subseteq\widetilde\fg_1^{\perp_f}=\fg^{\perp_f}=\fg^f$. But
$\widetilde{\fz}$ is an ideal, so $\widetilde{\fz}\cap [\fg ,\fg ]$
is an ideal contained in $\fg^f$.  As $\fn$ is the maximal ideal with this property,
we see that $\widetilde{\fz}\cap [\fg ,\fg ]\subseteq\fn$.
\end{proof}

The following lemma shows that the construction of the quasi-polarization given above is compatible
with the reduction process of \S\ref{ss:reduction}.

\begin{lem}
With the same notation as above, assume that $\fg$ has nilpotence class $c$
such that $c!$ is invertible on $\fg$, let $\Ga=\Exp\fg$, let $\rho$ denote the
irreducible representation of $\Ga$ corresponding to the orbit of $f$, and let
$N$, $Z$ and $A$ be defined as in \S\ref{ss:reduction}. Then:
\begin{enumerate}[$($i$)$]
\item $N=\Exp\fn$;
\item $\Exp\widetilde{\fz}$ is the preimage of $Z$ in $\Ga$;
\item $A=\Exp\bigl(\widetilde{\fz}\cap\widetilde{\fz}^{\perp_f})$;
\item $\Exp(\fg_1)$ is the stabilizer in $\Ga$ of the $1$-dimensional character of $A$ induced by $f$.
\end{enumerate}
\end{lem}

\begin{proof}
The proof of (i)-(iii) is straightforward. It is also straightforward to check that
the stabilizer in $\Ga$ of the $1$-dimensional character of $A$ induced by $f$ equals
$\Exp(\fg'_1)$, where $\fg'_1=(\widetilde{\fz}\cap \widetilde{\fz}^{\perp_f})^{\perp_f}$. Finally,
$\fg_1:=\widetilde{\fz}+\widetilde{\fz}^{\perp_f}=\fg'_1$ (both  $\fg_1$ and  $\fg'_1$ contain the kernel
$\fg^f$ of our alternating pairing $B_f$, and it is easy to show that $\fg'_1/\fg^f =\fg_1/\fg^f$).
\end{proof}


\section{Vergne's construction of polarizations}\label{a:vergne}

The goal of this appendix is to present in detail a construction of
polarizations for completely solvable Lie algebras, due to Mich\`ele Vergne.
The key results are Theorems \ref{t:vergne1} and \ref{t:vergne2} (see also Remark \ref{r:not-every-polarization}). A related
result is Theorem \ref{t:classification}, which classifies
vector spaces equipped with a complete flag of subspaces and an alternating
bilinear form. We also give a reformulation of Vergne's construction
(cf.~\S\ref{aa:recursive}) which explains why it is very natural from the viewpoint of representation theory of nilpotent groups.

\mbr

We remark that the results of this section, with the exception of those in
\S\ref{aa:classification-triples}, have natural analogues for finite abelian
groups. In order to formulate and prove them one has to replace finite
dimensional vector spaces $V$ by finite abelian groups $A$; the dual space
$V^*$ in the sense of linear algebra by $A^*=\Hom(A,\cst)$; an alternating
bilinear form on $V$ by an alternating bi-additive map $B:A\times A\to\cst$;
subspaces of codimension $1$ by maximal proper subgroups; and arguments using
induction on the dimension of $V$ by arguments using induction on the order of
$A$. We leave the details as an exercise; see also
\cite{orbmethod}.

\subsection{Linear algebra}\label{aa:linear-algebra} Fix a
finite dimensional vector space $V$ over an arbitrary field $k$, equipped with
an alternating\footnote{We recall that this means that $B(v,v)=0$ for all $v\in V$. If $\operatorname{char}k=2$, this condition is \emph{stronger} than requiring $B$ to be skew-symmetric, i.e., $B(v,w)=-B(w,v)$ for all $v,w\in V$.} bilinear form $B:V\times V\to k$. If $W\subset V$ is a subspace,
we denote by $W^\perp$ its orthogonal complement in $V$ with respect to $B$.
One has $W^{\perp\perp}=W+\Ker B$, where $\Ker B=V^\perp$ is the kernel of $B$.
A subspace $W\subseteq V$ is said to be {\em isotropic} if $W\subseteq
W^\perp$. In this situation, $B$ induces an alternating bilinear form on
$W^\perp/W$. The following result is obvious:
\begin{lem}\label{l:d1}
Let $W\subseteq V$ be isotropic. If $W\subseteq W'\subseteq V$ and $W'$ is
isotropic, then $W'\subseteq W^\perp$ and $W'/W$ is an isotropic subspace of
$W^\perp/W$. Thus we obtain a bijection between the set of isotropic subspaces
of $V$ containing $W$ and the set of isotropic subspaces of $W^\perp/W$.
\end{lem}

\sbr

Even though this terminology is usually reserved for the situations where one
works with a symplectic form (i.e., a \emph{nondegenerate} alternating bilinear form), for the sake of brevity we will define a {\em
Lagrangian subspace} of $V$ to be a subspace which is maximal among the
isotropic subspaces of $V$, with respect to inclusion.
\begin{lem}\label{l:d2}
A subspace $L\subseteq V$ is Lagrangian if and only if $L^\perp=L$.
\end{lem}
\begin{proof}
By Lemma \ref{l:d1}, $L$ is Lagrangian if and only if $L\subseteq L^\perp$ and
$L^\perp/L$ has no nonzero isotropic subspaces. The latter condition means that
$L^\perp/L=(0)$ (otherwise any $1$-dimensional subspace of $L^\perp/L$ is
isotropic).
\end{proof}
\begin{lem}\label{l:d3}
Let $V'\subset V$ be a subspace of codimension $1$, and $L'$ a Lagrangian
subspace of $V'$ with respect to $B\big\vert_{V'}$.
\begin{enumerate}[$($i$)$]
\item A Lagrangian subspace $L\subseteq V$ contains $L'$ if and
only if $L\cap V'=L'$.
\item There is exactly one such $L$, namely,
$L=L'^\perp=L'^{\perp\perp}=L'+\Ker B$. \\ $($Here, $L'^\perp$ denotes the
orthogonal complement of $L'$ in $V$.$)$
\end{enumerate}
\end{lem}
\begin{proof}
(i) If $L\subseteq V$ is Lagrangian and $L\supset L'$, then $L\cap V'$ is an
isotropic subspace of $V'$ containing $L'$. But $L'$ is maximal among isotropic
subspaces of $V'$, so $L\cap V'=L'$.

\sbr

\noindent
(ii) By Lemma \ref{l:d1}, the Lagrangian subspaces of $V$ containing $L'$
correspond bijectively to the Lagrangian subspaces of $L'^\perp/L'$. But $L'$
is Lagrangian in $V'$, so $L'^\perp\cap V'=L'$, so $\dim(L'^\perp/L')\leq 1$.
Therefore $L'^\perp/L'$ has exactly one Lagrangian subspace, namely,
$L'^\perp/L'$. So there is exactly one Lagrangian subspace $L\subseteq V$
containing $L'$, namely, $L=L'^\perp$. As $L'^\perp$ is Lagrangian,
$L'^\perp=L'^{\perp\perp}$. Finally, $L'^{\perp\perp}=L'+\Ker B$.
\end{proof}

\mbr

Assume now that $V$ is equipped with an increasing filtration by subspaces
$V_i$ such that $V_0=(0)$, $V_n=V$ for some $n\geq 0$, and
$\dim(V_i/V_{i-1})\leq 1$ for all $i\geq 1$ (then for any subspace
$\widetilde{V}\subset V$, these conditions hold for the induced filtration
$\widetilde{V}_i:=\widetilde{V}\cap V_i$; this is the only reason why we impose
the inequality $\dim(V_i/V_{i-1})\leq 1$ rather than the equality
$\dim(V_i/V_{i-1})=1$). Put $B_{V_i}:=B\big\vert_{V_i}$.

\begin{thm}[Vergne's Theorem I]\label{t:vergne1}
There is exactly one Lagrangian subspace $L\subseteq V$ such that $L\cap V_i$
is Lagrangian in $V_i$ for every $i$. Namely, $L$ equals
\begin{equation}\label{e:L}
L(V_\bullet,B):=\sum_{i=1}^n \Ker B_{V_i}.
\end{equation}
\end{thm}
\begin{proof}
Lemma \ref{l:d3} allows one to proceed by induction on $n$.
\end{proof}

\begin{rems}\label{r:linalg}
\begin{enumerate}[(i)]
\item It is clear that $L(V_\bullet,B)$ is
isotropic in $V$, but the fact that it is Lagrangian is nontrivial. As
explained above, it follows from Lemma \ref{l:d3}. It also follows from Remark
\ref{r:L} and Proposition \ref{p:good} below.
\item Suppose that $\dim(V_i/V_{i-1})=1$ for some $i$. Then it is
easy to see that either $\Ker B_{V_{i-1}}$ is a codimension $1$ subspace of
$\Ker B_{V_i}$, or $\Ker B_{V_i}$ is a codimension $1$ subspace of $\Ker
B_{V_{i-1}}$.
\end{enumerate}
\end{rems}

The reader may prefer to skip the next subsection and go directly to
\S\ref{aa:vergne-cocycles}.

\subsection{A classification
theorem}\label{aa:classification-triples} In view of the discussion above, it
is natural to study triples of the form $(V,B,V_\bullet)$ consisting of a
finite dimensional vector space $V$ over $k$, an alternating bilinear form $B$
on $V$, and a complete flag of subspaces of $V$,
\[
V_\bullet\ : \qquad (0)=V_0\subset V_1 \subset \dotsb \subset V_n = V \qquad
(\dim V_j=j).
\]
We say that two such triples, $(V,B,V_\bullet)$ and $(V',B',V'_\bullet)$, are
{\em isomorphic} if there exists a $k$-linear isomorphism $\phi:V\to V'$ such
that $B'(\phi(v),\phi(w))=B(v,w)$ for all $v,w\in V$ and $\phi(V_j)=V'_j$ for
$0\leq j\leq n$. In this subsection we prove
\begin{thm}\label{t:classification} There exists a bijection
between the set of isomorphism classes of triples $(V,B,V_\bullet)$ with $\dim
V=n$, and the set of all involutions of the set $\{1,2,\dotsc,n\}$, i.e.,
elements $\sg\in S_n$ $($the symmetric group on $n$ letters$)$ such that
$\sg^2=1$.
\end{thm}
Let us define a basis $(e_1,\dotsc,e_n)$ of $V$ to be {\em good} with respect
to $(V_\bullet,B)$ if
\begin{enumerate}[1)]
\item for every $i$, the subspace $V_i$ is generated by
$e_1,\dotsc,e_i$; and
\item for each $i$ there exists at most one $j$ such that
$B(e_i,e_j)\neq 0$.
\end{enumerate}
Given such a basis, let $K$ be the set of indices $i$ such that $B(e_i,e_j)=0$
for all $j\in\{1,\dotsc,n\}$. For $i\not\in K$, let $\sg(i)$ denote the unique
$j$ such that $B(e_i,e_j)\neq 0$. It is clear that $\sg(i)\neq i$, that
$\sg(i)\not\in K$, and that $\sg(\sg(i))=i$. Defining $\sg(k)=k$ for all $k\in
K$, we obtain an involution of $\{1,2,\dotsc,n\}$ whose set of fixed points is
exactly $K$.

\begin{rem}\label{r:L}
Given a good basis, put $A=\{i \st \sg(i)>i\}$. Then the subspace
$L(V_\bullet,B)$ defined by \eqref{e:L} is spanned by the vectors $e_i$ where
$i\in K\cup A$. This immediately implies that $L(V_\bullet,B)$ is Lagrangian.
\end{rem}

\begin{prop}\label{p:good}
\begin{enumerate}[$($i$)$]
\item A good basis for $(V_\bullet,B)$ exists.
\item The involution $\sg$ is uniquely determined by
$(V_\bullet,B)$. In fact, it is determined by the relative position of the flag
$V_\bullet$ and the $($incomplete$)$ flag formed by the orthogonal complements
of the subspaces $V_i$ in $V$.
\end{enumerate}
\end{prop}
\begin{proof}
(i) By induction, we may assume the existence of a good basis
$(e'_1,\dotsc,e'_{n-1})$ for $V_{n-1}$. Put $K_{n-1}=\{i\leq n-1\st
B(e'_i,e'_j)=0 \text{ for all } j\leq n-1\}$. We can choose $e_n\in V\setminus
V_{n-1}$ so that $B(e'_i,e_n)=0$ for $i\not\in K_{n-1}$. Indeed, it suffices to
pick an arbitrary $e'\in V\setminus V_{n-1}$ and define
\[
e_n = e' + \sum_{i\not\in K_{n-1}} \frac{B(e',e'_i)}{B(e'_i,e'_{\tau(i)})}\cdot
e'_{\tau(i)},
\]
where $\tau$ is the involution of $\{1,\dotsc,n-1\}$ corresponding to the basis
$(e'_1,\dotsc,e'_{n-1})$. Now if $B(e'_i,e_n)\neq 0$ for at most one $i\in
K_{n-1}$, then $(e'_1,\dotsc,e'_{n-1},e_n)$ is already a good basis for $V$ and
we are done. If not, choose $t\in K_{n-1}$ with $B(e'_t,e_n)\neq 0$, and for
each $i\in K_{n-1}\setminus\{t\}$ replace $e'_i$ with
\[
e_i = e'_i + \frac{B(e_n,e'_i)}{B(e'_t,e_n)}\cdot e'_t.
\]
It is clear that this operation produces a good basis for $V$.

\sbr

\noindent
(ii) Let $r_{ij}$ be the rank of the map $V_i\to V_j^*$ induced by $B$. Put
$\eps_{ij}=r_{ij}-r_{i-1,j}-r_{i,j-1}+r_{i-1,j-1}$. Each $\eps_{ij}$ equals $0$
or $1$, and $\eps_{ij}=1$ if and only if $i\not\in K$ and $j=\sg(i)$.
\end{proof}

Theorem \ref{t:classification} follows easily from this proposition. Namely,
parts (i) and (ii) allow us to associate a well defined involution $\sg\in S_n$
to a triple $(V,B,V_\bullet)$, which of course depends only on the isomorphism
class of the triple. Moreover, if $(e_1,\dotsc,e_n)$ is a good basis for
$(V,B,V_\bullet)$, then by rescaling we may assume that $B(e_i,e_{\sg(i)})=1$
whenever $\sg(i)>i$, which clearly implies that knowing $\sg$ allows us to
recover $(V,B,V_\bullet)$ up to an isomorphism.

\subsection{Vergne's theorem for
$2$-cocycles}\label{aa:vergne-cocycles}
\begin{thm}[Vergne's Theorem II]\label{t:vergne2}
In the situation of Theorem \ref{t:vergne1}, suppose that $V$ is a Lie algebra
over $k$, each $V_i$ is a Lie ideal of $V$, and $B$ is a $2$-cocycle, i.e.,
\begin{equation}\label{e:2cocycle}
B\bigl([x,y],z\bigr) + B\bigl([y,z],x\bigr) + B\bigl([z,x],y\bigr) = 0 \qquad
\forall\, x,y,z\in V.
\end{equation}
Then $L=L(V_\bullet,B)$ is a Lie subalgebra of $V$.
\end{thm}
\begin{proof}
It is easy to deduce from \eqref{e:2cocycle} that $[\Ker B_{V_i},\Ker
B_{V_j}]\subseteq \Ker B_{V_i}$ for all $1\leq i\leq j\leq n$.
\end{proof}

If $B$ is the coboundary of $f\in V^*$, i.e., if $B(x,y)=f\bigl([x,y]\bigr)$,
then Theorem \ref{t:vergne2} says that $L(V_\bullet,B)$ is a polarization of
$V$ at $f$.  The existence of a complete flag of ideals $V_\bullet$ means by
definition that $V$ is completely solvable\footnote{We recall that every
nilpotent Lie algebra is completely solvable; every completely solvable Lie
algebra is solvable; and over an algebraically closed field every solvable Lie
algebra is completely solvable.}. So one obtains the following

\begin{cor}[M.~Vergne]
Let $V$ be a finite-dimensional Lie algebra over a field. If $V$ is completely
solvable, then every $f\in V^*$ admits a polarization.
\end{cor}

\begin{rem}\label{r:not-every-polarization}
It is known (even when $V$ is nilpotent) that not every polarization for $f$ can be obtained by Vergne's construction from some complete flag of ideals of $V$.
\end{rem}

\begin{rem}
The main reason why we prefer to formulate Theorem \ref{aa:vergne-cocycles} for
any 2-cocycle $B$ rather than in the particular case where $B$ is a coboundary
is that the proof becomes more transparent. In fact, Theorem
\ref{aa:vergne-cocycles} follows from this particular case because if $$0\to
k\to\hat V\to V\to 0$$ is the central extension corresponding to $B$ then the
pullback of $B$ to $V$ is a coboundary.

\mbr

On the other hand, Theorem \ref{aa:ass} below, which is an analog of Theorem
\ref{aa:vergne-cocycles} for associative algebras, does not seem to follow from
the case where $B$ is a coboundary.
\end{rem}

Theorem \ref{t:vergne2} has the following analog for associative algebras,
which we learned from Carlos A.~M.~Andr\'e's article \cite{andre} (see
Proposition 5.1 in the electronic version).

\begin{thm}\label{aa:ass}
In the situation of Theorem \ref{t:vergne1}, assume that $V$ is an associative
algebra over $k$, each $V_i$ is a two-sided ideal of $V$, and $B$ satisfies the
identity
\begin{equation}\label{e:stronger}
B(xy,z) + B(yz,x) + B(zx,y) = 0 \qquad \forall\, x,y,z\in V.
\end{equation}
Then $L$ is a $($multiplicatively closed$)$ subalgebra of $V$.
\end{thm}
\begin{proof}
It is easy to deduce from \eqref{e:stronger} that for all $1\leq i\leq j\leq n$
we have $(\Ker B_{V_i})\cdot(\Ker B_{V_j})\subseteq\Ker B_{V_i}$ and $(\Ker
B_{V_j})\cdot(\Ker B_{V_i})\subseteq\Ker B_{V_i}$.
\end{proof}

\begin{rems}
\begin{enumerate}[(i)]
\item
Condition \eqref{e:stronger} implies \eqref{e:2cocycle}.
\item
The space of alternating forms $B$ on $V$ satisfying \eqref{e:stronger}
contains the subspace of forms $B(x,y)=f(xy-yx)$ for $f\in V^*$. The
corresponding quotient space is dual to the degree $2$ cyclic
homology\footnote{We warn the reader that the space dual to cyclic homology is
{\em quite different} from cyclic cohomology.} of the associative algebra $V$.
\end{enumerate}
\end{rems}
In the next subsection we will need the following fact. The proof is completely
obvious.
\begin{lem}\label{l:d4}
\begin{enumerate}[$($i$)$]
\item Suppose that $V$ is a Lie algebra and $B$ is a $2$-cocycle.
Then the orthogonal complement of an ideal $I\subseteq V$ with respect to $B$
is a Lie subalgebra of $V$.
\item Suppose that $V$ is an associative algebra and $B$ satisfies
\eqref{e:stronger}. Then the orthogonal complement of a two-sided ideal
$I\subseteq V$ is a multiplicative subalgebra of $V$.
\end{enumerate}
\end{lem}

\subsection{Recursive procedure}\label{aa:recursive} Let
$V_\bullet$ and $B$ be as in Theorem \ref{t:vergne1}. Let $k$ be an index such
that $V_k$ is isotropic. Put $\widetilde{V}:=V_k^\perp$,
$\widetilde{V}_i:=\widetilde{V}\cap V_i$,
$\widetilde{B}:=B\big\vert{\widetilde{V}}$.
\begin{lem}\label{l:d5}
\begin{enumerate}[$($i$)$]
\item In this situation, $L(V_\bullet,B)=L(\Vt_\bullet,\Bt)$.
\item If $V$ is a Lie algebra, each $V_i$ is an ideal, and $B$ is
a $2$-cocycle, then $\Vt$ is a Lie subalgebra of $V$, each $\Vt_i$ is a Lie
ideal of $\Vt$, and $\Bt$ is a $2$-cocycle for $\Vt$.
\item If $V$ is an associative algebra, each $V_i$ is a two-sided
ideal of $V$, and $B$ satisfies \eqref{e:stronger}, then $\Vt$ is a
multiplicative subalgebra of $V$, each $\Vt_i$ is a two-sided ideal of $\Vt$,
and $\Bt$ satisfies \eqref{e:stronger} with $V$ replaced by $\Vt$.
\end{enumerate}
\end{lem}
\begin{proof} (ii) and (iii) follow immediately from (i) and Lemma
\ref{l:d4}. Let us prove (i). It is clear that $L(V_\bullet,B)\supseteq V_k$,
so $L(V_\bullet,B)\subseteq V_k^\perp=\Vt$. For every $i$, it follows that
$L(V_\bullet,B)\cap\Vt_i=L(V_\bullet,B)\cap\Vt\cap V_i=L(V_\bullet,B)\cap V_i$
is Lagrangian in $V_i$ by Theorem \ref{t:vergne1}, and hence, \emph{a
fortiori}, is Lagrangian in $\Vt_i$. Using Theorem \ref{t:vergne1} again, we
conclude that $L(V_\bullet,B)=L(\Vt_\bullet,\Bt)$.
\end{proof}

\mbr

Lemma \ref{l:d5} gives a recursive procedure of finding $L(V_\bullet,B)$. To
see this, note that if $B=0$, then $L(V_\bullet,B)=B$, and if $B\neq 0$, then
there exists $k$ such that $V_k$ is isotropic but $V_k\not\subseteq\Ker B$, so
$\Vt:=V_k^\perp$ is not equal to $V$. (The minimal $k$ such that
$V_k\not\subseteq\Ker B$ has this property.)

\mbr

Now assume that $V$ is a Lie algebra, each $V_i$ is an ideal, and
$B(x,y)=f\bigl([x,y]\bigr)$ for some $f\in V^*$. Then the above recursive
construction of $L(V_\bullet,B)$ can be reformulated as follows. If $f\in V^*$
is $V$-invariant (i.e., if $f$ is a Lie algebra character), then
$L(V_\bullet,B)$ is $V$ itself. Otherwise there exists $k$ such that
$f_k:=f\big\vert_{V_k}$ is a Lie algebra character, but $f_k$ is not invariant
with respect to the action of $V$ on $V_k^*$. Let $\Vt\subset V$ be the
stabilizer of $f_k$. Now replace $V$ by $\Vt$, $f$ by $f\big\vert_{\Vt}$, and
$V_i$ by $V_i\cap\Vt$. When formulated this way, the recursive procedure is
completely parallel to the classical one used to prove that any irreducible
unitary representation of a unipotent group over $\bR$ (or of a finite
nilpotent group) is induced from a $1$-dimensional character of a subgroup.


\section{Equivariant derived
categories}\label{a:equivariant-derived}

\subsection{Derived categories of constructible
complexes}\label{aa:derived-constructible} Let $k$ be an arbitrary
field, and $X$ a separated scheme of finite type over $k$. Given a
prime $\ell$ different from $\operatorname{char}k$, one knows how
to associate to $X$ a triangulated category $D^b_c(X,\ql)$, the
\emph{bounded derived category of constructible complexes of
$\ell$-adic sheaves on $X$}. Originally this was done for $k$
finite or algebraically closed (see, e.g., \cite{deligne-weil}).
For an arbitrary \emph{perfect} field $k$, the construction of
$D^b_c(X,\ql)$ is given in \cite{ekedahl}. Finally, if $k$ is not
perfect, we simply put
\[
D^b_c(X,\ql) := D^b_c(X\otimes_k k^{perf},\ql),
\]
where $k^{perf}$ is a perfect closure of $k$ (see
\cite{greenberg}), i.e., the minimal perfect subfield of
$\overline{k}$ containing $k$. (The motivation for this definition
is that $D^b_c(X,\ql)$ should only depend on the \'etale topos of
$X$, which does not change if we pass from $X$ to $X\tens_k
k^{perf}$.)

\mbr

In this article we will write $\sD(X)$ in place of $D^b_c(X,\ql)$
to simplify the notation. Also, we will refer to an object of
$\sD(X)$ simply as an \emph{$\ell$-adic complex on $X$}. The prime
$\ell$ different from the characteristic of the ground field will
be fixed throughout the discussion.

\mbr

If $M\in\sD(X)$ and $d\in\bZ$, we will denote by $M[d]$ the
complex obtained by shifting $M$ $d$ places to the left. The
functor $M\longmapsto M[d]$ is an autoequivalence of $\sD(X)$.
Another important autoequivalence is provided by the Tate twist.
It is defined as follows. For each $r\in\bN$, we have the locally
constant \'etale sheaf $\mu_{\ell^r}$ on $X$ given by
\[
\mu_{\ell^r}(U) = \bigl\{ \xi\in\Ga(U,\cO_U) \,\bigl\lvert\,
\xi^{\ell^r}=1 \bigr\} \quad\text{for every \'etale } U\rar{}X.
\]
Note that $\mu_{\ell^r}$ is naturally a sheaf of
$\bZ/\ell^r\bZ$-modules. Moreover, we have the sheaf morphisms
$\mu_{\ell^r}\rar{}\mu_{\ell^{r-1}}$ given by
$\xi\longmapsto\xi^\ell$, which are compatible with the module
structures and the obvious quotient maps
$\bZ/\ell^r\bZ\rar{}\bZ/\ell^{r-1}\bZ$. Hence we obtain a
projective system $\bigl\{\mu_{\ell^r}\bigr\}_{r\geq 1}$, which
determines a lisse $\bZ_\ell$-sheaf
$\mu_{\ell^\infty}\equiv\bZ_\ell(1)$. Extending scalars to
$\bQ_\ell$ or $\ql$, we obtain the sheaves $\bQ_\ell(1)$ and
$\ql(1)$ on $X$, called the \emph{$($first$)$ Tate twists}.

\mbr

The sheaves $\bQ_\ell(1)$ and $\ql(1)$ are invertible (with
respect to the usual tensor product on $\ell$-adic sheaves), and
their inverses are denoted by $\bQ_\ell(-1)$ and $\ql(-1)$,
respectively. If $n\in\bN$, we define
$\bQ_\ell(n)=\bQ_\ell(1)^{\tens n}$ and
$\bQ_\ell(-n)=\bQ_\ell(-1)^{\tens n}$, and we put
$\bQ_\ell(0)=(\bQ_\ell)_X$, the constant $\bQ_\ell$-sheaf of rank
$1$ on $X$. Similarly, $\ql(n)$ is defined for all $n\in\bZ$. More
generally, given $M\in\sD(X)$, we put $M(n)=M\tens_{\ql} \ql(n)$;
this is called the \emph{$n$-th Tate twist of $M$}. Note that if $k$ is algebraically closed, then each of the \'etale sheaves
$\mu_{\ell^r}$ in the definition above is already constant, and
hence $\bZ_\ell(1)\cong(\bZ_\ell)_X$,
$\bQ_\ell(1)\cong(\bQ_\ell)_X$, $\ql(1)\cong(\ql)_X$ are all
constant as well. However, even in this case the Tate twists do not have \emph{canonical} trivializations.

\subsection{Grothendieck's six functors}\label{aa:six-definitions}
We now briefly review a useful formalism for working with the
categories of $\ell$-adic complexes. However, we warn the reader
that for simplicity we drop the decorations ``$R$'' and ``$L$''
from the notation for all the functors we consider. Thus, for
instance, we denote by $\tens$ what is usually denoted by
$\overset{L}{\tens}$, and so on.

\mbr

Let $k$ and $\ell$ be as above, and let $f:X\rar{}Y$ be a morphism
of separated schemes of finite type over $k$. One can define four
exact functors:
\[
f_*:\sD(X)\rar{}\sD(Y), \qquad f^*:\sD(Y)\rar{}\sD(X),
\]
\[
f_!:\sD(X)\rar{}\sD(Y) \quad\text{and}\quad
f^!:\sD(Y)\rar{}\sD(X),
\]
called the \emph{pushforward}, \emph{pullback}, \emph{pushforward
with proper supports}, and \emph{extraordinary pullback},
respectively. The pairs $(f^*,f_*)$ and $(f_!,f^!)$ are adjoint
pairs of functors; in other words, $f^*$ is \emph{left adjoint} to
$f_*$, and $f^!$ is \emph{right adjoint} to $f_!$. We note also
that there is always a natural morphism of functors
$f_!\rar{}f_*$, which is an isomorphism whenever $f$ is
\emph{proper}. On the other hand, if $f$ is \emph{smooth of
relative dimension\footnote{This means that $f$ is flat and its
geometric fibers are regular and have dimension $d$ everywhere.}
$d$}, then $f^!\cong f^*[2d](d)$.

\mbr

There are also two exact bifunctors,
\[
\tens : \sD(X)\times\sD(X)\rar{}\sD(X) \quad\text{and}\quad
\boxtimes : \sD(X)\times\sD(Y)\rar{}\sD(X\times Y),
\]
called the \emph{tensor product} and the \emph{exterior tensor
product}, respectively. They are related to each other as follows.
Let $p:X\times Y\rar{}X$ and $q:X\times Y\rar{}Y$ be the two
projections. Then $M\boxtimes N\cong p^*(M)\tens q^*(N)$
canonically for all $M\in\sD(X)$, $N\in\sD(Y)$. On the other hand,
consider the special case $X=Y$, and let $\De_X:X\rar{}X\times X$
denote the diagonal embedding. Then $M\tens
N\cong\De_X^*(M\boxtimes N)$ for $M,N\in\sD(X)$.

\mbr

The bifunctor $\tens$ makes $\sD(X)$ a symmetric monoidal category
with unit object $(\ql)_X$. Moreover, it turns out that $\sD(X)$
is \emph{closed}, i.e., it has inner Homs. In other words, we have
a bifunctor
\[
\hom:\sD(X)^{op}\times\sD(X)\rar{}\sD(X)
\]
and trifunctorial isomorphisms
\[
\Hom(K\tens M,N) \cong \Hom(K,\hom(M,N)) \quad\text{for all }
K,M,N\in\sD(X).
\]
The ``six functors'' referred to in the title of this subsection
are $f_*,f^*,f_!,f^!,\tens,\hom$.

\subsection{Verdier duality}\label{aa:verdier} In the situation
above, let $X$ be a separated scheme of finite type over $k$, and
let $a:X\rar{}\Spec k$ denote the structure morphism. We define
$\bK_X:=a^!\ql$ and call it the \emph{dualizing complex} of $X$.
In the special case where $X$ is smooth of relative dimension $d$
over $k$, it follows that $\bK_X\cong\ql[2d](d)$. We define the
\emph{Verdier duality functor}
\[
\bD_X : \sD(X)^{op}\rar{}\sD(X)
\]
by $\bD_X(M)=\hom(M,\bK_X)$. There is a natural morphism of
functors $\Id_{\sD(X)}\rar{}\bD_X\circ\bD_X$ which is known to be
an isomorphism; in particular, $\bD_X$ is a triangulated\footnote{This means that $\bD_X$ commutes with shifts and takes distinguished triangles to distinguished triangles.}
anti-autoequivalence of $\sD(X)$. If $X$ is smooth of relative
dimension $d$ over $k$, we see that
$\bD_X(M)\cong\hom\bigl(M,(\ql)_X\bigr)[2d](d)$ for all
$M\in\sD(X)$.

\subsection{Formalism of the six functors}\label{aa:six-formalism}
The six functors introduced above enjoy a number of useful
properties. We will list a few of them that are implicitly used in
the computations appearing in this paper. Everywhere below we let
$k$, $\ell$ be as above, and we let $f:X\rar{}Y$ denote a morphism
of separated schemes of finite type over $k$.
\begin{enumerate}[(1)]
\item \emph{Duality and tensor product}: we have $\bD_X(M\tens
N)\cong\hom(M,\bD_X N)$, bifunctorially with respect to
$M,N\in\sD(X)$.
\item Each of the functors $f_*$, $f^*$, $f_!$, $f^!$ commutes
with Verdier duality; also, if $M\in\sD(X)$, then
$\bD_X(M(n))\cong(\bD_X M)(-n)$ for all $n\in\bZ$.
\item \emph{Projection formula}: we have $f_!\bigl((f^* L)\tens
M\bigr) \cong L\tens f_!M$, bifunctorially with respect to
$L\in\sD(Y)$, $M\in\sD(X)$.
\item We have
\[
\hom_Y(f_! M,N)\cong f_*\hom_X(M,f^! N),
\]
bifunctorially with respect to $M\in\sD(X)$, $N\in\sD(Y)$, and
\[
f^!\hom_Y(K,L)\cong\hom_X(f^*K,f^!L),
\]
bifunctorially with respect to $K,L\in\sD(Y)$.
\item \emph{Proper base change theorem}: Given a cartesian diagram
of separated schemes of finite type over $k$
\[
\xymatrix{
  Y' \ar[r]^{f'} \ar[d]_{\pi'} &   X' \ar[d]^\pi   \\
  Y \ar[r]^f &    X
   }
\]
there is a natural isomorphism $f^*\pi_!M\rar{\simeq}\pi'_!
f^{\prime *}M$ for every $M\in\sD(X')$.
\item \emph{Smooth base change theorem}: Given a cartesian diagram
as above where $f$ is \emph{smooth}, there is a natural
isomorphism $f^*\pi_*M\rar{\simeq}\pi'_* f^{\prime *}M$ for every
$M\in\sD(X')$.
\item \emph{Verdier duality theorem}: We have an isomorphism of
functors
\[
f_*\circ\bD_X \cong \bD_Y\circ f_! \,:\, \sD(X)^{op} \rar{}
\sD(Y).
\]
\item Similarly, we have $f^*\circ\bD_Y\cong\bD_X\circ f^!$,
$\bD_Y\circ f_*\cong f_!\circ\bD_X$ and $\bD_X\circ f^*\cong
f^!\circ\bD_Y$.
\end{enumerate}

\subsection{Convolution}\label{aa:convolution} Convolution of functions on a
finite group defined in \S\ref{aa:convol-finite} has an analogue
for $\ell$-adic complexes on an algebraic group.\footnote{To be more precise, it has two analogues: ``convolution with compact supports'' $*_!$, which is the only kind used in our paper and which we therefore call simply ``convolution'' and denote by $*$, and ``convolution without compact supports'' $*_*$, defined by $M*_*N=\mu_*(M\boxtimes N)$.} Let $G$ be an
algebraic group over a field $k$ (i.e., a reduced group scheme of
finite type over $k$). Then $G$ is automatically separated, so
$\sD(G)$ is defined. Let $\mu:G\times G\rar{}G$ denote the
multiplication morphism. With this notation, we define the
\emph{convolution bifunctor} for $G$,
\[
* : \sD(G)\times\sD(G) \rar{} \sD(G),
\]
by the formula
\[
M*N = \mu_!(M\boxtimes N).
\]
It makes $\sD(G)$ a monoidal category (which is not braided unless
$G$ is commutative) with unit object $\e:=1_!(\ql)=1_*(\ql)$,
the delta-sheaf at the identity element of $G$, where by abuse of
notation we write $1:\Spec k\rar{}G$ for the morphism which
defines the identity in $G$, and $\ql$ denotes the constant sheaf
of rank $1$ on $\Spec k$.

\subsection{Categories of equivariant objects}\label{aa:categories-equivariant}
Let $G$ be an algebraic group over a field $k$ acting on a separated variety
$X$ (i.e., a separated reduced scheme of finite type) over $k$. Write
$\al:G\times X\to X$ for the action morphism and $\pi:G\times X\to X$ for the
projection. If $\ell$ is a prime different from $\operatorname{char}k$, as
usual, one would like to have the notion of the ``equivariant derived
category'' $\sD^{BL}_G(X)$. The correct approach to the construction of $\sD^{BL}_G(X)$
is developed in the book \cite{ber-lunts}, which explains our notation; we say a few words about it in
\S\ref{aa:ber-lun} below. In this subsection we introduce a more naive
definition which turns out to be equivalent to the one studied
by Bernstein and Lunts when $G$ is unipotent (but not for other types of
algebraic groups), see Proposition \ref{p:conn-unip-equivar}.

\begin{defin}\label{d:equiv-derived}
We define $\sD_G(X)$ to be the category of pairs of the form $(M,\phi)$, where
$M\in\sD(X)$ and $\phi:\al^*M\rar{\simeq}\pi^*M$ is an isomorphism in
$\sD(G\times X)$ satisfying the following two conditions.
\begin{enumerate}[$($a$)$]
\item If $\mu:G\times G\to G$ is the product in $G$ and
$\pi_{23}:G\times G\times X\to G\times X$ is the projection along the first
factor $G$, then
\[
\pi_{23}^*(\phi)\circ(\id_G\times\al)^*(\phi) = (\mu\times\id_X)^*(\phi).
\]
More precisely, the two compositions of natural isomorphisms
\[
(\id_G\times\al)^*\al^* M \cong (\mu\times\id_X)^*\al^* M \xrar{\
(\mu\times\id_X)^*(\phi)\ } (\mu\times\id_X)^*\pi^* M \cong \pi_{23}^*\pi^* M
\]
and
\[
(\id_G\times\al)^*\al^* M \xrar{\ (\id_G\times\al)^*(\phi)\ }
(\id_G\times\al)^*\pi^* M \cong \pi_{23}^*\al^* M \xrar{\ \pi_{23}^*(\phi)\ }
\pi_{23}^*\pi^* M
\]
are equal.
\item If $1:\Spec k\to G$ is the morphism defining the identity in $G$ and
$1\times\id_X:X\to G\times X$ is the induced $k$-morphism, so that
$\al\circ(1\times\id_X)=\id_X=\pi\circ(1\times\id_X)$, then
$(1\times\id_X)^*(\phi)$ is the natural isomorphism
\[
(1\times\id_X)^*\al^*M \cong M \cong (1\times\id_X)^*\pi^* M
\]
\end{enumerate}
A morphism of objects $(M,\phi)\to(N,\psi)$ in $\sD_G(X)$ is defined as a
morphism $\nu:M\to N$ in $\sD(X)$ satisfying
$\phi\circ\al^*(\nu)=\pi^*(\nu)\circ\psi$.
\end{defin}

\mbr

In general, the category $\sD_G(X)$ is not triangulated (however, it is
triangulated when $G$ is unipotent by Proposition \ref{p:conn-unip-equivar}).
Nevertheless, it obviously inherits many structures from $\sD(X)$. For example,
we have the shift functors $(M,\phi)\mapsto(M[n],\phi[n])$ on $\sD_G(X)$ for
all $n\in\bZ$. It is also obvious that if $Y$ is another separated variety over
$k$ with a $G$-action and $f:X\to Y$ is a $G$-equivariant morphism, then
$f^*:\sD(Y)\to\sD(X)$ can be lifted to a functor $f^*:\sD_G(Y)\to\sD_G(X)$. On
the other hand, we have cartesian diagrams
\[
\xymatrix{
  G\times X \ar[rr]^{\id_G\times f} \ar[d]_{\al_X} & & G\times Y \ar[d]^{\al_Y} \\
  X \ar[rr]^f & & Y
   }
 \qquad\text{and}\qquad
\xymatrix{
  G\times X \ar[rr]^{\id_G\times f} \ar[d]_{\pi_X} & & G\times Y \ar[d]^{\pi_Y} \\
  X \ar[rr]^f & & Y
   }
\]
where $\al_X$, $\al_Y$ are the action morphisms and $\pi_X$, $\pi_Y$ are the
projections, so the proper base change theorem (see \S\ref{aa:six-formalism})
implies that $f_!:\sD(X)\to\sD(Y)$ lifts to a functor
$f_!:\sD_G(X)\to\sD_G(Y)$. Furthermore, $\al_Y$ and $\pi_Y$ are smooth
morphisms (because $G$ is reduced), so the smooth base change theorem (see
\S\ref{aa:six-formalism}) implies that $f_*:\sD(X)\to\sD(Y)$ also lifts to a
functor $f_*:\sD_G(X)\to\sD_G(Y)$.

\mbr

These remarks imply that if $\fg$ is another algebraic group over $k$, and we
are given an algebraic action of $G$ on $\fg$ by group automorphisms, then the
convolution bifunctor $*:\sD(\fg)\times\sD(\fg)\to\sD(\fg)$ defined in
\S\ref{aa:convolution} lifts to a bifunctor
$\sD_G(\fg)\times\sD_G(\fg)\to\sD_G(\fg)$, which we also denote by $*$.

\mbr

Finally, we have a forgetful functor $\sD_G(X)\to \sD(X)$ given by
$(M,\phi)\mapsto M$. If $G$ is connected and unipotent, this functor is exact
and fully faithful by Proposition \ref{p:conn-unip-equivar}.

\subsection{Relation with the Bernstein-Lunts approach}\label{aa:ber-lun}

In this subsection we show that for unipotent groups $G$ our naive definition
of $\sD_G(X)$ agrees with the ``scientific'' one introduced by Bernstein and
Lunts in \cite{ber-lunts}. If $G$ is an algebraic group over a field $k$ and
$X$ is a separated variety over $k$, we will write $\sD_G^{BL}(X)$ for the
$G$-equivariant bounded derived category of constructible complexes of
$\ell$-adic sheaves on $X$ defined in \emph{op.~cit.} It is a triangulated
category equipped with an exact forgetful functor $\ff:\sD_G^{BL}(X)\to\sD(X)$.

\begin{prop}\label{p:conn-unip-equivar}
Let $G$ be a unipotent group over a field $k$, and $X$ a separated variety over
$k$. Let $\ell$ be a prime different from $\operatorname{char}k$, and let
$\sD_G(X)$ be as in Definition \ref{d:equiv-derived}.
\begin{enumerate}[$($a$)$]
\item There is a natural equivalence of categories between $\sD_G(X)$ and
$\sD_G^{BL}(X)$ which commutes with shifts and is compatible with the forgetful
functors $\sD_G(X)\to\sD(X)$, $\sD_G^{BL}(X)\to\sD(X)$. In particular,
$\sD_G(X)$ is triangulated.
\item If $G$ is connected, the forgetful functor $\sD_G(X)\to\sD(X)$ is fully
faithful.
\end{enumerate}
\end{prop}

The proof will be included in the final version of the article.


\section{Duality for perfect commutative unipotent groups}\label{a:duality}

In this appendix we return to the setup of \S\ref{ss:perfect-duality} and explain the Serre duality theory for perfect connected
commutative group schemes over a perfect field of positive characteristic. This
theory was already mentioned in \S\ref{ss:perfect-duality}. It is used in the
orbit method for unipotent groups of small nilpotence class, because it allows
one to define the dual of the Lie ring scheme associated to such a group and
study the geometric properties of the coadjoint action. We follow \cite{saibi}, \S1.4, since that work is written in a language best suited for our purposes. We note that Saibi, in turn, refers to Begueri's work \cite{begueri}.

\subsection{Perfect schemes}\label{aa:perfect} To define the duality
functor, one realizes $\cC_k$ as a full subcategory of the category of perfect
group schemes over $k$. We explain this theory below.

\mbr

Recall that a scheme $S$ over $\bF_p$ is said to be {\em perfect} if the
absolute Frobenius morphism $\Phi=\Phi_{S,p}:S\to S$ defined in
\S\ref{ss:frobenii} is an isomorphism. Let $\Sch_p$ be the category of
$\bF_p$-schemes and $\Schp_p\subset\Sch_p$ the full subcategory of perfect
schemes. The embedding $\Schp_p\hookrightarrow\Sch_p$ has a right adjoint
$\Sch_p\to\Schp_p$, called the {\em perfectization functor} and denoted
$S\mapsto S^{perf}$. As a topological space, $S^{perf}=S$, and the structure
sheaf of $S^{perf}$ is the inductive limit of
$\cO_S\rar{\Phi^*}\cO_S\rar{\Phi^*}\dotsb$. Equivalently,
$S^{perf}=\limproj\bigl(S\overset{\Phi}{\longleftarrow}S\overset{\Phi}
{\longleftarrow}\dotsb\bigr)$. We refer the reader to \cite{greenberg} for more
details.

\mbr

A $k$-scheme is said to be {\em perfect} if it is perfect as an
$\bF_p$-scheme. The perfectization of any $k$-scheme is a perfect
$k$-scheme. A {\em perfect group scheme over $k$} is a group
object in the category of perfect $k$-schemes. This is the
same\footnote{If $k$ is not perfect, then a group object in the
category of perfect $k$-schemes is a group scheme not over $k$,
but over its perfect closure $k^{perf}$.} as a group scheme over
$k$ which is perfect as a scheme.

\mbr

If $G$ is a group scheme over $k$, then $G^{perf}$ is a perfect group scheme
over $k$. The relative Frobenius morphism $G\rar{}G^{(p)}$ induces an
isomorphism between the perfectizations. So we get a functor $G\mapsto
G^{perf}$ from the category $\cC_k$ defined in \S\ref{ss:perfect-duality} to
that of perfect group schemes over $k$. It is fully faithful. The category of
connected commutative {\em unipotent quasi-algebraic groups} over $k$, which we
denote by $\cC'_k$, is the essential image of this functor.

\subsection{Definition of $G^*$}\label{aa:definition-dual} We still keep the
assumption that $k$ is perfect\footnote{It is harmless because the
category $\cC_k$ defined in \S\ref{ss:perfect-duality} does not
change if $k$ is replaced by $k^{perf}$.}. If $G$ is any
commutative group scheme over $k$, we first define $G^*$ as a
functor on the category of perfect schemes $S$ over $k$ by the
formula
\begin{equation}\label{e:starstar}
G^*(S) = \Ext^1 (G\times S, \bQ_p/\bZ_p ),
\end{equation}
where $\Ext^1$ is computed in the category of commutative group schemes over
$S$ and $\bQ_p/\bZ_p$ is viewed as a discrete group scheme over $S$.

\mbr

One can prove (see \cite{begueri} and \cite{saibi}, Th\'eor\`eme 1.4.1) that if
$G\in\cC'_k$, then $G^*$ is representable by an object of $\cC'_k$, which we
also denote by $G^*$. Moreover, the functor $G\mapsto G^*$ is exact and
involutive on the category $\cC'_k$. Note that the perfectness of $S$ is
essential for formula \eqref{e:starstar} to hold, see Remark (iii) below.

\begin{rems}\label{rems:def-dual}
\begin{enumerate}[(i)]
\item The RHS of \eqref{e:starstar} equals $\Ext^1(G\times S,p^{-n}\bZ_p/\bZ_p)$
for any $n\in\bN$ such that $G$ is annihilated by $p^n$.
\item If $G\in\cC_k$, then for any perfect $k$-scheme $S$, one has $G(S)=G^{perf}(S)$
and
\[
\Ext^1(G\times S,\bQ_p/\bZ_p)=\Ext^1(G^{perf}\times S,\bQ_p/\bZ_p).
\]
So \eqref{e:starstar} still holds if $G$ or $G^*$ is considered as an object of
$\cC_k$.
\item Define the {\em coperfectization} $X^{coperf}$ of an $\bF_p$-scheme $X$ to
be the ind-scheme
\[
X^{coperf} := ``\limind{}" \bigl( X\xrar{\ \Fr_{X/k}\ } X^{(p)}\xrar{\
\Fr_{X^{(p)}/k}\ } \dotsb \bigr),
\]
where $``\limind{}"$ denotes the inductive limit in the category of sheaves of
sets on the category of schemes over $k$ with the fppf topology. If $G$ is a
group scheme over $k$, then $G^{coperf}$ is a group object in the category of
$k$-ind-schemes. We claim that if $G\in\cC_k$, then the RHS of
\eqref{e:starstar}, viewed as a functor on the category of {\em all}
$k$-schemes $S$, is ind-representable by $(G^*)^{coperf}$, where $G^*$ is
considered as an object of $\cC_k$. This is a formal consequence of the
equality \eqref{e:starstar} for $S^{perf}$ and the equality $\Ext^1(G\times
S,\bQ_p/\bZ_p)=\Ext^1(G\times S^{perf},\bQ_p/\bZ_p)$.
\end{enumerate}
\end{rems}

\subsection{Finite fields}\label{aa:finite-fields} Now let $k=\bF_q$. We will explain
why the commutative unipotent quasi-algebraic group $G^*$ defined
by \eqref{e:starstar} solves the problem posed in
\S\ref{ss:basic-examples}. Namely, the group $G^*(\bF)$ can be naturally identified with the group $\widehat{G(\bF)}$ (we recall that the latter is the direct limit of the Pontryagin duals $G(\bF_{q^n})^*$ with respect to the base change maps defined in \S\ref{ss:characters-abelian-groups}). This follows from the next proposition. We were unable to find exactly this result in the literature,
therefore we provide its proof; however, the idea is essentially
contained in \cite{lang}, even though Lang did not formulate it in
the same way.
\begin{prop}\label{p:lang}
Let $G$ be any connected commutative algebraic group over $\bF_q$. Let $A$ be a
finite abelian group. Then the morphism $\Hom(G(\bF_q),A)\to\Ext^1(G,A)$
arising from the short exact sequence
\[
 0 \rar{} G(\bF_q) \rar{} G \xrar{\Phi_q-\id} G \rar{} 0
\]
is an isomorphism.
\end{prop}
\begin{proof}
We have $\Hom(G,A)=0$ because $G$ is connected. To show that the morphism
$\Ext^1(G,A)\to\Ext^1(G,A)$ induced by $\Phi_q-\id:G\to G$ equals $0$, we note
that for any extension $0\to A\to \widetilde{G} \to G\to 0$, the morphism
$\Phi_q-\id:\widetilde{G}\to\widetilde{G}$ annihilates $A$.
\end{proof}

\subsection{Another point of view}\label{aa:avoid-perfect}
In practice, it is sometimes a little more convenient to work with
usual algebraic groups rather than with perfect quasi-algebraic
groups. The only drawback of this approach is that the Serre dual
is no longer defined up to canonical isomorphism. In this
subsection we briefly explain how to define all the objects we
wish to work with without passing to perfectization.

\mbr

As before, let $k$ be a perfect field of characteristic $p>0$.
\begin{defin}[cf. \cite{saibi}, D\'efinition 1.5.1(ii)]\label{d:dual-pair}
A \emph{dual pair} of unipotent $k$-groups is a triple $(G,G',\cE)$, where $G$
and $G'$ are connected commutative unipotent algebraic groups over $k$, and
$\cE$ is a $\bQ_p/\bZ_p$-torsor over $G'\times G$ with the property that the
induced $\bQ_p/\bZ_p$-torsor $\cE^{perf}$ over $(G')^{perf}\times G^{perf}$ is
the universal family of central extensions of $G^{perf}$ by $\bQ_p/\bZ_p$
parameterized by $(G')^{perf}$, and in particular identifies $(G')^{perf}$ with
$(G^{perf})^*$, the Serre dual of $G^{perf}$, defined as in
\S\ref{aa:definition-dual}.
\end{defin}

\begin{rems}\label{r:dual-pair}
\begin{enumerate}[(1)]
\item The definition above makes sense because if $H$ is a connected algebraic
group and $A$ is a finite discrete abelian group, the forgetful functor from
the groupoid of central extensions of $H$ by $A$ to the groupoid of $A$-torsors
over $H$ is fully faithful. Thus, for an $A$-torsor over $H$, being a central
extension is a property rather than an extra structure.
\item If $G$ is a connected commutative unipotent algebraic group over $k$,
then there always exists a dual pair $(G,G',\cE)$. Indeed, one can take $G'$ to
be any algebraic group over $k$ with $(G')^{perf}\cong(G^{perf})^*$, which
exists because $(G^{perf})^*\in\cC_k'$, and use a general fact (\cite{sga4},
expos\'e VIII): if $X$ is a scheme over $k$ and $A$ is a finite discrete
abelian group, the natural functor from the groupoid of $A$-torsors over $X$ to
the groupoid of $A$-torsors over $X^{perf}$ is an equivalence of categories.
\item The dual pair above is non-unique because there are many different choices for
$G'$: for instance, one can always replace $G'$ by $(G')^{(p)}$. However, once
$G'$ is chosen, the previous remark shows that the torsor $\cE$ is uniquely
determined.
\end{enumerate}
\end{rems}

\mbr

We now prove a general result.
\begin{prop}\label{p:Representative}
Let $G$ be an affine algebraic group over $k$, let $X$ be an affine scheme of
finite type over $k$, and let $Z=X^{perf}$. Given an action of $G^{perf}$ on
$Z$, there exists an affine scheme $Y$ of finite type over $k$, an action of
$G$ on $Y$, and an isomorphism $Y^{perf}\cong Z$ which induces the given action
of $G^{perf}$ on $Z$.
\end{prop}
\begin{proof}
Let us write $G=\Spec A$ and $X=\Spec B$, where $A$ and $B$ are finitely
generated $k$-algebras. Note that $A$ is a Hopf algebra. Moreover, the
$G^{perf}$-action on $Z$ amounts to a $k$-algebra homomorphism
$\nu:B^{perf}\rar{}A^{perf}\tens_k B^{perf}$ which makes $B^{perf}$ a comodule
over $A^{perf}$. It is well known and easy to check that any comodule over any
coalgebra is a filtered union of finite dimensional sub-comodules. In
particular, since $B$ is finitely generated over $k$, there exists a finite
collection of elements $b_1,\dotsc,b_r\in B^{perf}$ such that the image of $B$
in $B^{perf}$ is contained in $C:=k[b_1,\dotsc,b_r]$ and
\[
\nu(b_j) = \sum_{i=1}^N a_{ij}\tens b_i, \qquad 1\leq j\leq r, \qquad \text{for
some } a_{ij}\in A^{perf}.
\]
We have $A\subset A^{perf}$ (recall that our algebraic groups are reduced by
assumption), and by definition, there exists $m\in\bN$ such that
$a_{ij}^{p^m}\in A$ for all $1\leq i\leq N$, $1\leq j\leq r$. Thus $\nu$
restricts to a $k$-algebra homomorphism $\nu_1:C^{p^m}\rar{}A\tens_k C^{p^m}$.
Now put $Y=\Spec C^{p^m}$. Then $Y$ is an affine scheme of finite type over
$k$, and $\nu_1$ defines an action of $G$ on $Y$. In addition, since $C^{p^m}$
contains the image of $B^{p^m}$ in $B^{perf}$, it is clear that the inclusion
$C^{p^m}\into B^{perf}$ induces an isomorphism
$(C^{p^m})^{perf}\rar{\simeq}B^{perf}$, and the proof is complete.
\end{proof}

\begin{cor}\label{c:representative-dual}
Let $G$ be a connected unipotent group over $k$ whose nilpotence class is
$<p=\operatorname{char}k$, and let $\fg=\Log G$, as in
\S\ref{ss:lazard-unipotent}. Then there exist a dual pair $(\fg,\fg',\cE)$ of
unipotent $k$-groups and an action of $G$ on $\fg'$ by algebraic group
automorphisms which induces the canonical action of $G^{perf}$ on
$(\fg')^{perf}\cong(\fg^{perf})^*$. Moreover, $\cE$ has a natural
$G$-equivariant structure with respect to the induced action of $G$ on
$\fg'\times\fg$.
\end{cor}
\begin{proof}
The first statement is a special case of Proposition \ref{p:Representative}. The second one follows formally from the first one. Indeed, since $G^{perf}$ acts on $\fg^{perf}$ by group automorphisms, it follows from the universal property of $\cE^{perf}$ that it has a canonical $G^{perf}$-equivariant structure. This structure amounts to an isomorphism between $\act^*\cE^{perf}$ and $\pr^*\cE^{perf}$ satisfying a certain compatibility condition, where
$\act:G^{perf}\times(\fg')^{perf}\times\fg^{perf}\rar{} (\fg')^{perf}\times\fg^{perf}$ is the action morphism and $\pr:G^{perf}\times(\fg')^{perf}\times\fg^{perf}\rar{} (\fg')^{perf}\times\fg^{perf}$ is the projection. However, by the first statement of this corollary, $\act$ is the perfectization of a morphism
$G\times\fg'\times\fg\rar{}\fg'\times\fg$, and of course, $\pr$ is the
perfectization of the projection $G\times\fg'\times\fg\rar{}\fg'\times\fg$.
Using the fact stated in Remark \ref{r:dual-pair}(2), we see that $\cE$ has a
unique $G$-equivariant structure which induces the canonical
$G^{perf}$-equivariant structure on $\cE^{perf}$.
\end{proof}


\section{Fourier-Deligne transform}\label{a:fourier-deligne}

\subsection{Fourier transforms}\label{aa:fourier} The Fourier-Deligne transform was
introduced by P.~Deligne in \cite{deligne}. Its definition can
perhaps be motivated by the usual Fourier transform. If $f$ is an
$L^1$ function on $\bR$ with respect to the Lebesgue measure, its
Fourier transform is a function on ``a different copy of $\bR$''
defined by
\[
\widehat{f}(\xi) = \int_{-\infty}^\infty f(x) e^{-i\xi x}\, dx.
\]
Observe that $e^{-i\xi x}$ is a ``universal unitary character of
$\bR$'', in other words, the map $\xi\longmapsto\bigl(x\mapsto
e^{-i\xi x}\bigr)$ identifies $\bR$ with the Pontryagin dual of
$\bR$ (as a topological group). Moreover, the passage from $f$ to
$\widehat{f}$ can be thought of as consisting of three steps:
\begin{enumerate}[1)]
\item Consider $f$ as a function of two variables, $\xi$ and $x$,
i.e., pull $f$ back by the second projection
$p:\bR\times\bR\rar{}\bR$ (where the coordinate on the first copy
of $\bR$ is denoted by $\xi$ and that on the second copy is
denoted by $x$).
\item Multiply $p^*(f)$ by the ``universal character'' $e^{-i\xi
x}$.
\item Integrate the result along the fibers of the first
projection $p':\bR\times\bR\rar{}\bR$.
\end{enumerate}

\mbr

Each of these operations has an obvious analogue in the world of
$\ell$-adic complexes. The only problem is that ``integration
along the fibers'' has two analogues for (complexes of) sheaves:
the $!$-pushforward and the $*$-pushforward (see
\S\ref{aa:six-definitions}). Fortunately, these two choices give
equivalent definitions of the Fourier-Deligne transform. This is
one of its key properties (see \S\ref{aa:FD-properties} below),
and is the main result of \cite{saibi}.

\begin{defin}\label{d:fourier-deligne}
Let $k$ be a perfect field of characteristic $p>0$, let $\ell$ be a prime
different from $p$, and let $(G,G',\cE)$ be a dual pair of unipotent $k$-groups
as in Definition \ref{d:dual-pair}. Fix an injective homomorphism of abelian
groups $\ze:\bQ_p/\bZ_p\into\ql^\times$, and let $\cL$ denote the local system
on $G'\times G$ associated to $\cE$ via the homomorphism $\ze$. Consider the
projections
\[
G' \stackrel{\pr'}{\longleftarrow} G'\times G \rar{\pr} G.
\]
The \emph{Fourier-Deligne transform} associated to the data $(G,G',\cE,\ze)$ is
the functor
\[
\cF : \sD(G) \rar{} \sD(G')
\]
defined by
\[
\cF(M) = \pr'_! \bigl(\pr^*(M)\tens\cL\bigr)[\dim G].
\]
\end{defin}
\begin{rems}\label{r:fourier-deligne}
\begin{enumerate}[(1)]
\item Strictly speaking, we have not one Fourier-Deligne transform, but a
family of such transforms, parameterized by all possible embeddings
$\ze:\bQ_p/\bZ_p\into\ql^\times$. However, throughout this paper we fix $\ze$
once and for all, and so we speak of ``the'' Fourier-Deligne transform.
\item It was already known to Deligne \cite{deligne} that $\cF$ is an
equivalence of categories.
\item We get an isomorphic functor if we replace $\pr'_!$ with $\pr'_*$ in the
definition above: see Theorem \ref{t:fourier-supports}. The first proof of this
fact was published by Saibi \cite{saibi}.
\item The shift $[\dim G]$ appears in the definition above to ensure that $\cF$
takes perverse sheaves to perverse sheaves: see Corollary
\ref{c:fourier-perverse}.
\end{enumerate}
\end{rems}

\subsection{Perverse sheaves}\label{aa:perverse} Before discussing the properties of
the Fourier-Deligne transform, we briefly recall the definitions
and a few basic facts from the theory of perverse sheaves; we
refer the reader to \cite{bbd} for the details and for more
information.

\mbr

We only work with perverse sheaves for the \emph{self-dual} (or \emph{middle})
\emph{perversity}. Let $k$ be a field and $\ell$ a prime different from
$\operatorname{char}k$, as in Appendix \ref{a:equivariant-derived}. Given a
separated scheme $X$ of finite type over $k$, we define a $t$-structure on the
triangulated category $\sD(X)$, called the \emph{perverse $t$-structure}, as
follows. The full subcategory $^p\sD^{\leq 0}(X)\subseteq\sD(X)$ consists of
all complexes $K$ such that $\dim\Supp \cH^i(K)\leq -i$ for all $i\in\bZ$,
where $\cH^i(K)$ denotes the $i$-th cohomology sheaf of $K$, and
$\supp\cH^i(K)$ is its support, which is a constructible subset of $X$. In
particular, if $K\in{}^p\sD^{\leq 0}(X)$, then $\cH^i(K)=0$ for all $i>0$. We
also define
\[
^p\sD^{\geq 0}(X) = \bD_X\bigl({}^p\sD^{\leq 0}(X)\bigr), \qquad
\Perv(X) = {}^p\sD^{\geq 0}(X)\cap{}^p\sD^{\leq 0}(X).
\]
(This is not exactly the same as, but is equivalent to, the
definition given in \cite{bbd}.) The full subcategories
${}^p\sD^{\geq 0}(X)$ and ${}^p\sD^{\leq 0}(X)$ define a
$t$-structure on $X$ (see \emph{op.~cit.}, \S2.2), and $\Perv(X)$
is its heart. It is an abelian category whose objects are called
\emph{perverse sheaves} on $X$. By definition, $\Perv(X)$ is
stable under $\bD_X$. One knows that every perverse sheaf on $X$
has finite length (\emph{op.~cit.}, Th\'eor\`eme 4.3.1(i)). The
following result (\emph{op.~cit.}, Th\'eor\`eme 4.1.1) is
essentially due to M.~Artin:
\begin{thm}\label{t:artin}
If $f:X\rar{}Y$ is an affine morphism of separated schemes of
finite type over $k$, the functor $f_*:\sD(X)\rar{}\sD(Y)$ takes
${}^p\sD^{\leq 0}(X)$ into ${}^p\sD^{\leq 0}(Y)$.
\end{thm}
In view of the Verdier duality theorem (see
\S\ref{aa:six-formalism}), we immediately obtain
\begin{cor}[\emph{op.~cit.}, Corollaire 4.1.2]\label{c:artin}
Under the same assumptions, the functor $f_!:\sD(X)\rar{}\sD(Y)$
takes ${}^p\sD^{\geq 0}(X)$ into ${}^p\sD^{\geq 0}(Y)$.
\end{cor}

\mbr

\subsection{Properties of the Fourier-Deligne transform}\label{aa:FD-properties}
Until Proposition \ref{p:fourier-equivariant} we fix $(G,G',\cE,\ze)$ as in
Definition \ref{d:fourier-deligne}, and we will write $d=\dim G=\dim G'$.
\begin{thm}[P.~Deligne, \cite{deligne}]\label{t:fourier-equivalence}
The functor $\cF:\sD(G)\rar{}\sD(G')$ is an equivalence of categories, and a
quasi-inverse functor is given by $N\longmapsto (\iota^*\cF' N)(d)$, where
$(d)$ denotes the $d$-th Tate twist, $\iota:G\to G$ is the inversion map
$x\mapsto -x$, and $\cF':\sD(G')\rar{}\sD(G)$ is also a Fourier-Deligne
transform, defined by
\[
\cF'(N) = \pr_!\bigl( (\pr^{\prime*} N)\tens\cL \bigr)[d].
\]
\end{thm}
A published proof of this result is given in \cite{saibi}, Th\'eor\`eme
2.2.4.1.
\begin{thm}[M.~Saibi]\label{t:fourier-supports}
If $\cF_*:\sD(G)\rar{}\sD(G')$ is defined by $\cF_*(M)=\pr'_* \bigl( (\pr^*
M)\tens\cL \bigr)[d]$, then the functors $\cF_*$ and $\cF$ are naturally
isomorphic.
\end{thm}

This theorem was conjectured by Deligne, and is one of the main results of
\cite{saibi}. However, we will sketch a proof which is much shorter than the
one appearing in {\em op.~cit.}, following an argument explained to us by
Dennis Gaitsgory\footnote{We thank Dennis for allowing us to reproduce his proof. It appears that a similar idea was discovered independently by several other mathematicians.}. In view of Theorem \ref{t:fourier-equivalence}, to prove that
the functor $\cF_*$ is isomorphic to $\cF$, it suffices to show that $\cF_*$ is
right adjoint to the functor $N\longmapsto (\iota^*\cF' N)(d)$. However, this
can be done by a straightforward computation which is based on the fact that if $A$, $B$ are functors between arbitrary categories such that the composition $A\circ B$ is defined, and $A'$, $B'$ are functors right adjoint to $A$ and $B$, respectively, then $B'\circ A'$ is right adjoint to $A\circ B$.

\mbr

For the reader's convenience, we provide a detailed computation. Given $M\in\sD(G)$ and $N\in\sD(G')$, we have a chain of bifunctorial isomorphisms
\begin{eqnarray*}
\hm{(\iota^*\cF'N)(d),M} &\cong&
\hm{\pr_!((\pr^{\prime*}N)\tens\cL)[d],(\iota_*M)(-d)} \\
&\cong& \hm{((\pr^{\prime*}N)\tens\cL)[d],\pr^!(\iota_*M)(-d)} \\
&\cong& \hm{((\pr^{\prime*}N)\tens\cL)[d],\pr^*(\iota_*M)[2d]} \\
&\cong& \hm{(\pr^{\prime*}N)\tens\cL,\pr^*(\iota_*M)[d]} \\
&\cong& \hm{(\pr^{\prime*}N)\tens\cL,(\id_{G'}\times\iota)_*\pr^*(M)[d]} \\
&\cong& \hm{(\id_{G'}\times\iota)^*((\pr^{\prime*}N)\tens\cL),\pr^*(M)[d]} \\
&\cong& \hm{(\id_{G'}\times\iota)^*(\pr^{\prime*}N)\tens(\id_{G'}\times\iota)^*\cL,\pr^*(M)[d]} \\
&\cong& \hm{(\pr^{\prime*}N)\tens\cL^{-1},\pr^*(M)[d]} \\
&\cong& \hm{\pr^{\prime*}N,(\cL\tens\pr^*(M))[d]} \\
&\cong& \hm{N,\pr'_* ( (\pr^* M)\tens\cL )[d]}.
\end{eqnarray*}
The justifications of these isomorphisms are mostly trivial. The third one
follows from the fact that $\pr$ is a smooth morphism of relative dimension
$d$. The seventh one, which amounts to
$(\id_{G'}\times\iota)^*\cL\cong\cL^{-1}$, follows easily from the universal
property of the $\bQ_p/\bZ_p$-torsor $\cE$ over $G'\times G$. This completes
the proof of Theorem \ref{t:fourier-supports}.

\begin{cor}[Fourier-Deligne transform commutes with duality]\label{c:fourier-duality}
There exist functorial isomorphisms
\[
\cF\bigl(\bD_G^- M\bigr) \cong \bigl(\bD_{G'}\cF(M)\bigr)(-d)
\]
for all $M\in\sD(G)$, where $\bD_G^-=\iota^*\circ\bD_G=\bD_G\circ\iota^*$, as
in Conjecture \ref{conj:1} of \S\ref{ss:conjectures}, $\iota:G\rar{}G$ is the inversion map, and $\bD_G$
is the Verdier duality functor, defined in \S\ref{aa:verdier}.
\end{cor}
\begin{proof}
This is a straightforward exercise using Theorem \ref{t:fourier-supports}.
\end{proof}
\begin{cor}[Fourier-Deligne transform preserves perversity]\label{c:fourier-perverse}
The functor \\ $\cF:\sD(G)\rar{}\sD(G')$ takes perverse sheaves on $G$ to
perverse sheaves on $G'$.
\end{cor}
\begin{proof}
It is well known and easy to see that the functors $\pr^*[d]$ and $-\tens\cL$
take perverse sheaves to perverse sheaves. Moreover, since $\pr'$ is an affine
morphism, it follows from Theorem \ref{t:artin} that $\cF_*$ takes $\Perv(G)$
into $^p\sD^{\leq 0}(G')$, and it follows from Corollary \ref{c:artin} that
$\cF$ takes $\Perv(G)$ into $^p\sD^{\geq 0}(G')$. But $\cF\cong\cF_*$ by
Theorem \ref{t:fourier-supports}, so $\cF$ takes $\Perv(G)$ into $^p\sD^{\leq
0}(G')\cap{}^p\sD^{\geq 0}(G')=\Perv(G')$, as desired.
\end{proof}

\mbr

Another fact about the Fourier-Deligne transform that we need is given in Proposition \ref{p:fourier-convolution} below. It is independent of the previous theorems; indeed, it only uses
the following property of the torsor $\cE$:
\begin{lem}\label{l:property-E}
Let $\al:G\times G\to G$ denote the group operation in $G$, and let
$p_{12},p_{13}:G'\times G\times G\to G'\times G$ be the projections
$(g',g_1,g_2)\mapsto(g',g_1)$ and $(g',g_1,g_2)\mapsto(g',g_2)$, respectively.
Then there is a canonical isomorphism $(\id_{G'}\times\al)^*\cE\cong
p_{12}^*(\cE)\otimes p_{23}^*(\cE)$ of $\bQ_p/\bZ_p$-torsors on $G'\times
G\times G$. If $\iota:G\to G$ denotes the inversion map, there is also a
canonical isomorphism $(\id_{G'}\times\iota)^*\cE\cong\cE^\vee$, where
$\cE^\vee$ denotes the dual $\bQ_p/\bZ_p$-torsor\footnote{The statement of the
lemma uses the fact that if $A$ is a discrete abelian group, then the category
of $A$-torsors over any scheme $X$ is a rigid monoidal category.}.
\end{lem}
\begin{proof}
As before, we may replace $G$ with $G^{perf}$, $G'$ with $G^*=(G^{perf})^*$,
and $\cE$ with the universal torsor $\cE$ on $G^*\times G$. By definition,
$\cE$ is then the universal central extension of $G$ by $\bQ_p/\bZ_p$, i.e., if
we think of $G^*\times G$ as a group scheme over $G^*$ (where the group
operation on $G^*$ is ignored), then $\cE$ is a central extension of $G^*\times
G$ by $\bQ_p/\bZ_p$ which represents the universal class in $\Ext^1(G^*\times
G,\bQ_p/\bZ_p)$. But now the statement is obvious, because if we put
$H=G^*\times G$, then we have a natural identification $H\times_{G^*}H\cong
G^*\times G\times G$, so that the group operation in $H$ is identified with
$(\id_{G^*}\times\al):H\times_{G^*}H\to H$, the two projections
$H\times_{G^*}H\to H$ are identified with $p_{12}$ and $p_{13}$, and the
inversion map $H\to H$ is identified with $(\id_{G^*}\times\iota)$.
\end{proof}

\begin{prop}\label{p:fourier-convolution}
The Fourier-Deligne transform $\cF:\sD(G)\rar{}\sD(G')$ can be naturally
upgraded to an equivalence of symmetric monoidal categories, where $\sD(G)$ is
equipped with the monoidal structure given by the convolution defined in
\S\ref{aa:convolution} and $\sD(G')$ is equipped with the monoidal structure
given by $M\odot N=(M\tens N)[-d]$.
\end{prop}
\begin{proof}
We only check that there are natural bifunctorial isomorphisms
$\cF(M*N)\cong(\cF(M)\tens\cF(N))[-d]$ for all $M,N\in\sD(G)$; the rest
(compatibility with units and associativity and commutativity constraints) is
straightforward. First we set up some notation. Recall the projections
\[
G' \xleftarrow{\ \pr'\ } G'\times G \xrightarrow{\ \pr\ } G \qquad \text{and}
\qquad p_1,p_2:G\times G\rar{} G.
\]
As in Lemma \ref{l:property-E}, we will write $\al:G\times G\to G$ for the
group operations and $p_{12},p_{13}:G'\times G\times G\to G'\times G$ for the
projections. We also introduce the projections
\[
\ti{\pr}:G'\times G\times G \rar{} G\times G  \qquad \text{(along the first
factor)},
\]
\[
\ti{p_1},\ti{p_2}: G'\times G\times G \rar{} G \qquad \text{(onto the second
and third factors)},
\]
\[
\text{and} \qquad \ti{\pr}':G'\times G\times G \rar{} G' \qquad \text{(onto the
first factor)}.
\]
Let us now fix $M,N\in\sD(G)$. By definition,
\[
\cF(M*N) = \pr'_!(\pr^*(M*N)\tens\cL)[d] = \pr'_!\bigl( \pr^*\al_!(p_1^*M\tens
p_2^*N)\tens\cL \bigr)[d].
\]
We have a cartesian diagram
\[
\xymatrix{
  G'\times G\times G \ar[rr]^{\id_{G'}\times\al} \ar[d]_{\ti{\pr}}  &  & G'\times
  G \ar[d]^{\pr} \\
  G\times G \ar[rr]^{\al}  &  &  G
   }
\]
Hence the proper base change theorem (\S\ref{aa:six-formalism}) implies that
there is a natural isomorphism
\[
\pr'_!\bigl( \pr^*\al_!(p_1^*M\tens p_2^*N)\tens\cL \bigr)[d] \cong
\pr'_!\bigl( ( (\id_{G'}\times\al)_! \ti{\pr}^*(p_1^*M\tens p_2^*N) ) \tens \cL
\bigr) [d].
\]
Since $p_1\circ\ti{\pr}=\ti{p_1}$ and $p_2\circ\ti{\pr}=\ti{p_2}$, we have
\[
\ti{\pr}^*(p_1^*M\tens p_2^*N) \cong \ti{p_1}^*M \tens \ti{p_2}^* N.
\]
Therefore, by the projection formula (\S\ref{aa:six-formalism}),
\[
\bigl( (\id_{G'}\times\al)_! \ti{\pr}^*(p_1^*M\tens p_2^*N) \bigr) \tens \cL
\cong (\id_{G'}\times\al)_! \bigl( \ti{p_1}^*M \tens \ti{p_2}^* N \tens
(\id_{G'}\times\al)^*\cL \bigr).
\]
But the first statement of Lemma \ref{l:property-E} implies that
$(\id_{G'}\times\al)^*\cL\cong p_{12}^*\cL\tens p_{13}^*\cL$. Putting the
previous computations together and using the equality
$\pr'\circ(\id_{G'}\times\al)=\ti{\pr}'$, we obtain the first isomorphism in
the following chain:
\begin{eqnarray*}
\cF(M*N) &\cong& \ti{\pr}'_! \bigl( \ti{p_1}^*M \tens \ti{p_2}^*N \tens
p_{12}^*\cL \tens p_{13}^*\cL \bigr) [d] \\
&\cong& \ti{\pr}'_! \bigl( p_{12}^*(\pr^*M\tens\cL) \tens p_{13}^*(\pr^*N \tens
\cL) \bigr) [d] \\
&\cong& \pr'_! p_{12!} \bigl( p_{12}^*(\pr^*M\tens\cL) \tens
p_{13}^*(\pr^*N\tens\cL) \bigr) [d] \\
&\cong& \pr'_! \bigl( (\pr^*M\tens\cL) \tens p_{12!}p_{13}^* (\pr^*N\tens\cL)
\bigr) [d] \\
&\cong& \pr'_! \bigl( (\pr^*M\tens\cL) \tens \pr^{\prime *} \pr'_!
(\pr^*N\tens\cL) \bigr)[d] \\
&\cong& \pr'_!(\pr^*M\tens\cL) \tens \pr'_!(\pr^*N\tens\cL)[d] \\
&\cong& (\cF(M)\tens\cF(N))[-d].
\end{eqnarray*}
In this chain the second isomorphism follows from $\ti{p_1}=\pr\circ p_{12}$
and $\ti{p_2}=\pr\circ p_{13}$, the third one follows from $\ti{\pr}'=\pr'\circ
p_{12}$, the fourth and sixth isomorphisms follow from the projection formula
(\S\ref{aa:six-formalism}), the seventh one follows from the definition, and
the fifth one follows from the proper base change theorem
(\S\ref{aa:six-formalism}) applied to the cartesian diagram
\[
\xymatrix{
  G'\times G\times G \ar[rr]^{p_{13}} \ar[d]_{p_{12}}  &  & G'\times
  G \ar[d]^{\pr'} \\
  G'\times G \ar[rr]^{\pr'}  &  &  G
   }
\]
This completes the proof.
\end{proof}

\mbr

Finally, we will need to know that the Fourier-Deligne transform works in the
equivariant setting, i.e., it can be defined on the categories $\sD_G(\fg)$
used in Section \ref{s:char-sheaves-orbmethod}.
\begin{prop}\label{p:fourier-equivariant}
Let $G$ be a connected unipotent group over $k$ of nilpotence class $<p$, let
$\fg=\Log G$, let $(\fg,\fg',\cE)$ be a dual pair of unipotent $k$-groups
satisfying the conclusion of Corollary \ref{c:representative-dual}, and let
$\ze:\bQ_p/\bZ_p\into\ql^\times$ be as before. The Fourier-Deligne transform
$\cF:\sD(\fg)\rar{}\sD(\fg')$ defined in \S\ref{aa:fourier} can be naturally
lifted to an equivalence of monoidal categories
$\cF:\sD_G(\fg)\rar{\sim}\sD_G(\fg')$, where $\sD_G(\fg)$ is equipped with
convolution $($cf. \S\ref{aa:categories-equivariant}$)$ and $\sD_G(\fg')$ is
equipped with the monoidal structure defined by $M\odot N=(M\tens N)[-\dim G]$.
\end{prop}
This result is straightforward. One defines the lift of $\cF$ by the same
formula as in Definition \ref{d:fourier-deligne}, using the fact that the local
system $\cL$ on $\fg'\times\fg$ has a $G$-equivariant structure coming from
that on $\cE$ together with the remarks in \S\ref{aa:categories-equivariant}.
Then one easily verifies that the proofs of all the results in this subsection
go through in the equivariant setting by checking that each step is compatible
with the $G$-equivariant structures. The details are left to the reader.


\section{Some counterexamples related to the orbit
method}\label{a:counterexamples}

Let us briefly fix our conventions and describe some of the counterexamples
presented in this appendix. We use the setup and notation of Section
\ref{s:orbmethod}. In particular, $\Ga$ will denote a $p$-group of nilpotence
class $<p$ and $\fg=\Log (\Ga)$ will denote the corresponding Lie ring. Recall
from \S\ref{ss:fourier} that the composition
\begin{equation}\label{e:compos}
\Phi:\bC \Ga \cong \Meas(\Ga) \xrar{\log_*} \Meas(\fg) \rar{\cF} \Fun(\fg^*)
\end{equation}
restricts to an algebra isomorphism between the center $\Z \Ga$ of $\bC \Ga$
and the algebra $\Fun(\fg^*)^\Ga$ of $\Ga$-invariant functions on $\fg^*$. One
of the goals of this appendix is to show that the composition \eqref{e:compos}
{\em is not an isomorphism of $\Z\Ga$-modules}, where the $\Z\Ga$-module
structure on $\Fun(\fg^*)$ is induced by the isomorphism
$\Z\Ga\rar{\simeq}\Fun(\fg^*)^\Ga$.

\mbr

If $f\in\fg^*$ and $\Om\subset\fg^*$ is the $\Ga$-orbit of $f$, we denote by
$\rho_\Om$ the irreducible representation of $\Ga$ associated to $\Om$ via the
orbit method. We write $\Ga^f$ for the stabilizer of $f$ in $\Ga$ and $\fg^f$
for the Lie subring of $\fg$ corresponding to $\Ga^f$; one can check (cf.
\cite{kirillov,orbmethod}) that $\fg^f$ coincides with the kernel of the
natural alternating bilinear form $B_f:\fg\times\fg\to\cst$ associated to $f$,
so the restriction of $f$ to $\fg^f$ is a Lie homomorphism\footnote{This also
follows from the $\Ga^f$-invariance of $f:\fg^f\to\bC^{\times}$.}
$\fg^f\to\bC^{\times}$. Therefore it defines a group homomorphism
$\Gamma_f\to\bC^{\times}$. We present a counterexample below showing that in
general {\em the representation of $\Gamma$ induced by this group homomorphism
is not isomorphic to a multiple of the representation $\rho_\Om$}.

\subsection{Statements to disprove}\label{aa:assertions}
Our goal is to show that each of the following statements is false in general:

\begin{enumerate}[(1)]
\item The representation of $\Ga$ obtained from inducing the
$1$-dimensional representation of $\Ga^f$ defined by $f$ is isomorphic to a
multiple of $\rho_\Om$.
\item Every Lie subring
$\fh\subseteq\fg$ satisfying $f([\fh,\fh])=\{1\}$ is contained in a
polarization of $\fg$ at $f$.
\item The composition \eqref{e:compos} is an isomorphism of $\Z \Ga$-modules,
where $\Fun(\fg^*)$ is a $\Z \Ga$-module via the algebra isomorphism $\Z
\Ga\rar{\simeq}\Fun(\fg^*)^\Ga$ mentioned above. In other words,
$\Phi(xy)=\Phi(x)\Phi(y)$ whenever $x\in\Z \Ga$ and $y\in\bC \Ga$.
\end{enumerate}

\sbr

\noindent
On the other hand, it is a simple exercise to prove that all three statements
hold when the nilpotence class of $\fg$ is at most $2$.

\subsection{A counterexample to (1) and (2)}\label{aa:1-and-2}
Let $\chi_f$ denote the $1$-dimensional character of $\Ga^f$ defined by $f$. It
follows from the orbit method \cite{kirillov,orbmethod} that the irreducible
summands of the induced representation $\Ind_{\Ga^f}^\Ga(\chi_f)$  correspond
to those coadjoint orbits in $\fg^*$ that intersect $f+(\fg^f)^\perp$
nontrivially. Thus, (1) is equivalent to $f+(\fg^f)^\perp\subseteq\Om$, which
by counting is equivalent to $f+(\fg^f)^\perp=\Om$. To find a situation where
this fails, let $q$ be odd, and let $\fg$ be the Lie algebra over $\bF_q$ with
basis $x,y,z,t$ and commutation relations
\[
[x,y]=z, \quad [x,z]=t, \quad [y,z]=0,\quad [t,\fg]=(0).
\]
Define $f:\fg\to\cst$ by $f=\psi\circ p$, where $p:\fg\to\bF_q$ is defined by
$p(x)=p(y)=p(z)=0$, $p(t)=1$, and $\psi:\bF_q\to\cst$ is a nontrivial additive
character of $\bF_q$. One can easily check that $f+(\fg^f)^\perp\neq\Om$ in
this case: for instance, $\exp(\ad x)(y)=y+t/2$, so that $\Ad^*(e^x)(f)-f$ does
not vanish on $y$, whereas $y\in\fg^f$.

\mbr

This construction also provides a counterexample to statement (2) of
\S\ref{aa:assertions}. Namely, $\fh:=\Span(x,t)$ is an abelian subalgebra of
$\fg$, and we claim that $\fh$ is not contained in any polarization of $\fg$ at
$f$. To prove this, notice that any polarization at $f$ must contain $\fg^f$,
whereas $\fh$ and $\fg^f$ together generate $\fg$ as a Lie algebra, because
$y\in\fg^f$. Since $f([\fg,\fg])\neq\{1\}$, we see that $\fh$ cannot be
contained in any polarization of $\fg$ at $f$.

\subsection{Properties equivalent to (3)}\label{aa:equiv}
Recall that our property (3) says that the map $\Phi$ defined by
\eqref{e:compos} is an isomorphism of $\Z \Ga$-modules, i.e.,
$\Phi(xy)=\Phi(x)\Phi(y)$ whenever $x\in\Z \Ga$ and $y\in\bC \Ga$.  We will
formulate some properties equivalent to (3), and then we will show that they do
not hold in general.

\mbr

For every orbit  $\Om\subset\fg^*$ let $I(\Omega)$ be the kernel
of the corresponding representation of $\bC \Ga$; so  $I(\Omega)$
is a two-sided ideal of $\bC \Ga$. Let
$J(\Omega)\subset\Fun(\fg^*)$ be the ideal of functions that
vanish on $\Om$ and $K(\Omega)\subset\Fun(\fg^*)$ the ideal of
functions supported on $\Om$. Remark \ref{r:idemp} implies that
(3) is equivalent to each of the following properties:

\sbr

\begin{enumerate}
\item[(4)] for every orbit $\Om$ the composition \eqref{e:compos} takes $I(\Omega)$ onto $J(\Omega)$;
\item[(5)] $\Phi^{-1}(J(\Omega))$ is a left ideal of $\bC \Ga$ for every orbit $\Om$;
\item[(6)] $\Phi^{-1}(K(\Omega))$ is a left ideal of $\bC \Ga$ for every orbit $\Om$.
\end{enumerate}

\sbr

\noindent
Note that since $\Phi^{-1}(J(\Omega))$ and $\Phi^{-1}(K(\Omega))$ are stable
under $\Ga$-conjugation, one can replace the word  ``left" in (5) and (6) by
``right" or ``two-sided".

\mbr

Let $\Perm(\Om)$ denote the permutation representation of $\Ga$ associated to
the $\Ga$-action on $\Om$. Since $\Phi$ commutes with $\Ga$-conjugation
property (4) implies that $\bC \Ga /I(\Omega)$ and $\Fun(\fg^*)/J(\Omega)$ are
isomorphic as $\Ga$-modules, which is equivalent to the following property:

\mbr

(7) The $\Ga$-module $\Perm(\Om)$ is isomorphic to $\rho_\Om\tens\rho_\Om^*$.

\mbr

So a counterexample to (7) would also be a counterexample to (3)-(6). In the
next subsection we will construct a counterexample to (7) in which $\fg$ has
nilpotence class $4$. Then we will show that (7) holds whenever the nilpotence
class of $\fg$ is at most $3$, while (3)-(6) can already fail when $\fg$ has
nilpotence class $3$.

\subsection{Counterexample to (7) }\label{aa:counterexample-3}
This subsection is close in spirit to a work by F.~du Cloux (see
\S\ref{aa:Cloux}). We will only consider the case where $\fg$ is a
Lie algebra over $\bF_q$, rather than an arbitrary Lie ring. In
this case $\fg^*$ will be identified with the underlying additive
group of the dual vector space $\Hom_{\bF_q}(\fg,\bF_q)$. (The
identification depends on a choice of a nontrivial additive
character $\bF_q\to\cst$. This choice is irrelevant in what
follows.) We consider a coadjoint orbit $\Om\subset\fg^*$ as a
closed subvariety. In particular, the tangent bundle $T\Om$ is
defined. However, we do not need $T\Om$ as a variety, but only as
a set, namely, the set of pairs $(f,\la)$, where $f\in\Om$ and
$\la\in T_f\Om\cong\fg/\fg^f$; the latter will be thought of as a
subspace of $\fg^*$, namely, the image of the linear map
$\fg\to\fg^*$ induced by the form $B_f$.

\mbr

We want to compare the character of $\rho_\Om\tens\rho_\Om^*$ and that of
$\Perm(\Om)$. The former is the Fourier transform of the measure $\mu$ on
$\fg^*$ defined by $\mu=(1/\abs{\Om})\cdot(\de_\Om*\de_{-\Om})$. Equivalently,
$\mu(\la)=\abs{\pi^{-1}(\la)}/\abs{\Om}$, where $\la\in\fg^*$ and $\pi
=\pi_{\Omega}:\Om\times\Om\to\fg^*$ is the subtraction map $(f,g)\mapsto f-g$.
The latter is the Fourier transform of the measure $\mut$ on $\fg^*$ defined by
$\mut(\la)=\abs{\pit^{-1}(\la)}/\abs{\Om}$, where $\pit
=\pit_{\Omega}:T\Om\to\fg^*$ is the natural map, $(f,\la)\mapsto\la$. (Note
that $\pit_{\Omega}$ is a ``degeneration'' of $\pi_{\Omega}$ because one can
think of a point of $T\Om$ as a pair of infinitely close points of $\Om$.)

\mbr

Thus $\rho_\Om\tens\rho_\Om^*\cong\Perm(\Om)$ if and only if
\begin{equation}\label{e:eqcard}
\abs{\pi_{\Omega}^{-1}(\la)}=\abs{\pit_{\Omega}^{-1}(\la)} \mbox{ for all }
\la\in\fg^*.
\end{equation}
Observe now that the maps $\pi$ and $\pit$ make sense for any smooth closed
subvariety $Y$ of an affine space $A$ over $\bF_q$ (they take values in the
vector space $V$ corresponding to $A$), and we can ask the same question: is it
true that $\abs{\pi^{-1}(\la)}=\abs{\pit^{-1}(\la)}$ for all $\la\in V$? If $Y$
is the affine space itself the answer is ``yes".

\begin{exer} \label{ex:parab}
If $q$ is odd and $Y$ is a parabola in the affine plane, the answer is also
``yes".
\end{exer}

Usually the answer is no. E.g., the answer is negative for the following curve
$Y$.

\begin{exer}\label{ex:Veronese}
Let  $Y:=\{(t,t^2,t^3)\in\bA^3\,\big\vert\, t\in\bF_q\}$, and suppose that
$\operatorname{char}(\bF_q)\geq 5$. Then the maps $\pi$ and $\pit$ have
different images.
\end{exer}

Now we can give a counterexample to (7). We assume that
$\operatorname{char}\bF_q\geq 5$. Let $\fg$ be a semidirect product of a
$1$-dimensional subalgebra $\fc\subset\fg$ generated by an element $v\in\fc$
and a $4$-dimensional abelian ideal $\fa\subset\fg$ so that $(\ad v)^4=0$ and
$(\ad v)^3\ne 0$ (i.e., $\ad v:\fa\to\fa$ is a generic nilpotent operator). Let
$\Om'\subset\fa^*$ be a generic $C$-orbit, $C:=\Exp\fc$; here ``generic" means
that $\Om'$ is not contained in the kernel of $(\ad^*v)^3:\fa^*\to\fa^*$. Let
$\Om$ be the preimage of $\Om'$ in $\fg^*$; it is easy to see that $\Om$ is a
$\Ga$-orbit.
\begin{exer}\label{ex:Veronese2}
Equality $\eqref{e:eqcard}$ does not hold for $\Om$.
\end{exer}

\noindent Hint: $\Om'$ can be identified with the curve from Exercise \ref{ex:Veronese}.

\medskip

In the construction above it is {\em not} enough to take $\dim\fa=3$ (in this
case $\Om'$ can be identified with the parabola from Exercise \ref{ex:parab}).
More generally, we have the following
\begin{thm}
If $\fg$ has nilpotence class $\leq 3$ then $\fg$ satisfies property $(7)$,
i.e., $\rho_\Om\tens\rho_\Om^*\cong\Perm(\Om)$ for every coadjoint orbit
$\Om\subset\fg^*$.
\end{thm}
\begin{proof}
Let $\ad^*: \fg\to\End\fg$ be the coadjoint action, i.e., $\ad^* x:= -(\ad
x)^*$. Recall that for $f\in\fg^*$, we denote by $\fg^f$ its stabilizer in
$\fg$. Now fix $f$, let $\Om$ denote its $\Ga$-orbit, and let $\pi$, $\pit$ be
as above. It is easy to see that $\abs{\fg^f}^2\cdot\abs{\pi^{-1}(\la)}$ is the
number of solutions of the equation
\begin{equation}\label{e:d1}
\bigl( e^{\ad^* x} - e^{\ad^* y} \bigr)(f)=\la, \qquad x,y\in\fg,
\end{equation}
and $\abs{\fg^f}^2\cdot\abs{\pit^{-1}(\la)}$ is the number of solutions of the
equation
\begin{equation}\label{e:d2}
e^{\ad^* z} \bigl(\ad^* u (f)\bigr) = \la, \qquad z,u\in\fg.
\end{equation}
Thus it suffices to establish a bijection between the sets of solutions of
these two equations. We leave it to the reader to check that the change of
variables $x=\log(e^z e^{u/2})$, $y=\log(e^z e^{-u/2})$ transforms \eqref{e:d2}
into \eqref{e:d1}; here we are using the equality $e^{\ad^* z}\ad^* u =
e^{\ad^* z} \bigl( e^{\ad^* u /2} - e^{-\ad^* u/2}\bigr)$, which holds because
$(\ad^* u)^3=0$. This change of variables is invertible: $u=\log(e^{-y}
e^{-x})$, $z=\log(e^y e^{u/2})=\log\bigl(e^y (e^{-y} e^x)^{1/2}\bigr)$.
\end{proof}

\subsection{Counterexample to (3)--(6) of nilpotence class 3.}\label{aa:counterexample-2}
Statement (3) was formulated in \S\ref{aa:assertions} and the equivalent
statements (4)--(6) were formulated in \S\ref{aa:equiv}. We will show that they
already fail in a counterexample similar to the one given in
\S\ref{aa:counterexample-3}, but with $\dim\fa=3$ instead of $\dim\fa=4$.

\mbr

First of all, property (6) can be reformulated as follows:

\sbr

\begin{enumerate}
\item[(6$'$)] for every $\la\in\fg^*$ and $\ga\in\Ga$ the function
$$x\mapsto\la (\log (\ga e^x)), \qquad x\in\fg$$
is a linear combination of the characters from the $\Ga$-orbit $\Om_{\la}\ni
\la$.
\end{enumerate}

\sbr

\noindent
Now suppose that

\sbr

\begin{enumerate}[(i)]
\item
$\fg$ is a semidirect product of a $1$-dimensional subalgebra $\fc\subset\fg$
generated by an element $v\in\fc$ and an abelian ideal $\fa\subset\fg$,
\item
the restriction of $\la$ to $(\ad v)^2 (\fa )$ is nontrivial.
\end{enumerate}

\sbr

\noindent
Then (6$'$) cannot hold. For if it does, then the function
\begin{equation}\label{e:laCH1}
x\mapsto\la (\log (e^v e^x)), \qquad x\in\fa
\end{equation}
must be a linear combination of characters of $\fa$ of the form
\begin{equation}\label{e:laCH2}
x\mapsto\la (e^{t\ad v}x), \qquad t\in\bF_q\, .
\end{equation}
On the other hand, since $\fa$ is abelian the Campbell-Hausdorff formula tells
us that the function \eqref{e:laCH1} is proportional to the character
\begin{equation}\label{e:laCH3}
\la (x+\sum_{i>0} c_i (\ad v)^i x), \qquad x\in\fa ,
\end{equation}
where $c_i\in\bZ \bigl[\frac{1}{(p-1)!}\bigr]$ are certain universal coefficients; in particular,
$c_1=1/2$ and $c_2=1/12$. Since $c_2\ne c_1^2/2$, the character \eqref{e:laCH3}
does not have the form \eqref{e:laCH2}.

\subsection{Du Cloux's theorem} \label{aa:Cloux}
The results of \S\ref{aa:counterexample-3} above are similar in spirit to a
theorem by F.~du Cloux \cite{du-cloux}. Its statement is as follows. Let $\fg$
be a finite dimensional nilpotent Lie algebra over a field $k$ of
characteristic $0$, let $G=\Exp\fg$ be the corresponding unipotent group over
$k$, let $\fg^*=\Hom_k(\fg,k)$, let $S(\fg)$ denote the symmetric algebra of
$\fg$ (viewed as the algebra of polynomial functions on $\fg^*$), and let
$U(\fg)$ denote the universal enveloping algebra of $\fg$. Thus $S(\fg)$ is the
analogue of the commutative algebra $\Fun(\fg^*)$ considered above, and
$U(\fg)$ is the analogue of the group algebra $\bC\Ga$. If $\Om\subset\fg^*$ is
a $G$-orbit, let $J(\Om)$ denote the ideal of functions in $S(\fg)$ that vanish
on $\Om$, and let $I(\Om)\subset U(\fg)$ denote the kernel of the
representation of $U(\fg)$ associated to $\Om$ by Kirillov's theory. Assume
moreover that $\dim\Om=2$. Then {\it $S(\fg)/J$ and $U(\fg)/I$ are isomorphic
as $\fg$-modules with respect to the adjoint action if and only if $\Om$ has
degree $\leq 2$ as a subvariety of $\fg^*$.}



\begin{thebibliography}{55}

\bibitem[And95]{andre-1} C.A.M.~Andr\'e, {\em Basic characters of
the unitriangular group}, J. of Algebra \textbf{175} (1995),
287--319.

\bibitem[And98]{andre} C.A.M.~Andr\'e, {\em Irreducible characters
of finite algebra groups}, in: ``Matrices and Group Representations, Coimbra,
1998'', Textos Mat. S\'er B \textbf{19}, Univ. Coimbra, Coimbra, 1999,
pp.~65--80, arXiv: {\tt math.RT/9811132}

\bibitem[And02]{andre-3} C.A.M.~Andr\'e, {\em Basic characters of
the unitriangular group $($for arbitrary primes$)$}, Proc. Amer.
Math. Soc. \textbf{130} (2002), 739--765.

\bibitem[Ba86]{basmanova} S.A.~Basmanova, {\em The relation
between irreducible representations of a finite $p$-group and the Lie algebra
associated with it}, Vestnik Moskov. Univ. Ser. I Mat. Mekh. 1986, no. 5,
58--60, 100--101.

\bibitem[Be80]{begueri} L.~Begueri, {\em Dualit\'e sur un corps local \`a corps
r\'esiduel alg\'ebriquement clos}, M\'em. Soc. Math. France (N.S.) 1980/81,
no.~4.

\bibitem[BBD82]{bbd} A.A.~Beilinson, J.~Bernstein and P.~Deligne,
\emph{Faisceaux Pervers}, in: ``Analyse et topologie sur les
espaces singuliers (I)'', Ast\'erisque \textbf{100}, 1982.

\bibitem[BCD72]{solvable} P.~Bernat, C.~Conze, M.~Duflo,
N.~L\'evy-Nahas, M.~Rais, P.~Renouard and M.~Vergne, ``Repr\'esentations des
Groupes de Lie R\'esolubles.'' Monographies de la Soc. Math. de France
\textbf{4} (1972).

\bibitem[BL94]{ber-lunts} J.~Bernstein and V.~Lunts, ``Equivariant
sheaves and functors'', Lect. Notes in Math. \textbf{1578},
Springer-Verlag, Berlin, 1994.

\bibitem[Bo06a]{base-change} M.~Boyarchenko,
{\em Base change maps for unipotent algebra groups}, Preprint, January 2006,
arXiv: {\tt math.RT/0601133}


\bibitem[Bo06b]{charsheaves-orbmethod} M.~Boyarchenko, {\em
Character sheaves and the orbit method}, notes for a talk at the
AHA, LP, CFT-MM conference (CIRM, Luminy, France), June 2006,
available upon request.

\bibitem[BDx]{idempotents} M.~Boyarchenko and V.~Drinfeld, {\em
Idempotents in monoidal categories}, in preparation.

\bibitem[BSx]{orbmethod} M.~Boyarchenko and M.~Sabitova, {\em The
orbit method for $p$-groups and pro-$p$-groups}, in preparation.

\bibitem[De76]{deligne} P.~Deligne, Letter to D.~Kazhdan, November
29, 1976 (unpublished).

\bibitem[De80]{deligne-weil} P.~Deligne, {\em La conjecture de
Weil II}, Publ. Math. IHES \textbf{52} (1980), 137--252.

\bibitem[DI05]{diaconis-isaacs} P.~Diaconis and I.M.~Isaacs,
{\em Supercharacters and superclasses for algebra groups}, Preprint, 2005,
available at {\tt
http://www-stat.stanford.edu/$\sim$cgates/PERSI/papers/supercharacters.pdf}

\bibitem[dC81]{du-cloux} F.~du Cloux, {\em Non isomorphisme entre $U({\fg})/I$
et $S({\fg})/J$}, C. R. Acad. Sci. Paris S\'er. I Math. \textbf{293} (1981),
no.~1, 5--8.

\bibitem[Ek90]{ekedahl} T.~Ekedahl, {\em On the adic formalism},
in: ``The Grothendieck Festschrift, Vol. II'', 197--218, Progr.
Math. \textbf{87}, Birkh\"auser Boston, Boston, MA, 1990.

\bibitem[Gr65]{greenberg} M.J.~Greenberg, {\em Perfect closures of rings and schemes},
Proc. AMS \textbf{16} (1965), 313--317.

\bibitem[Gu73]{gutkin} E.A.~Gutkin, {\em Representations of
algebraic unipotent groups over a self-dual field}, Funkts. Analiz i Ego
Prilozheniya \textbf{7} (1973), 80.

\bibitem[Ha04]{halasi} Z.~Halasi, {\em On the characters and
commutators of finite algebra groups}, Jour. of Algebra \textbf{275} (2004),
481--487.

\bibitem[Isa95]{isaacs} I.M.~Issacs, {\em Characters of groups
associated with finite algebras}, Jour. of Algebra \textbf{177} (1995),
708--730.

\bibitem[Ka77]{kazhdan} D.~Kazhdan, {\em Proof of Springer's hypothesis},
Israel J. Math. \textbf{28} (1977), no.~4, 272--286.

\bibitem[Khu98]{kh} E.I.~Khukhro, ``$p$-Automorphisms of Finite
$p$-groups'', Lond. Math. Soc. Lect. Note Series \textbf{246}, Cambridge
University Press, 1998.

\bibitem[Ki62]{kirillov} A.A.~Kirillov, {\em Unitary representations of
nilpotent Lie groups}, Uspehi Mat. Nauk \textbf{17} (1962), no. 4 (106),
57--110.

\bibitem[Ki95]{kir-triangular} A.A.~Kirillov, {\em Variations on
the triangular theme}, in: ``Lie Groups and Lie Algebras:
E.B.~Dynkin's Seminar'', 43--73, Amer. Math. Soc. Transl. Ser. 2,
\textbf{169}, Providence, RI, 1995.

\bibitem[La56]{lang} S.~Lang, {\em Algebraic groups over finite fields},
Amer. J. Math. \textbf{78} (1956), 555--563.

\bibitem[La83]{langlands} R.P.~Langlands, {\em Les d\'ebuts d'une formule des
traces stable}, Publ. Math. Univ. Paris VII \textbf{13}, Paris, 1983.

\bibitem[Laz54]{lazard} M.~Lazard, {\em Sur les groupes nilpotents et anneaux
de Lie}, Ann. Sci. Ecole Norm. Sup. (3) \textbf{71} (1954), 101--190.

\bibitem[LP81]{lion-perrin} G.~Lion and P.~Perrin, {\em Extension des
repr\'esentations de groupes unipotents $p$-adiques. Calculs d'obstructions},
in: ``Noncommutative harmonic analysis and Lie groups (Marseille, 1980)'',  pp.
337--356, Lecture Notes in Math. \textbf{880}, Springer, Berlin-New York, 1981.

\bibitem[Lu03]{lusztig} G.~Lusztig, {\em Character sheaves and
generalizations}, in: ``The unity of mathematics'' (In honor of the ninetieth
birthday of I.M.~Gelfand, Editors:  P.~Etingof, V.~Retakh, I.~M.~Singer),
443--455, Progr. Math. \textbf{244}, Birkh\"auser Boston, Boston, MA, 2006,
arXiv: {\tt math.RT/0309134}

\bibitem[Lu85]{char-sheaves} G.~Lusztig, {\em Character sheaves
I-V}, Adv. in Math. \textbf{56}, \textbf{57}, \textbf{59}, \textbf{61}
(1985,1986).

\bibitem[Ma96]{matveev} A.V.~Matveev, {\em On a connection between irreducible
representations of Lie algebras and irreducible representations of $p$-groups},
Mat. Sb. \textbf{187} (1996), no. 7, 93--96.

\bibitem[NT89]{n-t} H.~Nagao and Y.~Tsushima, ``Representations of finite groups''.
Translated from the Japanese. Academic Press, Inc., Boston, MA, 1989.

\bibitem[Pre95]{previtali} A.~Previtali, {\em On a conjecture concerning
character degrees of some $p$-groups}, Arch. Math. (Basel) \textbf{65} (1995),
no.~5, 375--378.

\bibitem[Sa96]{saibi} M.~Saibi, {\em Transformation de Fourier-Deligne sur les
groupes unipotents}, Ann. Inst. Fourier (Grenoble) \textbf{46} (1996), no.~5,
1205--1242.

\bibitem[Se60]{grps-proalg} J.-P.~Serre, {\em Groupes
proalg\'ebriques}, Publ. Math. IHES \textbf{7} (1960).

\bibitem[Se77]{serre-reps} J.-P.~Serre, ``Linear representations of finite
groups,'' Springer-Verlag, New York-Heidelberg, 1977.

\bibitem[SGA4]{sga4} M.~Artin, A.~Grothendieck and J.-L.~Verdier, ``SGA 4:
Th\'eorie des Topos et Cohomologie \'Etale des Sch\'emas'', Lecture Notes in
Math. \textbf{269}, \textbf{270} and \textbf{305}, Springer-Verlag, 1972, 1973.

\bibitem[SGA$4\frac{1}{2}$]{sga4.5} P.~Deligne, with J.-F.~Boutot, L.~Illusie and
J.-L.~Verdier, ``SGA $4\frac{1}{2}$: Cohomologie \'Etale'', Lecture Notes in
Math. \textbf{569}, Springer, Heidelberg, 1977.

\bibitem[SV96]{suslin-voevodsky} A.~Suslin and V.~Voevodsky, {\em Singular homology of
abstract algebraic varieties}, Invent. Math. \textbf{123} (1996),
no.~1, 61--94.

\bibitem[Ver70]{vergne} M.~Vergne, {\em Construction de sous-alg\`ebres
subordonn\'ees \`a un \'el\'ement du dual d'une alg\`ebre de Lie r\'esoluble},
C. R. Acad. Sci. Paris S\'er. A-B \textbf{270} (1970), A173--A175.

\bibitem[Wa79]{waterhouse} W.C.~Waterhouse,
``Introduction to affine group schemes'', Graduate Texts in Mathematics
\textbf{66}, Springer-Verlag, New York-Berlin, 1979.

\bibitem[Yan01]{ning-yan} N.~Yan, {\em Representation theory of
the finite unipotent linear groups}, PhD thesis, University of
Pennsylvania, 2001 (unpublished).

\end{thebibliography}
\end{document}